\documentclass[a4paper,12pt]{book}
\usepackage{import}

\usepackage[utf8]{inputenc}
\usepackage[T1]{fontenc}

\usepackage[a4paper,top=5cm,bottom=4cm,left=3cm,right=3cm,marginparwidth=2cm, headheight=110pt]{geometry}

\usepackage{mathrsfs}

\usepackage{amsmath, amsthm, amsfonts, amscd, amssymb}

\usepackage{makeidx}
\makeindex

\usepackage{commath}

\usepackage[all,cmtip]{xy}

\usepackage{tikz-cd}
 \usetikzlibrary{positioning}

\usepackage{mathtools}

\usepackage[utf8]{inputenc}

\usepackage{pdfpages}

\usepackage{graphicx}
\usepackage[colorinlistoftodos]{todonotes}
\usepackage[colorlinks=true, allcolors=blue]{hyperref}

\usepackage[utf8]{inputenc}
\usepackage{xcolor,graphicx}
\usepackage{transparent}
\usepackage{eepic}
\definecolor{verdeuam}{cmyk}{.78,.01,.83,.02}
\definecolor{azullateral}{cmyk}{.97,.76,.01,.03}

\usepackage[utf8]{inputenc}

\usepackage{cleveref}

\usepackage{afterpage}

\newcommand{\ove}[1]{\overline{#1}}

\newcommand{\PSL}{\mathrm{PSL}}

\newcommand{\Mat}{\mathrm{Mat}}

\newcommand{\C}{\mathbb{C}}

\newcommand{\N}{\mathbb{N}}
\newcommand{\Q}{\mathbb{Q}}
\newcommand{\R}{\mathbb{R}}

\newcommand{\Z}{\mathbb{Z}}

\renewcommand{\P}{\mathbb{P}}

\DeclareMathOperator{\Aut}{Aut}

\DeclareMathOperator{\Hom}{Hom}
\DeclareMathOperator{\id}{id}

\DeclareMathOperator{\im}{im}

\DeclareMathOperator{\Gal}{Gal}

\renewcommand{\b}{b_1^{(2)}}

\DeclareMathOperator{\ab}{ab}
\newcommand{\bC}{\overline{C}}
\DeclareMathOperator{\Der}{Der}
\newcommand{\al}{\alpha}
\newcommand{\be}{\beta}
\newcommand{\ga}{\gamma}
\newcommand{\Ga}{\Gamma}
\newcommand{\La}{\Lambda}
\newcommand{\p}{\F_p}
\newcommand{\bp}{\beta_1^{ \mathtt{mod \, p}}}
\newcommand{\abr}{r_{\ab}}
\DeclareMathOperator{\Frac}{Frac}
\newcommand{\CC}{\mathcal{C}}
\DeclareMathOperator{\rk}{rk}
\newcommand{\Qo}{Q_{\text{ore}}}
\DeclareMathOperator{\ore}{ore}
\newcommand{\n}{\unlhd}
\newcommand{\no}{\unlhd_{o}}
\newcommand{\lo}{\leq_{o}}
\newcommand{\ag}{G^{\text{abs}}}
\newcommand{\af}{\phi^{\text{abs}}}
\newcommand{\ap}{\psi^{\text{abs}}}
\newcommand{\hp}{\hat{p}}
\newcommand{\cdp}{\text{cd}_p}
\newcommand{\pH}{\p[[{\bf H}]]}
\newcommand{\pF}{\p[[{\bf F}]]}

\newcommand{\sagecode}[1]{}
\newcommand{\D}{\mathcal{D}}

\newcommand{\mm}{\mathfrak{m}}
\newcommand{\F}{\mathbb{F}}
\newcommand{\FF}{{\bf F}}
\newcommand{\GG}{{\bf G}}
\newcommand{\HH}{{\bf H}}
\newcommand{\NN}{{\bf N}}
\newcommand{\KK}{{\bf K}}

\renewcommand{\subset}{\subseteq}
\newcommand{\sub}{\subseteq}

\newcommand{\tens}[1]{%
  \mathbin{\mathop{\otimes}\limits_{#1}}%
}

\newcommand{\hn}[1]{%
  \mathbin{\mathop{*}\limits_{#1}}%
}

\renewcommand{\hat}{\widehat}
\newcommand{\ti}{\widetilde}

\renewcommand{\epsilon}{\varepsilon}

\newcommand{\lrar}{\longrightarrow}

\newcommand{\rar}{\rightarrow}

\newcommand{\lan}{\langle}
\newcommand{\ran}{\rangle}
\newcommand{\inc}{\xhookrightarrow{}}

\usepackage{scalerel} 
\usepackage{graphicx} \parskip 1ex

\renewcommand\Hat[1]{\arraycolsep=0pt\relax%
\begin{array}{c} 
\stretchto{ 
\scaleto{ 
\scalerel*[\widthof{\ensuremath{#1}}]{\kern-.5pt\bigwedge\kern-.5pt} 
{\rule[-\textheight/2]{1ex}{\textheight}} 
 }{\textheight} %
 }{0.5ex}\\ 
 #1\\ 
 \rule{-1ex}{0ex} 
\end{array} 
}

\usepackage{verbatim} 
\DeclareFontFamily{U}{mathx}{\hyphenchar\font45} \DeclareFontShape{U}{mathx}{m}{n}{ <5> <6> <7> <8> <9> <10> <10.95> <12> <14.4> <17.28> <20.74> <24.88> mathx10 }{} \DeclareSymbolFont{mathx}{U}{mathx}{m}{n} \DeclareFontSubstitution{U}{mathx}{m}{n} \DeclareMathAccent{\widecheck}{0}{mathx}{"71} \DeclareMathAccent{\wideparen}{0}{mathx}{"75}

\newtheorem{thm}{Theorem}[section]
\theoremstyle{definition}
\newtheorem{defi}[thm]{Definition}
\newtheorem{lemma}[thm]{Lemma}
\newtheorem{lem}[thm]{Lemma}

\newtheorem{remark}[thm]{Remark}

\newtheorem{eg}[thm]{Example}

\newtheorem{prop}[thm]{Proposition}
\newtheorem{cor}[thm]{Corollary}
\newtheorem{notation}[thm]{Notation}
\newtheorem{nota}[thm]{Notation}

\newtheorem*{thm*}{Theorem}

\include{thebibliography}

\include{thebibliography}
\begin{document}


\thispagestyle{empty}

\begin{picture}(0,0)
\put(-50,-380)
{\scalebox{1.1}{
\begin{minipage}{1\textwidth}
\noindent{\begin{minipage}[t][.9\textheight][t]{.98\textwidth}
{\begin{minipage}[t][90pt][t]{.6\textwidth}
{\begin{minipage}[t][50pt][b]{.45\textwidth}
\centerline{\transparent{1}{\includegraphics[width=.9\textwidth]{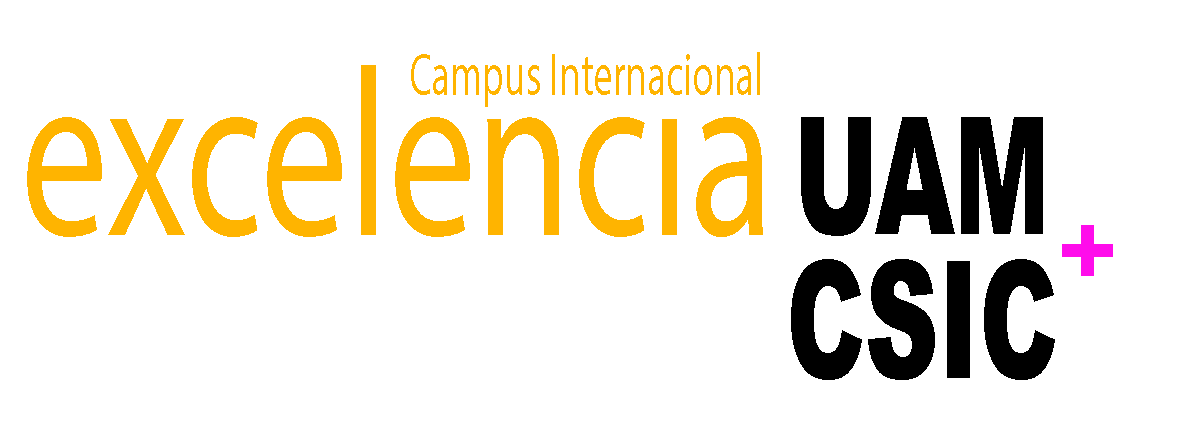}}}
\end{minipage}}

\vfill

{\begin{minipage}{1.3\textwidth}\parindent=0pt
\textcolor{verdeuam}{\rule{1\textwidth}{1\baselineskip}}

\textcolor[gray]{.8}{\rule{1\textwidth}{1\baselineskip}}

\raisebox{.9\baselineskip}{\textcolor{verdeuam}{\rule{1\textwidth}{2pt}}}

\end{minipage}}
\end{minipage}}
\colorbox{verdeuam}{\begin{minipage}[t][120pt][c]{.4\textwidth}
	\centerline{\transparent{1}{\includegraphics[width=.9\textwidth]{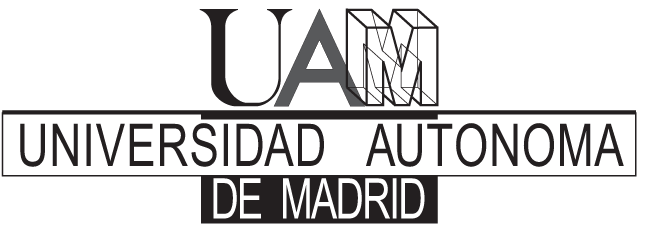}}}
\end{minipage}}

\vspace{1.3cm}

\begin{center}
            
        \Large
        \textbf{Embeddings into pro-$p$ groups and the construction of parafree groups}
            
            
        \vspace{1.7cm}
        \large  
        Ismael Morales López\\
        \vspace{0.7cm}
        \normalsize 
        Supervised by Andrei Jaikin-Zapirain\\
        \vspace{0.4cm}
        \normalsize 
        Universidad Autónoma de Madrid
            
            
        \vspace{1.8cm}
            
            
            
    \end{center}

\begin{minipage}{.36\textwidth}
\transparent{1}{\includegraphics[width=.9\textwidth]{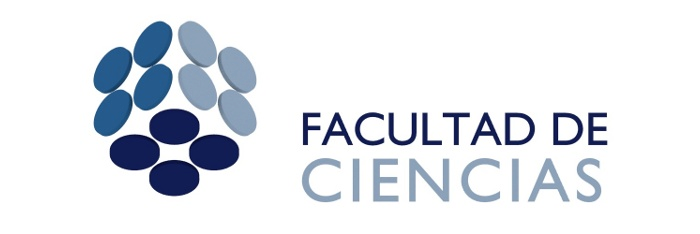}}
\end{minipage}\hfill
\begin{minipage}{.7\textwidth}
{\small This master thesis was supported by the Spanish Ministry of  }\\
{\small Science and Innovation through the grant JAE Intro SOMdM of the academic year 2020-2021.} 

\end{minipage}
\end{minipage}}%

\end{minipage}}
}
\end{picture}

\clearpage

\begin{center}
    {\Large \bf Abstract}
\end{center}
The main objects of study in this thesis are parafree groups, although the theory of pro-$p$ groups plays a significant role. A finitely generated group $G$ is termed {\bf parafree} if
\begin{itemize} 
\item[(i)]  there exists a free group $F$ which has the same isomorphism types of nilpotent quotients as  $G$; and
\item[(ii)] $G$ is residually nilpotent.
\end{itemize}




There are two original results in this thesis. On the one hand, we can completely describe which free products of finitely generated groups, with abelian amalgams, are parafree.  Secondly, we can  characterise more precisely  which abelian HNN extensions of finitely generated groups are parafree.

The latter is not an explicit description, although it can still describe, as a corollary, which abelian HNN extensions of finitely generated groups  $G$ with $d\left(G/[G, G]\right)=2$ are parafree.

\clearpage


\thispagestyle{empty}
\begin{center}
    {\Large \bf Acknowledgements}
\end{center}

The author wishes to thank his supervisor, Andrei Jaikin-Zapirain, for his generosity in sharing many exciting and original ideas during these last years. 

\clearpage

\tableofcontents

\mainmatter
\chapter{Introduction} \label{introduction}

In this thesis we study groups that resemble free groups in terms of their residual, cohomological and structural properties. 

There has been much research in group theory around comparing infinite groups by the study of some of their quotients. A famous open problem is the following. 

{\bf Remeslennikov's Question:} Given a finitely generated free group $F$ and a finitely generated residually finite group $G$ having the same isomorphism types of finite quotients, does it follow that $F\cong G$?

One can pose this question in more general terms. A {\bf variety of groups $\CC$}  is a collection of groups closed under taking subgroups, quotients and direct products. A group $G$ is said to be {\bf residually-$C$} if for every $1\neq g\in G$ there exists a normal $N\unlhd G$ such that $g\notin N$ and $G/N\in C$. This way, one can ask whether every two finitely generated residually-$\CC$ groups having the same isomorphism types of quotients belonging to $\CC$ are necessarily isomorphic.   

The problem of Remeslennikov  is conjectured to have  positive answer. However, for other varieties of groups, such as the variety of nilpotent groups or finite $p$-groups; the analogous question is known to have negative answer. 

We say that two groups $G$ and $H$ have the same {\bf nilpotent genus} if they have the same isomorphism types of nilpotent quotients. Even if it does not follow that $G\cong H$, one can still ask which properties must $G$ and $H$ share. 

Bridson and Reid \cite{Bri14} solved a few problems posed by Baumslag regarding the previous question; and these results exhibit the possible divergence between groups $G$ and $H$ having the same nilpotent genus. For example, they show that there exist two pairs of finitely generated residually-(torsion-free nilpotent) groups $H\inc G$ with the same nilpotent genus such that: for the first pair, $G$ has solvable conjugacy problem while $N$ does not; and, for the second pair, $G$ is finitely presented while $H_2(H, \Q)$ is infinite-dimensional.

The case we study in this thesis corresponds to one of them being free. 
A finitely generated group $G$ is said to be {\bf parafree} if it is residually nilpotent and has the same isomorphism types of nilpotent quotients as some free group.  Equivalently, a finitely generated residually nilpotent group $G$ is parafree if there exists a free group $F$ such that, for all $k$, $G/\ga_k G\cong F/\ga_k F$.

Baumslag introduced this class of groups. He gave many examples of non-free parafree groups \cite{Bau67} and he also studied  some of their key structural properties \cite{Bau69}. As opposed to general pairs of finitely generated residually-(torsion-free nilpotent) groups of the same nilpotent genus, parafree groups are expected to share many properties of free groups. We can mention a few. The Parafree conjecture, which appears in Kirby's list \cite[Problem 3.78]{Kir95}, expects $H_2(G, \Z)$ to vanish if $G$ is parafree. It is also conjectured that parafree groups have cohomological dimension at most two. However, these problems are particularly difficult keeping in mind we do not even know if (finitely generated) parafree groups are finitely presented. 

In this thesis we do not address these questions. We are concerned with the construction of parafree groups.  Most known examples of parafree groups essentially fall into one of two classes of groups, namely
\begin{itemize}
    \item[(a)] amalgamated products of a parafree group and $\Z$; and
    \item[(b)] cyclic HNN extensions of free groups.
\end{itemize}
Our purpose is to describe  which abelian HNN extensions of parafree groups and which free products of parafree groups with abelian amalgams are still parafree. We should remark that two-generated subgroups of parafree groups are free, so  non-trivial abelian subgroups of parafree groups are isomorphic to $\Z$.

Given a group $G$, we denote by $\abr (G)$ the minimal number of generators of its abelianisation $G_{\ab}=G/[G, G]$. We name $\abr(G)$ the {\bf abelian rank} of $G$.

There are two original results in this thesis. In regard to amalgamated products of parafree groups, we prove the following (\cref{amalgamparafree}). 

\begin{thm*}[{\bf Parafree amalgamated products of cyclic amalgam}]  Let $U$ and $V$ be finitely generated groups. Let $1\neq u\in U$ and let $1\neq v\in V$. Consider the following amalgamated product of cyclic amalgam
\[W=U\underset{u=v}{*}V\cong \frac{U* V}{\lan \lan uv^{-1}\ran \ran}.\]
Then $W$ is parafree if and only if the three following conditions hold:
\begin{enumerate}
    \item The groups $U$ and $V$ are parafree.
    \item The element $uv^{-1}$ of $U* V$ is not a proper power in the abelianisation.
    \item At least one of $u$ or $v$ is not a proper power in $U$ or $V$, respectively.
\end{enumerate} 
In this case, $W$ is parafree of abelian rank $\abr(U)+\abr(V)-1$. 
\end{thm*}

In regard to HNN extensions of parafree groups, we prove the following (\cref{hnnparafree}). 
\begin{thm*}[{\bf Parafree HNN cyclic extensions}]  Let $U$ be a finitely generated group. Let $u, v\in U\setminus \{1\}$. Consider the following cyclic HNN extension of $U$ \[W=\frac{U* \lan t\ran}{\lan \lan tut^{-1}v^{-1}\ran \ran}. \]
Then $W$ is  parafree if and only if the four following conditions hold:
\begin{enumerate}
    \item The group $U$ is parafree.
    \item The image of $uv^{-1}$ is not a proper power in $U_{\ab}$.
    \item At least one of $u$ or $v$ is not a proper power in $U$.
    \item The image of $u$ is non-trivial in some finite nilpotent quotient of $U$.
\end{enumerate}
In this case,  $W$ is  parafree of the same abelian rank as $W$. 
\end{thm*}

The fourth property of the previous theorem may be, in general, difficult to check. However, for parafree groups of abelian rank 2, this can be reformulated in very simple terms (\cref{hnnparafree2}).

\begin{thm*}[{\bf Parafree cyclic HNN extensions of groups with abelian rank 2}]
Let $U$ be a finitely generated group of abelian rank 2. Let $u, v\in U\setminus \{1\}$. Consider the following cyclic HNN extension of $U$ \[W=\frac{U* \lan t\ran}{\lan \lan tut^{-1}v^{-1}\ran \ran}. \]
Then $W$ is  parafree if and only if the four following conditions hold:
\begin{enumerate}
    \item The group $U$ is parafree.
    \item The image of $uv^{-1}$ is not a proper power in $U_{\ab}$.
    \item At least one of $u$ or $v$ is not a proper power in $U$.
    \item  The images of $u$ and $v$ generate a subgroup isomorphic to $\Z^2$ in $U_{\ab}$. 
\end{enumerate}
In this case,  $W$ is  parafree of abelian rank 2. 
\end{thm*}

The methods that we use to prove these results mainly rely on the theory of group rings of residually-(torsion-free nilpotent) groups and a dimension theory on their modules, for whose required computations and estimations we shall use the theory of $L^2$-Betti numbers. A combination of these ideas lead to a criterion for constructing embeddings of abstract groups into free pro-$p$ groups (\cref{pro-$p$ embedding lemma}. This criterion was developed by Jaikin-Zapirain \cite{And20} to embed finitely generated subgroups of free $\Q$-groups into free pro-$p$ groups. The aforementioned theorems for characterising parafree amalgamated products and HNN extensions are studied in more depth and extended to the more general context of fundamental groups in \cite{And21}, where we prove the following. 
\begin{thm*}[Corollary 1.4, \cite{And21}]  Let $(\mathcal G, \Gamma)$ be a  graph of   groups over a finite 
graph $\Gamma$ whose edge morphisms are injective. Let  $W=\pi(\mathcal G, \Gamma)$ be its  fundamental group. Assume that all  vertex subgroups  $\mathcal G(v)$ ($v\in V(\Gamma)$)   are finitely generated and all   edge subgroups $\mathcal G(e)$ ($e\in E(\Gamma)$)  are cyclic. Then $W$ is parafree if and only if the following four conditions hold.
 
 \begin{enumerate}
 \item All the vertex  subgroups $\mathcal G(v)$ ($v\in V(\Gamma)$)    are parafree.
 
 \item The abelianisation of $W$ is torsion-free   of rank 
 $$\abr(W)=\sum_{v\in V(\Gamma)}\abr(\mathcal G(v))-\sum_{e\in E(\Gamma)} \abr(\mathcal G(e)) -\chi(\Gamma),$$   where $\chi(\Gamma)=|V(\Gamma)|-|E(\Gamma)|$-1.

 \item All the centralisers of non-trivial elements in $W$ are cyclic.

 \item For  each non-trivial edge subgroup of  $\mathcal G(e)$  ($e\in E(\Gamma)$) there is a finite nilpotent quotient of $W$ where the image of this edge subgroup is non-trivial.
 \end{enumerate}

\end{thm*}

We would like to conclude this introduction by pointing out a recent development \cite{And212} which shows in particular that in order to give a positive solution to Remeslennikov's problem, it suffices to rule out non-free parafree groups. 

\clearpage
\section{Structure of the exposition}

\begin{itemize}
    \item {\bf Chapter 1: } We briefly introduce the problem of constructing parafree groups and we overview our two main results of the thesis. In \cref{notationconven}, we fix some conventions that will be followed during the whole text.
    \item {\bf Chapter 2:} We discuss some important classes of groups in terms of their structural and residual properties in order to compare them, later on, with parafree groups.
    \item {\bf Chapter 3:} Here we review two old applications of Bass-Serre theory on the structure of elements and subgroups of amalgamated products and HNN extensions.
    \item {\bf Chapter 4:} We introduce profinite groups and we give some tools to study pro-$p$ groups. We precise the notion of free product in the category of pro-$p$ groups to give more explicit descriptions of pro-$p$ completions of amalgamated products and HNN extensions in terms of the pro-$p$ completions of the factors during \cref{presentationsection}. 
    \item {\bf Chapter 5:} Here we collect some known examples of parafree groups. We give a pro-$p$ reformulation of the property of having the same nilpotent genus as a free group in \cref{characterisations}. In our applications, we work with this viewpoint. 
    \item {\bf Chapter 6:} Here we recall some standard properties about the augmentation ideals, with the aim of applying the main result of the chapter (\cref{kernelmodp}) to groups that arise as the construction of an amalgamated product or an HNN extension. 
    \item {\bf Chapter 7:} A more sophisticated and concrete version of the previous \cref{kernelmodp} requires developing further the theory of certain groups rings and a dimension theory on their modules. The main results are \cref{D12},  \cref{amaltor} and \cref{hnntor}.
    \item {\bf Chapter 8:} We introduce some techniques of $L^2$-Betti numbers that will be required to estimate dimensions of modules related to augmentation ideals. The main result is \cref{parabetti}.
    \item {\bf Chapter 9:} Here we start introducing the main tool for constructing embeddings of abstract groups into pro-$p$ groups (\cref{pro-$p$ embedding lemma}), which is a sophisticated version of \cref{kernelmodp} mixed with the methods of Chapters 7 and 8. This will be applied to the settings of amalgamated products and HNN extensions of parafree groups. The main results are theorems \ref{amalgamparafree} and \ref{hnnparafree}.
\end{itemize}
\section{Notation and conventions} \label{notationconven}
We remind a few basic definitions and we set several conventions  which are worth-recalling because some of them are  not standard.
\begin{enumerate}
    \item All of our rings $R$ will be associative and unitary. All ring homomorphisms will map $1\mapsto 1$. We denote by $R^{\times}$ or $R^*$ the multiplicative group of units of $R$. All of our modules will be left modules if there is no explicit mention. 
    
    \item Given a group $G$, we will denote by $G_{\ab}=G/[G, G]$ its abelianisation and by $G'=[G, G]$ its commutator subgroup (also known as the derived subgroup). The  $p$-abelianisation of $G$  is the quotient $G/G^{p}[G, G]$, which is isomorphic to both $H_1(G, \F_p)$ and $H^1(G, \F_p)$.
    \item Given a word $\omega\in F_n$, where $F_n$ is the free group considered with free generators $x_1, \dots, x_n$; we denote $\omega_{x_i}\in \Z$ to be defined in such a way that the canonical image of $\omega$ in the abelianisation $\Z^n$ of $F_n$ is $(\omega_{x_1}, \dots, \omega_{x_n})\in \Z^n$. 
    \item We denote by $k$ a commutative ring. We are particularly interested in $k=\F_p$ and $k=\Z$. 
    \item Given a set $S$, we denote by $kG^{(S)}$ the free $kG$ module with a basis $\{e_s\}_{s\in S}$ indexed by $S$.
    \item 
Given a group $G$, we denote by $G_k$ or $\ga_k G$ (resp. $G_{k, p}$ or $\ga_{k, p} G$) its {\bf lower central series} (resp. its {\bf $p$-lower central series}). These are defined recursively  as follows: $G_1=G_{1, p}=G$ and, for $n\geq 1$, 
\[G_{n+1}=[G_{n}, G], \, \, \, G_{n+1, p}=G_n^p[G_n, G].\]
If $G=F$ is free, we might denote $G_k$ by $G_{(k)}$ in order to avoid confusion, because we leave the notation $F_k$ for the free group on $k$ generators. 
We refer to the quotients $G/G_k$ (resp, $G/G_{k, p}$) as the lower central quotients (resp. $p$-lower central quotients) of $G$.
\item 
Whenever we say that a map is \textit{natural} we mean that it is functorial. This is of particular importance in homological arguments. On the other side, we say that a map, or other mathematical object, is \textit{canonical} if it does not rely on a particular choice of generators or data. Sometimes, a map can be both canonical and natural. For example, given $\R$-vector spaces $U$ and $V$, $U\tens{\R} V$ is canonically isomorphic to $V\tens{\R} U$; although this isomorphism is also natural.

\item 
The group $\Z_p$ is the group of $p$-adic integers. The letters  $\FF, \GG, \HH, \NN$ will denote pro-$p$ groups. We reserve $\FF$ to denote a free pro-$p$ group. The free group of rank $n$ will be denoted by $\FF_n$.
  \item Let $G$ be a group. We denote its {\bf rank} by
\[d(G)=\inf \{m: \mbox{there exists a generating set of $G$ with $m$ elements}\} \in \Z_{\geq 0}\cup \{\infty\}.\]
We denote its {\bf abelian rank} (usually named \textit{para-rank} or \textit{parafree rank}) by
\[\abr(G)=d(G_{\ab}).\]
If $\GG$ is a profinite group, we also denote 
\[d(\GG)=\inf \{m: \mbox{$\GG$ admits a topological generating set with $m$ elements}\}.\]
We say that $\GG$ is topologically finitely generated, or simply that it has finite rank, if $d(\GG)<\infty$.

\item Given a pro-$p$ group $\GG$ of finite rank, we denote by $\ga_k^{(p)}\GG$ its {\bf $p$-lower central series}.  These are defined recursively by $\gamma_{1}^{(p)}\GG=\GG$ and  $\gamma_{n+1}^{(p)}\GG=\GG^p\, [\GG, \gamma_{n}^{(p)}\GG]$ for $n\geq 1$. These are open normal subgroups of $\GG$. Similarly, we denote by $\ga_k\GG$  the {\bf lower central series} of $\GG$, which are defined recursively by $\ga_1\GG=\GG$ and $\ga_{k+1}\GG=[\ga_k\GG, \GG]$ for $k\geq 1$. These are closed normal subgroups of $\GG$. 

\item Given a group $G$ of the form $\F_p^n$,  $\Z^n$ or $F_n$; and a subset $S\sub G$; we say that $S$ is {\bf primitive} if it is part of a generating set of $n$ elements. If $G=\FF_n$ then we say $S\sub G$ is primitive if it is part of a topological generating set of $n$ elements.

\item Given a group $G$, we denote by $\Phi(G)$ its {\bf Frattini subgroup}.

\item Given an abstract group $G$, we denote its {\bf pro-$p$ completion} by $G_{\hp}$. 

\item Given a pro-$p$ group $\GG$, we will write $\HH\lo \GG$ if $\HH$ is an open subgroup, and $\NN\no \GG$ if $\NN$ is a normal open subgroup. 
\item A {\bf variety of groups} will always be denoted by $\CC$ and $\CC_p$ denotes the variety of finite $p$-groups. 
\end{enumerate}

\chapter{Some special classes of groups} \label{special}

In this chapter, we simply review a few characteristic properties of three classes of groups; namely free groups, surface groups, nilpotent and residually-(torsion-free nilpotent) groups. 

Free groups are our starting point. Their subgroup structure is particularly well-understood; and so are their residual and cohomological properties. 

\begin{thm}[Schreier's index-rank formula] \label{freeindex} Let $H$ be a subgroup of $F_n$ of finite index $k$. Then $H$ is a free group of rank equal to $k(n-1)+1$.
\end{thm}

This fundamental property can be proven by combinatorial methods (\cite[Proposition 3.9]{Lyn01}) and also using more sophisticated tools: by covering space theory (\cite[Chapter 4, Section 5]{May99}) and by Bass-Serre theory (\cite[Theorem 1.2, Chapter II]{Dic80}.
We will see that surface groups also enjoy a similar property. However, the feature of free groups we are most interested in is the following. 

\begin{thm}[Magnus] \label{freeres1} Free groups are residually-(torsion-free nilpotent). 
\end{thm}
\begin{proof}[Proof] Since free groups are residually-(finitely-generated free groups), it suffices to prove the statement for finitely generated free groups $F_n$. On the one hand, a result of Magnus states that $F_n/\ga_m F_n$ is torsion-free for all $m$. We can prove this inductively. The induction base $m=1$ is trivial and $F_n/\ga_{m+1} F_n$ fits into a short exact sequence 
\[1\lrar \frac{\ga_m F_n}{\ga_{m+1} F_n} \lrar \frac{F_n}{\ga_{m+1} F_n} \lrar \frac{F_n}{\ga_m F_n}\lrar 1. \]
The first group is a torsion-free abelian group freely generated by basic commutators. This is a consequence of the normal form of elements of $F_n/\ga_m F_n$ (described, using the collecting commutator process,  in \cite[Chapter 11]{Hal76}, for example). By the inductive hypothesis, $F_n/\ga_{m} F_n$ is also torsion-free, so $F_n/\ga_{m+1} F_n$ is torsion-free. 

It rests to verify that $F_n$ is residually-nilpotent. One argument consists on considering an embedding $F_n\inc \Q\lan \lan X_1, \dots, X_n\ran \ran^{\times}$, where the last group stands for the unit group of the ring $\Q\lan \lan X_1, \dots, X_n\ran \ran$ of power series in non-commuting indeterminates $X_1,\dots, X_n$, with coefficients in $\Q$. This embedding is defined by taking free generators $x_i$ to $1+X_i$, respectively. One would then notice that non-trivial elements of $\ga_m F_n$ embed into elements of the form $1+\sum_{i_1, \dots, i_k; k\geq m} a_{i_1, \dots, i_k} \, X_{i_1}\cdots X_{i_k}$; or, in other words, of the form $1+$ (power series supported in monomials of degree  $\geq m$). This is an standard argument originally due to Malcev and the last claim can be easily verified inductively (see \cite[Section 6]{Sta65} for details).

An alternative way to check that $F_n$ is residually nilpotent is using the fact that $F_n=\ga_1 F_n\supseteq \ga_2 F_n\supseteq \dots \supseteq \ga_{m+1} F_n\supseteq \dots $ is an strictly decreasing chain that verifies that $\ga_{m+1} F_n$ is characteristic in $\ga_m F_n$ for all $m$. By \cite[Theorem 5]{Tak51}, its intersection $\bigcap_m \ga_m F_n$ is trivial. 
\end{proof}

By \cref{nilresp} and the previous theorem, free groups are also residually-$p$ for every prime $p$. The two  ways to establish the residual nilpotence of free groups that we discussed during the proof of \cref{freeres1} are purely algebraic. One can use covering space theory methods to directly  prove that free groups are residually-$p$ for every prime $p$ (see \cite[Section 3.2]{Wil21}).

A deep defining theorem of free groups  due to John R. Stallings and   Richard Swan states that free groups are exactly the abstract groups of cohomological dimension 1. During \cref{propcohomology}, we will establish the analogous principle for pro-$p$ groups, which is significantly simpler.

\section{Surface groups}
Fundamental groups of connected and compact surfaces are called surface groups. From the classification of these surfaces, one can write down a classification of such groups. Furthermore, from their topological description as fundamental polygons, the classification of surface groups is precised in terms of their presentations. Non-closed surfaces (those with non-empty boundary) have free fundamental group and closed surfaces can be classified as follows. 

\begin{thm}\label{surfaceclass} The fundamental group of the closed orientable surface of genus $g\geq 0$, denoted by $\Sigma_g$, is 
\[\pi_1(\Sigma_g)\cong\lan x_1, \dots, x_g, y_1, \dots, y_g : [x_1, y_1]\cdots [x_g, y_g]=1\ran .\]
The fundamental group of the closed non-orientable surface of genus $g\geq 1$, denoted by $S_g$, is 
\[\pi_1(S_g)\cong\lan x_0, x_1, \dots, x_g : x_0^2x_1^2\cdots x_g^2=1\ran .\]
\end{thm}

\begin{remark} \label{surfaceab} Directly from their presentations, one can compute the abelianisation of a surface group. For example, let $G=\pi_1(\Sigma_g)$, then 
\[G/[G, G]\cong \Z^{2g}.\]
If $G=\pi_1(S_g)$, then 
\[G/[G, G]\cong \Z^{g}\times \Z/2.\]
In particular, their first Betti numbers (see \cref{Betti}) are $b_1(\pi_1(\Sigma_g))=2g$ and $b_1(\pi_1(S_g))=g$. 
\end{remark}

Interestingly, it is still fruitful to consider them as fundamental groups of a geometric object because this way one can understand their subgroups using covering space theory, as one does for free groups.

\begin{thm} \label{surfaceindex} Let $S$ be a closed surface and let $G$ be its fundamental group. Let $H$ be a subgroup of $G$. There are two cases:
\begin{itemize}
    \item If $H$ is has finite index in $G$ then $H\cong \pi_1(\hat{S})$, where $\hat{S}$ is a closed surface of Euler Characteristic $\chi (\hat{S})=|G: H|\, \chi (S)$. 
    \item If $H$ has infinite index in $G$, then $H$ is free. 
\end{itemize}
\end{thm}

Using covering space methods one can also verify residual properties of surface groups (see, for example, \cite{Hem72}).

\begin{prop} \label{surfaceresfinite} Surface groups are residually finite. 
\end{prop}
The previous result admits another proof in the orientable case since Fricke and Klein proved that $\pi_1(\Sigma_g)$ has a faithful representations in $\PSL_2(\C)$ and, by a result of Mal'cev, we know that $\PSL_2(\C)$ is residually finite. 

 Baumslag established only by purely algebraic means the following much stronger property in \cite{Bau62}. 

\begin{prop} Surface groups are residually free. 
\end{prop}

\section{Nilpotent groups}
Since parafree groups are residually nilpotent groups, it would be convenient to collect here some properties of this class of groups. 

For a thorough introductory development of the theory of nilpotent groups, we recommend the book of M. Hall \cite[Chapter 10]{Hal76}, which contains proofs for most of the results we are about to mention. The proofs of the residual properties of these groups use an standard inductive argument on the Hirsch length.  


\begin{prop} The following is true. \label{Hopfnilp}
\begin{enumerate}
    \item Finitely generated nilpotent groups are polycyclic.
    \item Polycyclic groups are residually finite. 
    \item Finitely generated residually finite groups are Hopfian. 
\end{enumerate}  
\end{prop}

\begin{remark} \label{Hopf1}
The Hopf property can be used to check whether a surjective map is injective. For example, if $G$ is a Hopfian group; and $f: G\rar H$ and $g: H\rar G$ are surjective group homomorphisms; then both $f$ and $g$ are isomorphisms. This is due to the fact that $g\circ f$ is a surjective endomorphism of $G$. By the Hopf property, $g\circ f$ must be injective. Since $f$ is surjective, the previous implies that $f$ and $g$ are injective. 
\end{remark}
We should still remark that finitely generated nilpotent groups are a much more particular class of groups than the class of finitely generated residually finite groups. In fact, they are closer to finitely generated abelian groups. In this sense, a different explanation for the Hopf property of finitely generated nilpotent groups $G$ can be given by means of the \textit{Noetherian condition}. R. Baer proved that any subgroup of $G$ is finitely generated \cite[Theorem 2.18]{Cle17}. From this, it directly follows that $G$ is Hopf.

The following is a consequence of the Burnside Basis Theorem for nilpotent groups. 
\begin{prop} \label{Basisnilp} Let $G$ be a nilpotent group. Then $[G, G]\subset \Phi(G)$. In particular, if a subgroup $H\leq G$ verifies that $H\, [G, G]=G$, then $H=G$.
\end{prop}

Parafree groups are not only residually nilpotent but also residually-(torsion-free nilpotent). This is due to the fact that, when $F$ is free, the quotients $F/\ga_n F$ are torsion-free nilpotent. 

In addition, this proves that parafree groups are residually-$p$ for every prime $p$.

\begin{prop}[Theorem 2.1, \cite{Gru57}] \label{nilresp} Finitely generated torsion-free nilpotent groups are residually-$p$ for every prime $p$. 
\end{prop}

The class of residually-(torsion-free nilpotent), which includes parafree groups, locally resemblances the family of groups  $\Z^n$ and, as such,  they fall into other interesting classes of groups. They are both orderable and locally indicable.

\begin{defi} A group $G$ is {\bf orderable} if there exists an order relation $\leq$ on $G$ with the property that, for all $r, s, t$, if $r\leq s$, then $rt\leq st$ and $tr\leq ts$. A  subset $S\subset G$ is said to be well-ordered if every subset of $S$ has a minimum. We denote this $\leq $ a {\bf group order} on $G$.
\end{defi}

There is a simpler way to codify a group ordering, with the concept of positive cone.

\begin{lem} \label{cone} A group $G$ is orderable if and only if there exists a subset $P\sub G$ such that $\{P, P^{-1}\}$ is a partition of $G\setminus \{e\}$, $P\cdot P\sub P$ and $gPg^{-1}\sub G$ for every $g\in G$. 
\end{lem}
\begin{proof}
Given an order $\leq$ on $G$, we can take the positive cone $P=\{g\in G: 1<g\}$. Reciprocally, given such subset $P$, we can define $x<y$ if $x^{-1}y\in P$. We write $x\leq y$ if $x=y$ or $x<y$. This defines an order on $G$.
\end{proof}

\begin{lem} \label{ordext} Suppose that a group $G$ has a central subgroup $N\leq G$ such that both $N$ and $G/N$ are orderable. Then $G$ is orderable. 
\end{lem}
\begin{proof}We use \cref{cone}. Consider the canonical projection $p: G\lrar G/N$. Given the positive cone $P_N$ for an ordering of $H$ and the positive cone $P_N'$ for an order on $G/H$, we can simply take $P=P_N\cup p^{-1}(P_N')$ as positive cone to define an order on $G$.
\end{proof}

We can now prove that finitely generated torsion-free residually nilpotent groups, and groups that are locally or residually of this way, are orderable. 

\begin{prop} \label{nilorder} Torsion-free  nilpotent groups are orderable. 
\end{prop}
\begin{proof} Let $G$ be torsion-free nilpotent. 
By an standard compactness argument, it suffices to check that $G$ is locally orderable. We will now prove the statement for finitely generated torsion-free nilpotent groups by induction on their Hirsch length. The only such groups of Hirsch length at most 1 are the trivial group and $\Z$, both orderable. Now let $H\leq G$ be a finitely generated subgroup. There is a finite normal series
\[1=G_0\unlhd G_1\unlhd \cdots \unlhd G_m=G, \]
where $G_k/G_{k-1}\cong \Z$  and $G_k/G_{k-1}\sub Z(G/G_{k-1})$ for all $1\leq k\leq m$. This $m$ is the Hirsch length of $H$ and this group is the central extension of $G_0\cong Z$ by $G/G_0$, both orderable by the induction hipothesis. So $H$ is orderable by \cref{ordext}.
\end{proof}

\begin{prop}\label{nilorder2} Residually-(torsion-free nilpotent) groups are orderable. 
\end{prop}
\begin{proof} 
Consider 
\[W=\{N\unlhd G: \mbox{$G/N$ is torsion-free nilpotent}\}.\]
We endow $W$ of a well order. For each $N\in W$, we endow $G/N$ of a group ordering. The group $\prod_{N\in W} G/N$ has a natural group order with respect to the lexicographic order on $W$ and the coordinate-wise order on each group $G/N$. The latter order exists by \cref{nilorder}. The restriction of this order to its subgroup \[G\inc \prod_{N\in W} G/N \]
proves that $G$ is orderable. 
\end{proof}

A simple observation is that these groups, namely residually-(torsion-free nilpotent) groups, are locally indicable\footnotemark.

\begin{defi} \label{locind}
We say that a group $\Ga$ is {\bf locally indicable} if every finitely generated nontrivial subgroup has an epimorphism onto $\Z$. 
\end{defi}

\footnotetext{In fact, orderable groups are locally indicable. However, we prefer not to mention this because local indicability follows directly from the structure theory of finitely generated nilpotent groups, as depicted in the proof of \cref{nilorder}.}
 
\chapter{Bass-Serre theory} \label{BassSerresection} 

We introduce one of the most fundamental and elementary tools of the Bass-Serre theory, namely the characterisation of free groups in terms of their action on graphs.

\begin{thm}[Reidemeister] \label{tree}
A group $G$ is free if and only if it acts freely on a tree. 
\end{thm}
\begin{proof} Omitted. This is a direct consequence of the structure theorem of groups acting on graphs. The reader is referred to \cite[Theorem 1.1, Chapter II]{Dic80}.
\end{proof}

A direct consequence is the following. 

\begin{cor}  A subgroup of a free group is free.
\end{cor}

A fundamental object in Bass-Serre theory is the notion of graph of groups $(\mathcal{G}, Y)$; its associated tree; and its associated fundamental group $\pi_1(\mathcal{G}, Y)$, 
which is a type of ``free construction''. With these tools, one can prove that many subgroups $H\leq \pi_1(\mathcal{G}, Y)$; such as those that act freely on the corresponding tree; are free. However, we will not need to work in such generality. We will simply revise some properties of two concrete examples, namely amalgamated products and HNN extensions. The families of free groups we previously referred to are collected in corollaries \ref{amalfree} and \ref{hnntree}. 

\subsubsection{Amalgamated products}

\begin{defi} Let $\theta: A\lrar G$ and $\omega: A\lrar H$ be two group monomorphisms. We define its corresponding amalgamated product, denoted $G\underset{A}{*} H$, by  the group 
\[G\underset{A}{*} H=  G*H\,  \Big/\, \lan \lan \theta(a)\, \omega(a)^{-1}, a\in A\ran\ran. \]
\end{defi}

We can give an explicit description of the elements of the group $G\underset{A}{*} H$.

\begin{prop}[Britton's lemma. Normal form of elements in amalgamated products] \label{brittonam} Let $\theta: A\lrar G$ and $\omega: A\lrar H$ be two group monomorphisms. Let $G_1$ (resp. $H_1$) be a set right-coset representatives of $\theta(A)$ (resp. $\omega(A)$) such that $1\in G_1$ (and $1\in H_1$). Then any element $g\in  G\underset{A}{*} H$ can be uniquely written in the form 
\[a \, g_1\, \cdots \, g_n,\]
where $a\in A$, and $1\neq g_i\in G_1\cup H_1 $ are such that they alternate from $H_1$ to $G_1$. In other words, $g_i\in H_1$ implies that $g_{i+1}\in G_1$, and vice versa. 
\end{prop}

A particular implication of this lemma is that both $G$ and $H$ canonically embed into $ G\underset{A}{*} H.$

\begin{cor} \label{brittonam1} Let $\theta: A\lrar G$ and $\omega: A\lrar H$ be two group monomorphisms. Let $n\geq 1$ and consider an element $g\in  G\underset{A}{*} H$ of the form \[g=g_0\,  g_1\, \cdots g_n,\]
where  $g_i\in G\cup H$, $g_k\notin H$ if $k>0$; and the $g_i$ alternate from $G$ to $H$. Then $g\neq 1$. 
\end{cor}

\begin{cor} \label{amalfree}
Let $F$ be a subgroup of $G \underset{A}{*} H$ that intersects trivially any conjugate of $G$ or $H$. Then $F$ is free. 
\end{cor}

\subsubsection{HNN extensions}
Here we introduce an important construction in group theory, named \textit{HNN extensions}. It was introduced by  Graham Higman, Bernhard Neumann, and Hanna Neumann. As the amalgamated products, they naturaly arise when taking the fundamental group of a certain topological construction. We will see that they also have a normal form theorem.

\begin{defi} Let $\theta: A\rar G$ be a group monomorphism. We define the {\bf HNN extension} of $G$ over $A$ with respect to $\theta$ as the group 
\[\lan G, t: tat^{-1}=\theta(a),\,  a\in A\ran \cong G*\lan t\ran \big/ \lan \lan tat^{-1}\theta(a)^{-1},\,  a\in A \ran\ran.\]
When the monomorphism $\theta$ is understood from the context, this is simply denoted by $G\hn{A}$ and the letter $t$ is called \textit{stable letter}.
\end{defi}

\begin{prop}[Britton's Lemma. Normal form of elements in HNN extensions] \label{brittonhnn} Let $\theta: A\rar G$ be a group monomorphism and let $A_{1}, A_{-1}$ be,  right-coset representatives in $G$ of the subgroups $A$ and $\theta(A)$, respectively. Any element $g\in G\, \hn{A}$ can be uniquely represented as a product
\[g_0\, t^{\epsilon_1}\, g_1\, \cdots \, t^{\epsilon_n}\, g_n,\]
where $g_0\in G$; $\epsilon_i\in \{-1, 1\}$; $g_i\in A_{1}\setminus \{1\}$ if $\epsilon_i=1$; $g_i\in A_{-1}\setminus \{1\}$ if $\epsilon_i=-1$; and $\epsilon_i=\epsilon_{i+1}$ if $g_i=1$. 
\end{prop}
In particular, $G$ embeds into $G\hn{A}$.

\begin{cor} \label{brittonhnn1} Let $\theta: A\rar G$ be a group monomorphism. Let $g\in G\, \hn{A}$ be an element of the form 
\[g_0\, t^{\epsilon_1}\, g_1\, \cdots \, t^{\epsilon_n}\, g_n,\]
where $g_0\in G$, $\epsilon_i\in \{-1, 1\}$; $g_i\notin A$ if $\epsilon_i=-\epsilon_{i+1}=1$; and $g_i\notin \theta(A)$ if $\epsilon_i=-\epsilon_{i+1}=-1$. Then $g\neq 1.$
\end{cor}

Some classical and important examples of HNN extensions are the Baumslag-Solitar groups $B(n, m)=\lan x, y| yx^n y^{-1}=x^m\ran$. 

\begin{eg} The group $B(1, m)\cong \lan x, y \, | \, yxy^{-1}=x^m\ran$ is isomorphic to $\Z[\frac{1}{m}]\rtimes \Z$, with $\phi: \Z\rar \Aut \Z[\frac{1}{m}] $ given by $\phi(n)(a)=nma$ for every $n\in \Z$, $a\in \Z[\frac{1}{m}]$. They are isomorphic via the isomorphism that maps $x\mapsto (1, 0)$ and $y\mapsto (0, 1)$. 
\end{eg}

\begin{eg} \label{BShopf} The groups $B(n, m)$, for $n, m\geq 2$ coprime, are more interesting. One can prove that in these groups are not Hopf. In fact, the endomorphism given by $x\mapsto x^n$ and $y\mapsto y$ is surjective, since it contains $x^n, x^m=yx^ny^{-1}$ and $y$, though it is not injective, since $1\neq [x, yxy^{-1}]$ belongs to the kernel. Since they are finitely generated and not Hopf, then  these groups are not residually finite and hence not parafree. 
\end{eg}

\begin{cor} \label{hnntree} Let $H\leq G\hn{A}$ be a subgroup such that $H$ intersects trivially all the conjugates of $G$. Then $H$ is free.  
\end{cor}

\chapter{Profinite groups} \label{profinitsection}

Let $I$ be set of subscripts with a partial order $\geq$ such that for every $i, j\in I$, there exists $k\in I$ such that $k\geq i$ and $ k\geq j$. An {\bf inverse system of groups} is a collection of groups $G_i$ together with a collection of maps $\{\phi_{ji}: j\geq i\in I\}$ such that 
\begin{enumerate}
    \item For every $j\geq i$, $\phi_{ji}: G_j\lrar G_i$ is a group homomorphism. 
    \item For every $i$, $\phi_{ii}=\id_{G_i}$.
    \item For every $k\geq j\geq i$, $\phi_{ki}=\phi_{ji}\circ \phi_{kj}$. 
\end{enumerate}
The {\bf inverse limit} of this inverse system of groups, denoted by $\varprojlim_{i\in I} G_i$, is defined by
\begin{equation}\label{prodefi} \varprojlim_{i\in I} G_i= \left\{(g_i)_i \in \prod_{i\in I} G_i: \mbox{$\phi_{ji}(g_j)=g_i$\, \, for all $j\geq i$} \right\}.\end{equation}
The group operation on $\varprojlim G_i$ is inherited from the natural group structure of $\prod_{i\in I} G_i$.  

A {\bf profinite group} is an inverse limit of finite groups and they arise naturally in infinite Galois theory. It can be seen that an automorphism of $\ove{\Q}$ restricts to an automorphism of each finite Galois extension $L/\Q$. Reciprocally, an automorphism of $\ove{\Q}$ is made out of elements of the groups $\Gal (L/\Q)$. However, $\Gal(\ove{\Q}/\Q)$ is not the whole $\prod_L \Gal(L/\Q)$, since the later collections of automorphisms should satisfy some consistency conditions in the form of (\ref{prodefi}) to be the restrictions of an automorphisms of $\ove{\Q}$. More precisely,  the absolute Galois group $\Gal (\ove{\Q}/\Q)$ of $\Q$ would be a profinite group, seen as the inverse limit of the finite groups $\Gal (L/\Q)$, where $L/\Q$ ranges over  finite Galois extensions with partial order $E\supseteq L$ and with restriction homomorphisms $r_{E L}: \Gal(E/\Q)\lrar \Gal (L/\Q)$. 

 From the point of view of category theory, the group $\varprojlim G_i$, endowed with the canonical projections $\phi_j: \varprojlim G_i\lrar G_j$, is defined by the following {\bf universal property}: For any group $G$ and any group homomorphisms $\psi_i: G\lrar G_i$ such that, for any $k\geq j$, the diagram 
 \begin{equation*} 
     \begin{tikzcd}
       & G_k\ar[dd, "\phi_{kj}"] \\
       G \ar[ur, "\psi_k"] \ar[dr, "\psi_j"] & \\
        & G_j
     \end{tikzcd}
 \end{equation*}
 is commutative; there exists a unique group homomorphism $\psi: G\lrar \varprojlim G_i$ such that the diagram 
 \begin{equation} \label{uniinverse}
     \begin{tikzcd}
       & & G_k\ar[dd, "\phi_{kj}"] \\
       G  \ar[r, dashrightarrow, "\psi"] \ar[urr, bend left, "\psi_k"] \ar[drr, bend right,  "\psi_j"] & \varprojlim G_i  \ar[ur, "\phi_k"] \ar[dr, "\phi_j"] & \\
       & & G_j
     \end{tikzcd}
 \end{equation}
  is commutative. This gives us a practical way to  construct homomorphisms between profinite groups from first principles. 

On the other side, profinite groups also have a rich structure as topological groups. We equip the groups $G_i$ with the discrete topology, the group $\prod G_i$ with the product topology and the inverse limit $\varprojlim G_i$ with the subspace topology. As an application of Tychonoff's theorem, these groups are compact. This topology on $\varprojlim G_i$ is named the {\bf profinite topology}; and, in the context of Galois theory, the Krull topology; which is crucial in the Fundamental Theorem of Galois theory.

Going back to the universal property of the diagram (\ref{uniinverse}), if the initial $G$ is a topological group and the maps $\psi_k$ are continuous; then the resulting $\psi: G\lrar \varprojlim G_i$ would be continuous, too.  

In the class of topological groups, there are three defining properties of the subclass of profinite groups, namely being Hausdorff, compact and totally disconnected. In particular,  these are far from most topological groups that arise in geometry.

The richness of profinite groups lies in the interface of their algebraic and topological features. Our main interest in these groups is that their structure captures {\bf residual properties} of abstract groups in the form of {\bf profinite invariants}.

We now extend upon the latter remark. A
 {\bf variety of groups} is a class of groups $\CC$ that verifies the following conditions.
\begin{enumerate}
    \item (Closed under isomorphism types) If   $C_1\cong C_2$ and $C_2\in \CC$, then $C_1\in \CC$.
    \item (Closed under direct products) If $C_1, C_2\in \CC$, then $C_1\times C_2\in \CC$. 
    \item (Closed under subgroups) If $C_1\leq C_2$ and $C_2\in \CC$, then $C_1\in \CC$.
    \item (Closed under quotients) If $C_1$ is a quotient of $C_2$ and $C_2\in \CC$, then $C_1\in \CC$. 
\end{enumerate}

A few  examples are the varieties of finite groups, finite $p$-groups, nilpotent groups and soluble groups.

Let $\CC$ be a variety. A {\bf pro-$\CC$ group}  is an inverse limit of groups belonging to $\CC$. We describe the most important examples of pro-$\CC$ groups. Let $G$ be an abstract group. Then the  {\bf pro-$\CC$ completion of $G$}, denoted by $G_{\hat{\CC}}$, is the inverse limit of the system of groups 
\[\{N\unlhd G: G/N\in \CC\},\]
with containment $\sub$ as partial order and canonical group homomorphisms $\phi_{N_1, N_2}: G/N_1\lrar G/N_2$ whenever $N_1\sub N_2$.

If $\CC$ is a variety of finite groups, then $G_{\hat{\CC}}$ is an inverse limit of finite groups in $\CC$ and it is also endowed with the profinite topology. If $\CC=\CC_p$ is the variety of finite $p$-groups, we denote its pro-$\CC$ completion by $G_{\hat{p}}=G_{\hat{\CC}}$ and we name it the {\bf pro-$p$ completion of $G$.}

One could also consider pro-$\CC$ completions of varieties which contain infinite groups, such as the variety of nilpotent groups; although these completions lack an interesting topological structure. 

We say that a group $G$ is residually-$\CC$ if the intersection of normal subgroups $N\unlhd G$ with $G/N\in \CC$ is trivial. In other words, $G$ is residually-$\CC$ if for every $1\neq g\in G$, there exists a group homomorphism $f: G\lrar C$, with $C\in \CC$, such that $f(g)\neq 1$.

\begin{prop} \label{Ccat} Let $\CC$ be a variety. The following is true:
 \begin{itemize} 
 \item There is a natural and canonical group homomorphism $\iota_{\hat{\CC}}: G\lrar G_{\hat{\CC}}$. This map injective if and only if $G$ is residually-$\CC$, and it is bijective if and only if $G\in \CC$. 
 \item If $h: G\lrar \KK$ group homomorphism and $\KK$ is a pro-$\CC$ group, then there exists a group homomorphism $ h_{\hat{\CC}}: G_{\hat{\CC}}\lrar \KK$ such that  the diagram 
  \begin{equation*}
     \begin{tikzcd}
      G \ar[d, "\iota_{\hat{\CC}}"] \ar[dr, bend left, "h"] & \\
      G_{\hat{\CC}} \ar[r, dashrightarrow, "h_{\hat{\CC}}"] & \KK 
     \end{tikzcd}
 \end{equation*}
 is commutative.
 \item If $f: G\lrar H$ is a group homomorphism, there exists a  group homomorphism $f_{\hat{\CC}}: G_{\hat{\CC}}\lrar H_{\hat{\CC}}$ such that the diagram
 \begin{equation*}
     \begin{tikzcd}
      G \ar[d, "\iota_{\hat{\CC}}"] \ar[r, "f"] & H\ar[d, "\iota_{\hat{\CC}}"]\\
      G_{\hat{\CC}} \ar[r, dashrightarrow, "f_{\hat{\CC}}"] & H_{\hat{\CC}}
     \end{tikzcd}
 \end{equation*}
 is commutative.
 \item If, in addition, $\CC$ is a variety of finite groups, then $\iota_{\hat{\CC}}$ has dense image and  the maps  $h_{\hat{\CC}}$ and $f_{\hat{\CC}}$ can be chosen to be continuous in exactly one way. 
 \end{itemize}
 \end{prop}
 When we work in the category of pro-$p$ groups, the maps $\iota_{\hat{\CC}}$, $h_{\hat{\CC}}$ and $f_{\hat{\CC}}$ of \cref{Ccat} will be denoted by $\iota_{\hat{p}}$, $h_{\hat{p}}$ and $f_{\hat{p}}$, respectively. 
 We observe that $G_{\hat{\CC}}$ is residually-$\CC$. 
 In order to inspect the isomorphisms types of quotients of $G$ belonging to a certain variety $\CC$; or to study whether $G$ is residually-$\CC$; it is natural to consider its pro-$\CC$ completion $G_{\hat{\CC}}$.

\begin{thm}[Dixon, Formanek, Poland, Ribes \cite{Dix82}]\label{proC} Let $\CC$ be a variety of finite groups. Two abstract groups $\Ga$ and $\Lambda$ have the same class of isomorphism types of quotients belonging to $\CC$ if and only if $\Ga_{\hat{\CC}}\cong \Lambda_{\hat{\CC}}$. 
\end{thm}

We remark that whenever we talk about morphisms between inverse limits, we talk about group homomorphisms; and, if these are also profinite, we additionally require morphisms to be continuous maps.  All group homomorphisms between profinite groups in this exposition will be naturally continuous.  Still, it is worth mentioning that sometimes this is not really an issue. For instance, since profinite groups are compact and Hausdorff, the inverse of a continuous bijective group homomorphism will always be  continuous. 

There are deeper principles in regard to the algebraic and topological properties of profinite groups. A profinite group $G$ is said to be {\bf topologically finitely generated} if there is a finite subset $S\sub G$ that generates a dense subgroup of $G$. A deep theorem due to Nikolov and Segal \cite{Nik06}, which we shall not use, states that all finite-index subgroups of a profinite group are open. As a consequence, a group homomorphism $f: G\lrar H$ between a profinite $G$ and $H$, where $G$ is finitely generated, must be continuous. This result is significantly easier in the category of $p$-groups and it is a classical observation of Serre.

\begin{prop} \label{nikprop} Let $\GG$ be a topologically finitely generated pro-$p$ group. Then any finite-index subgroup $\HH$ of $\GG$ has $p$-power index and it is open.
\end{prop}


We do will not develop the theory of profinite groups, though we would like to recall a few fundamental facts of this theory.  Since we are only going to work with pro-$p$ groups, we will still recollect some of the important features and tools of the theory of this particular subclass of profinite groups. We refer the reader to the books \cite{Dix99}, \cite{Rib00} and \cite{Wil98}  for a thorough treatment of this topic. 

We start recalling a fundamental observation of groups that are constructed as inverse limits. Group homomorphisms between inverse limits need not be induced from a homomorphism between their inverse systems. More precisely,  a group homomorphism $\varprojlim G_i\lrar \varprojlim H_i$  need not be induced from a sequence of group homomorphisms $G_i\lrar H_i$. Moreover, a group can be the inverse limit of many different inverse systems of groups. 

\begin{prop} \label{filtration}
Let $\GG$ be a profinite group. Let $\{\NN_n\}_n$ be a decreasing sequence of open subgroups with trivial intersection. Then the canonical map $ \GG\lrar \varprojlim_{n} \GG/\NN_n$ is an isomorphism. 
\end{prop}
\begin{proof}  The map is injective because $\bigcap \NN_n=1$. On the other side,  the image is closed since this is a map between Hausdorff compact topological sets. In addition, the image surjects into each factor $\GG/\NN_n$ of the inverse limit, so the image is also dense. This implies that $\phi $ is bijective, as we wanted. 
\end{proof}
Since we are going to consider {\bf closed subgroups and quotients of pro-$p$ groups}, we sketch how one can  make sense of them in the category of pro-$p$ groups. Let $\GG=\varprojlim G_i$ be a pro-$p$ group with projection maps $p_j: \GG\lrar G_j$. We can assume without loss of generality that each $p_j$ is surjective since, otherwise, the canonical map $\GG \lrar \varprojlim p_i (G_i)$ would be an isomorphism and the latter inverse limit has surjective projections $p_j$. Let $\HH\leq \GG$ be a closed subgroup and let $\NN\leq \GG$ be a normal closed subgroup. Then there are canonical continuous isomorphisms $\HH\lrar \varprojlim p_i(\HH)$ and  $\GG/\NN\lrar \varprojlim G_i/p_i(\NN)$, where $\HH$ is endowed with the subspace topology; $\GG/\NN$ is endowed with the quotient topology; and both $ \varprojlim p_i(\HH)$ and $ \varprojlim G_i/p_i(\NN)$ are endowed with the profinite topology.  This way,  both $\HH$ and $\GG/\NN$ are naturally pro-$p$ groups in the sense of inverse limits. 

These notions of subgroup and quotient enjoy all expected and desirable properties in the category of pro-$p$ groups. For example, there is a first isomorphism theorem. 

\begin{lemma}\label{first} Let $\GG$ and $\HH$ be pro-$p$ groups and let $f: \GG\lrar \HH$ be a continuous epimorphism.  For any closed $\NN\n \ker f$, the induced map $\GG/\NN\lrar \HH$ is continuous and verifies and makes  the following diagram
\begin{equation}
    \begin{tikzcd}
    \GG \ar[d] \ar[dr, "f"] &\\
    \GG/\NN\ar[r, dashrightarrow] & \HH
    \end{tikzcd}
\end{equation}
commutative. Moreover, if $\NN=\ker f$ then the induced $\GG/(\ker f)\lrar \HH$ is an isomorphism of pro-$p$ groups. 
\end{lemma}

We also remind that a subgroup of a pro-$p$ group is {\bf open} if and only if it is closed and has finite index. As a consequence, given an open normal subgroup $\HH$ of a pro-$p$ group $\GG$, the quotient $\GG/\HH$ is a finite $p$-group and it is endowed with the discrete topology. 

There is a particularly valuable case of \cref{filtration}. Let $\GG$ have finite rank.  Since $\GG$ is residually $p$, then we would have that $\GG$ is isomorphic to $ \varprojlim_{n} \GG/\gamma_{n}^{(p)} \GG $, where $\gamma_{n}^{(p)}\GG$ are the $p$-lower central series of $F$. These are defined recursively by $\gamma_{1}^{(p)}\GG=\GG$ and  $\gamma_{n+1}^{(p)}\GG=\GG^p\, [\GG, \gamma_{n}^{(p)}\GG]$ for $n\geq 1$. One can inductively use \cref{nikprop} to check that these are, in fact, open normal subgroups of $\GG$. Similarly, we denote by $\ga_k\GG$  the lower central series of $\GG$, which are defined recursively by $\ga_1\GG=\GG$ and $\ga_{k+1}\GG=[\ga_k\GG, \GG]$ for $k\geq 1$. We know that these are closed normal subgroups of $\GG$ (see \cite[Proposition 1.9 and Exercise 17]{Dix99}). 

This special filtration makes the topology of pro-$p$ groups special among profinite groups. However, the key aspect of pro-$p$ groups is the simple characterisation of their Frattini subgroup.

\begin{defi} Let $\GG$ be a pro-$p$ group. We define its {\bf Frattini subgroup $\Phi(\GG)$} to be the intersection of all its maximal  closed subgroups. 
\end{defi}

The interest of the Frattini subgroup, which is clearly normal and closed, is that it allows to reduce the algebraic-topological question of group generation to linear algebra. 

\begin{prop} \label{Fra1} Let $S$ be a subset of $\GG$. Then the following statements are equivalent. 
\begin{enumerate}
    \item $S$ topologically generates $\GG$.
    \item $S\Phi(\GG)$ topologically  generates $\GG$. 
    \item $S\Phi(\GG)/\Phi(\GG)$  topologically generates $\GG/\Phi(\GG)$. 
\end{enumerate}
In particular, $d(\GG)=d(\GG/\Phi(\GG))$. 
\end{prop}

\begin{prop} \label{Fra2} The Frattini subgroup is 
\[\Phi (\GG)= \ove{\GG^p\, [\GG, \GG]}, \]
so $\GG/\Phi(\GG)$ is a $p$-elementary abelian group. 
 Moreover, if $\GG$ is topologically finitely generated, then 
\[\Phi (\GG)= \GG^p \, [\GG, \GG]\, \, \, \mbox{and}\, \,  \, \GG/\Phi(\GG)\cong \F_p^{d(\GG)}.\]
\end{prop}

 Let $\Ga$ be a finitely generated abstract group with pro-$p$ completion $\Ga_{\hp}$. Given a finite $p$-group $Q$, there is a bijective correspondence
\begin{equation}\label{correspondence}
    \{\mbox{epimorphisms $\Ga\lrar P$}\}\longleftrightarrow \{\mbox{epimorphisms $\Ga_{\hp}\lrar P$}\}
\end{equation}
which assings to each epimorphism $f: \Ga\lrar P$ the epimorphsm $f_{\hp} :\Ga_{\hp}\lrar P$. Reciprocally, given an epimorphism $g: \Ga_{\hp}\lrar P$, it will be continuous by \cref{nikprop}. So the map $g\circ \iota_{\hp}: \Ga\lrar P$ is a surjective homomorphism.
In particular, given any $p$-elementary abelian group $A\cong \F_p^n$, the previous assignment produces a bijective correspondence 
\begin{equation}
    \{\mbox{epimorphisms $\Ga/\Ga^p \, [\Ga, \Ga]\lrar A$}\}\longleftrightarrow \{\mbox{epimorphisms $\Ga_{\hp}/\Ga_{\hp}^p\, [\Ga_{\hp}, \Ga_{\hp}] \lrar A$}\}.
\end{equation}
This implies that the dense canonical map $\Ga\lrar \Ga_{\hp}$ induces an isomorphism of elementary abelian $p$-groups 
\begin{equation} \label{proabiso} \Ga/\Ga^p \, [\Ga, \Ga]\lrar \Ga_{\hp}/\Ga_{\hp}^p\, [\Ga_{\hp}, \Ga_{\hp}]. \end{equation}

We encompass, as a corollary of propositions \ref{Fra1} and \ref{Fra2}, part of the utility of the Frattini group in the following statement. 

\begin{cor} \label{Fra12} Let $\GG$ have finite rank. Then $S\sub \GG$ generates $\GG$ topologically if and only if $S$ generates $\GG$ modulo $\GG^p \, [\GG, \GG]$. In particular, a  homomorphism $f: \GG\lrar \HH$ between finitely generated pro-$p$ groups $\GG$ and $\HH$ is surjective if and only if the induced 
\[f: \frac{\GG}{\GG^p \, [\GG, \GG]}\lrar \frac{\HH}{\HH^p \, [\HH, \HH]}\]
is surjective. 
\end{cor}

These results have surprising implications and constitute one of the reasons for which pro-$p$ groups are easier to work with than with abstract groups or more general profinite groups. 

 Notice that free groups $F$ do also have the property that $d(F)=\dim_{\p} F/F^p \, [F, F]$, but not for analogous reasons. As opposed to what happens in pro-$p$ groups, lifts of generators of $F/F^p \, [F, F]$ are not necessarily generators of $F$.

\section{Solving equations} \label{solving}
In the theory of pro-$p$ groups, we talk about topological generation, rather than abstract generation. As a consequence, elements of pro-$p$ groups cannot be generally described by finite words. In this section, our aim is to develop an elementary way of describing implicitly elements in pro-$p$ groups. This section is perfectly avoidable. It just suggests and alternative combinatorial way of thinking about elements in pro-$p$ groups.

Let $F_n$ be the free group on $X=\{x_1, \dots, x_n\}$ and let $\omega\in F_n$ be a word in $X$. Defining a group homomorphism from a one-relator group $\lan X| \omega\ran$ to another group $\GG$ is equivalent to finding $h_1, \dots, h_n\in \GG$ such that $\omega(x_1=h_1, \dots, x_n=h_n)=1$ in $\GG$. If $H$ is a topological group, there are methods to ensure non-explicit solutions $h_1, \dots, h_n$ of $\omega$. Suppose that $\GG$ is a pro-$p$ group. Take $X=\{x_1, x_2, x_3\}$ and consider the word $\omega=x_1[x_2, x_1][x_3, x_2]$. Take any $x_2=b\in \GG$ and $x_3=c\in \GG$. We want to find $a\in \GG$ that completes a solution $(a, b, c)$ of $\omega$. If $(a, b, c)$ is a solution, then 
\[a=[b, c][a, b]=[b, c][[b, c], [a, b]], b]= \cdots=f^{(n)}(a)=\cdots,\]
where $f^{(n)}$ is the $n$-th iteration of $f(z)=[b, c][z, b]$. It is tempting to take any initial $a_0\in \GG$ and define $a=\lim_n f^{(n)} (a_0)$ in $\GG$, since the required $a$ is a fixed point of $f$. However, we do not know whether the limit  $\lim_n f^{(n)} (a_0) $ exists or not. In this section, we introduce basis elementary tools to study when this sequence is convergent and when the resulting limit is independent of the initial $a_0\in \GG$. 

 There are mainly two ways to prove the existence of such $a$ by using the compactness of $\GG$. These are encoded in \cref{equations} and \cref{fixpoint}.

\begin{lemma} \label{equations}
Let $\omega$ be a word on $n$ letters. Let $\GG$ be the inverse limit $\varprojlim G_i$ of finite groups $G_i$. Let $a\in \GG$. Suppose that the equation $w(x_1, \dots, x_n)=p_i (a)$ has a solution in $G_i$ for every $i\in I$. Then the equation  $w(x_1, \dots, x_n)=a$ has a solution in $\GG$. Additionally, if the previous equation has unique solution modulo $G_i$ for each $i\in I$, then the solution is also unique in $\GG$. 
\end{lemma}

In particular, this allows us to extract roots in pro-$p$ groups. 
\begin{prop}[Divisibility of pro-$p$ groups] \label{pro-$p$ div} Let $\GG$ be a pro-$p$ group, let $a\in \GG$ and let $n$ be an integer not divisible by $p$. 
Then there exists a unique $b\in \GG$ such that $a=b^n$. 
\end{prop}

The application of \cref{equations} works to some extend and we are going to introduce another method for solving more general type of equations in pro-$p$ groups, provided a mild condition modulo $p$.

\begin{defi}
Let $\GG$ be a profinite group and let $\mathcal{N}$ be a collection of open normal subgroups with trivial intersection. For each $\NN\in \mathcal{N}$, we construct the continuous function $d_{\NN}: \GG\lrar \{0, 1\}$ which verifies that $d_{\NN}(x)=0$ if $x\in \NN$ and $d_{\NN}(x)=1$ if $x\notin \NN$. 
\end{defi}

\begin{prop} \label{prodistance}
Let $\GG$ be a profinite group of finite rank. Then the collection of all finite-index subgroups is countable. Let $\mathcal{N}$ be a subcollection of finite intersection. Then $\GG$ is isomorphic to  $\varprojlim_{\NN\in\, \mathcal{N}} \GG/\NN$. In addition, let us consider an ordering of $\mathcal{N}$ by $\N$. We denote $d_n= d_{U_n}$. Then the following function $d: \GG\times \GG\lrar [0, 1]$ is a metric for the topology of $\GG$:
\[d(x, y)=\sum_{n\in \N} 2^{-n}\, d_{n}(x y^{-1}).\]
\end{prop}

Once we have built an explicit metric for $\GG$, we will see how a fixed-point theorem from functional analysis can be used to solve equations in pro-$p$ groups $\GG$. 

\begin{defi}
Let $X$ be a metric space of distance $d$. We say that a function $f: X\lrar X$ is weakly contractive if it verifies that $d(f(x), f(y))<d(x, y)$ whenever $x\neq y$. 
\end{defi}

\begin{lemma}[Fixed point theorem in compact metric spaces] \label{fixpoint} 
Let $X$ be a compact metric space of distance $d$. Let $f: X\lrar X$ be a weakly contractive function.  Then $f$ has a unique fixed point. 
\end{lemma}
\begin{proof}
The function $q: X\rar \R$ defined by $x\mapsto d(x, f(x))$ is continuous since both $d$ and $f$ are clearly continuous. Let $m$ be its minimum, which is finite and is attained at some $x_0\in X$ by the compactness of $X$. If $x_0\neq f(x_0)$, we would have that $q(f(x_0))=d(f(x_0), f(f(x_0)))<d(x_0, f(x_0))=q(x_0)$, contradicting the definition of $x_0$. So $x_0$ is a fixed point. If $x_1$ and $x_2$ where different fixed points we would reach the contradiction $d(x_1, x_2)=d(f(x_1), f(x_2))<d(x_1, x_2)$. Thus $x_0$ is the unique fixed point. 
\end{proof}

\begin{prop} \label{Uchar}
Let $\GG$ be a topologically finitely generated pro-$p$ group, endowed with some metric from \cref{prodistance}. A function $f: \GG \lrar \GG$ is weakly contractive if and only if for every $x\neq y$ there exists an open normal subgroup $\NN\leq \GG$ such that $xy^{-1}\notin \NN$ and $f(x)f(y)^{-1}\in \NN $. 
\end{prop}

We now give many examples of contractive functions on pro-$p$ groups.

\begin{prop} \label{contractiveeg}
Let  $\omega$ be a word in $\{x_1, \dots, x_n\}$ such that $\omega_{x_1}$ is divisible by $p$. Then for every finitely generated pro-$p$ group $\GG$ and every $c_2, \dots, c_n\in \GG$, the  function $f: \GG\lrar \GG$ defined by $x\mapsto \omega (x, c_2, \dots, c_n)$ is weakly contractive.
\end{prop}

\begin{proof}
With the aim of applying \cref{Uchar}, we will ensure, for each $x\neq y$, the existence of an open normal subgroup $\NN$ such that, naming $d=xy^{-1}$, it is true that $d\notin \NN$, $d^p\in \NN$ and $[d, c]\in \NN$ for every $c\in \GG$.  Then it would be clear that $f(dy)f(y)^{-1}\equiv d^{\omega_{x_1}}\equiv 1 \mod \NN$. So $f(x)f(y)^{-1}\in \NN$ and $xy^{-1}\notin \NN$, as we wanted. 

Take $x\neq y \in \GG$. Since $d=xy^{-1}\neq 1$, there exists a continuous homomorphism $\phi : \GG\rar P$ to a finite $p$-group $P$ such that $\phi (d)\neq 1$.  Since $P$ is a finite $p$-group, there exists a unique $k$ such that $\phi(d) \in \gamma_{ k-1, p} P\setminus \gamma_{k, p} P$. Now consider the projection $p: P\rar P/\gamma_{p, k} P$. It is clear that $\NN=\ker( p\circ \phi)$ verifies the previously claimed conditions. 
\end{proof}

These families of contractive functions allow us to produce many examples of equations in $G$ with unique solutions. 

\begin{prop} \label{solvingeq}
Let $\omega$ be a word in $\{x_1, \dots, x_n\}$ such that $\omega_{x_1}$ is coprime with $p$. Then for every finitely generated pro-$p$ group $\GG$ and every $c_2, \dots, c_n\in \GG$, the equation $\omega (x, c_2, \dots, c_n)=1$ has a unique solution $x\in \GG$.  
\end{prop}
\begin{proof}
There exists a positive integer $m\in \Z$, coprime with $p$, such that $m\omega_{x_1}+1$ is divisible by $p$. As a consequence, the word $\omega'(x_1, \dots, x_n)=x_1\, \omega (x_1^{m}, x_2, \dots, x_n)$ lies under the assumptions of \cref{contractiveeg}. By \cref{pro-$p$ div}, the equation $\omega (x, c_2, \dots, c_n)=1$ has a unique solution $x\in \GG$ if and only if $y=\omega' (y, c_2, \dots, c_n)$ has a unique solution $y\in \GG$, and the later claim is ensured by \cref{fixpoint} and \cref{contractiveeg}. 
\end{proof}

We will later use this result to characterise which one-relator groups have free pro-$p$ completion. 

\section{Free products of pro-$p$ groups} 
We have two purposes for this subsection. One the one hand, we define what is meant for a pro-$p$ group to be free and we also discuss some methods for ensuring the pro-$p$ completion of a group to be free. On the other hand, we discuss the construction of the free product in the category of pro-$p$ groups. In the category of abstract groups, this would be the usual free product $*$, though in the later category this is more delicate. 

Let $\GG$ be a pro-$p$ group and let $n$ be a non-negative integer. We say that {\bf $\GG$ is a free pro-$p$ group of rank $n$} if it has a subset $X\sub \GG$  of size $n$ such that the pair $(\GG, X)$  verifies the  following universal property. 
Let $\KK$ be a pro-$p$ group and let $f_0: X\lrar \KK$ be set-theoretical map. Then there exists a unique continuous homomorphism $f: \GG \lrar \KK$ such that the diagram 
\begin{equation}\label{freex}
\begin{tikzcd}
 X\ar[dr, "f_0"] \ar[r, hookrightarrow] & \GG\ar[d, "f"]\\
  & \KK
\end{tikzcd}
\end{equation}
is commutative. In this case, we also say that {\bf $\GG$ is freely generated by $X$} or, simply, that $\GG$ is free on $X$.

\begin{prop}\label{freepropgroups} There exists a pro-$p$ group $\GG$ verifying the previous universal property (\ref{freex}) for some $X\sub \GG$ of size $n$ and, up to isomorphism, this group is unique. 
\end{prop}

\begin{proof} We start verifying the uniqueness with an standard argument.

{\bf Uniqueness:}  Suppose that we are given two pairs $(\GG, X)$ and $(\GG', X')$ verifying the universal property. We can take a bijection $\phi: X\lrar X'$. There exist continuous group homomorphisms $f: \GG\lrar \GG'$ and $f': \GG'\lrar \GG$ such that $f(x)=\phi(x)$, for all $x\in X$; and $f'(x')=\phi^{-1}(x')$, for all $x'\in X'$.  We can now verify that  $f'\circ f=\id_{\GG}$ and $f\circ f'=\id_{\GG'}$. Both are analogous. For the former, notice that $f'\circ f: \GG\lrar \GG$ is a continuous homomorphism whose restriction to $X$ is $f'( f(x))=f'(\phi(x))=\phi^{-1}(\phi(x))=x=\id (x)$. By the uniqueness of the universal property, this implies that $f'\circ f\equiv \id_{\GG}$.

{\bf Existence:} Let $F(X)$ be a free group on generators $X=\{x_1, \dots, x_n\}$. We denote $\FF(X)$ to be its pro-$p$ completion. For our argument, we simply need \cref{Ccat}. We know that the canonical map $\iota_{\hat{p}}: F(X)\inc \FF(X)$ is injective, since $F$ is residually-$p$. So we can identify $F$ with its isomorphic copy in $\FF$. We now prove that the pair $(\FF(X), X)$ verifies the universal property of \cref{freex}. Let $\KK$ be a pro-$p$ group and let $f_0: X\lrar \KK$ be a set-theoretic map. Since $F(X)$ is free with free generators $X$, there exists a unique group homomorphism $f_1: F(X)\lrar \KK$ such that $f_1(x)=f_0(x)$ for all $x\in X$. We can now consider its pro-$p$ extension $f=(f_1)_{\hat{p}}: \FF(X)\lrar \KK$, which is a continuous homomorphism. The following diagram
\begin{equation}
    \begin{tikzcd}
     X \ar[rrrr, red, hookrightarrow, bend left] \ar[rr, hookrightarrow] \ar[ddrrrr, red,  bend right,  "f_0"] & & F(X) \ar[ddrr, bend right, "f_1"] \ar[rr, "\iota_{\hat{p}}"] && \FF(X) \ar[dd, red, dashrightarrow, "f"] \\
     & & & &\\
     & & & & \KK 
    \end{tikzcd}
\end{equation}
is commutative. The subdiagram that is highlighted in red exhibits that the resulting $f: \FF(X)\lrar \KK$ has the required properties of (\ref{freex}). It would remain to prove that this map $f$ is unique. Let us consider another such map $f': \FF(X)\lrar \KK$ and its corresponding $f_1'=f'\circ \iota_{\hat{p}}: F(X)\lrar \KK$. By assumption, the following diagram 
\begin{equation}
    \begin{tikzcd}
     X \ar[rrrr,  hookrightarrow, bend left] \ar[rr, blue, hookrightarrow] \ar[ddrrrr, blue,  bend right,  "f_0"] & & F(X) \ar[ddrr, blue,  bend right, "f_1'"] \ar[rr, "\iota_{\hat{p}}"] && \FF(X) \ar[dd,  dashrightarrow, "f'"] \\
     & & & &\\
     & & & & \KK 
    \end{tikzcd}
\end{equation}
would be commutative. We study the subdiagram that is highlighted in blue. Since $X$ generates $F(X)$, the extensions $f_1'\equiv f_1$ must be identical. Secondly, the continuous homomorphisms $f, f': \FF(X)\lrar \KK$ would coincide in the dense subgroup  $\iota_{\hat{p}} (F(X))$. So $f'\equiv f$. This finishes the proof. 
\end{proof}

\begin{nota}
We denote by  $\FF_n$ the free pro-$p$ group of rank $n$. 
\end{nota} 
Analogously, we could have introduced the notion of free pro-$p$ group of any cardinal. However, we stick to the case of finitely generated pro-$p$ groups. The only free such groups are the $\FF_n$. Observe that $\FF_1$ is isomorphic to the group of $p$-adic integers $\Z_p$.

In terms of the universal property (\cref{freex}), we can rephrase the property of a  one-relator group having free pro-$p$ completion in terms of a problem of ensuring a unique solution to an equation related to the defining relator $\omega$. The latter reminds of \cref{solving}.

\begin{lemma} \label{pro-p-freeness} Let $\omega$ be a word in $X=\{x_1, \dots, x_n\}$. Suppose that for any finitely generated  pro-$p$ group $\GG$ and for all $c_2, \dots, c_n\in \GG$, the equation $\omega (x_1, c_2, \dots, c_n)=1$ has a unique solution $x_1\in \GG$. Then the one-relator group $\ti{G}=\lan X| \omega\ran$
verifies that its pro-$p$ completion $\ti{G}_{\hat{p}}$ is isomorphic to $\FF_{n-1}$. 
\end{lemma}

As a consequence, we can give a simple criterion for the pro-$p$ completion of a one-relator group to be free.

\begin{prop} \label{free pro-p completion}
Let $\omega$ be a word in $X=\{x_1, \dots, x_n\}$ and suppose that, for some $1\leq k \leq n$,  $\omega_{x_k}$ is coprime with $p$. Denote by $G=\lan X|\omega\ran$ the one-relator group of defining relation $\omega$. Then $G_{\hat{p}}\cong \FF_{n-1}$.
\end{prop}
\begin{proof}
This follows from \cref{pro-p-freeness} and \cref{solvingeq}.
\end{proof}

We leave stated two elementary lemmas on free pro-$p$ groups, which are direct consequences of the the above arguments and the properties of the Frattini subgroup, collected in propositions \ref{Fra1} and \ref{Fra2}.

\begin{cor}[Free bases of free pro-$p$ groups] \label{freebasis} Let $X$ be a subset of $\FF_n$.  Then $X$ freely generates $\FF_n$ if and only if the reduction map $X\lrar \FF/\Phi(\FF)$ is injective and  $\{x\, \Phi(\FF)\, : \, x\in X\}$ is a basis of the $\p$-vector space 
\[\frac{\FF}{\Phi(\FF)}\cong \F_p^n.\]
\end{cor}

\begin{cor}[Strong Hopf property] \label{Hopf2} Let $\FF$ and $\FF'$ be two finitely generated  pro-$p$ groups and suppose that $\FF'$ is free. Let $f: \FF\lrar \FF'$ be a continuous epimorphism. Then $d(\FF)\geq d(\FF')$, with equality if and only if $f$ is an isomorphism. 
\end{cor}
The proof of the last result does also require the Hopf property on finitely generated pro-$p$ groups, which is not hard to establish.

We have worked with a universal property in order to compare two pro-$p$ completions. We will refine this method in the next section. For the moment, we will introduce another important construction which generalises the free pro-$p$ group, namely the free product of pro-$p$ groups.

Let $\GG_1, \dots, \GG_n$ be pro-$p$ groups. We are interested in constructing their free product (coproduct) in the category of pro-$p$ groups. More precisely, we are interested in a pro-$p$ group $\GG$, endowed with $n$ canonical continuous homomorphisms $\phi_k: \GG_k\lrar \GG$, that verifies the following universal property: For any pro-$p$ group $\KK$ and any continuous homomorphisms $\psi_k: \GG_k\lrar \KK$, there exists a unique continuous homomorphisms $\psi: G\lrar \KK$ such that, for any $k$, the diagram 
\begin{equation}\label{unicop}
    \begin{tikzcd}
     \GG_k \ar[r, "\phi_k"] \ar[dr, "\psi_k"] & G\ar[d, dashrightarrow, "\psi"]\\
      & \KK
    \end{tikzcd}
\end{equation}
is commutative. That is, $\psi\circ \phi_k=\psi_k$. 

\begin{prop}  \label{coprod} Given pro-$p$ groups $\GG_1, \dots, \GG_n$, there exists a unique pro-$p$ group $\GG$ verifying the universal property of (\ref{unicop}). We denote this group $\GG$ by $\coprod_i \GG_i$ and by $\phi_k : \GG_k\lrar \coprod_i \GG_i$ the canonical maps. 
\end{prop}
\begin{proof} We first precise what we mean by saying that such $\GG$ is unique and we prove that, indeed, this object is unique. Then we give an explicit construction of it. 

{\bf Uniqueness:} Suppose that there exist two such objects $\GG$ and $\GG'$ endowed with canonical maps $(\phi_k)_k$ and $(\phi_k')_k$, respectively. We claim that there exists a group isomorphism $\psi: \GG\lrar \GG'$ such that $\psi\circ \phi_k=\phi_k'$ for every $k$. This is what we mean when we say that the object $\GG$ is unique. We are going to see how to produce this isomorphism $\psi$ from the universal property. 
Applying the universal property of (\ref{unicop}) that verifies $\GG$ to the data ($\KK=\GG'$, $\phi_k': \GG_k\lrar \GG'$) and the universal property that verifies $\GG'$ to the data $(\KK=\GG$, $\phi_k: \GG_k\lrar \GG)$, we get homomorphisms $\psi$ and $\psi'$ such that, for every $k$, the diagram
\begin{equation} \label{cop0}
    \begin{tikzcd}
     \GG_k \ar[rr, red, "\phi_k"] \ar[ddrr, red,  bend right, "\phi_k'"] & & \GG\ar[dd, red,  dashrightarrow, bend left, "\psi"] & & \\
      & & & & \\
       & & \GG' \ar[uu, blue,  dashrightarrow, bend left, "\psi'"] & & \ar[ll, blue,  "\phi_k'"] \GG_k,  \ar[uull, blue, bend right, "\phi_k"]  
    \end{tikzcd}
\end{equation}
is commutative. In other words, $\psi \circ \phi_k=\phi_k'$ and $\psi'\circ \phi_k'=\phi_k$. It directly follows that the diagram 
\begin{equation*}
    \begin{tikzcd}
     \GG_k \ar[rr, "\phi_k"] \ar[ddrr, "\phi_k"] & & \GG\ar[dd, dashrightarrow, bend right,  "\id"]\ar[dd, dashrightarrow, bend left, "\psi'\circ \psi"] \\
     & & \\
     & & \GG
    \end{tikzcd}
\end{equation*}
is commutative. By the uniqueness of the universal property of $\GG$ with respect to the data $(K=G, \phi_k: \GG_k\lrar \GG)$, this implies that $\psi'\circ \psi=\id$. Analogously, we prove that $\psi\circ \psi'=\id$, and the uniqueness is ensured. 

{\bf Existence:} Consider the free product $\ag=\GG_1* \cdots * \GG_n$ in the category of groups, with canonical maps $\af_k: \GG_k\lrar \ag$. We introduce the family 
\[\mathcal{N}_0 = \{N\n \ag : (\af_k)^{-1}(N)\no \GG_k, \,\,  \, \, \mbox{and\, \, \,  $\ag/N\in \CC_p$}\}.\]
We claim that for any $N_1, N_2\in \mathcal{N}_0$,   $ N_1\bigcap N_2\in \mathcal{N}_0$. In fact,  $N_1\bigcap N_2$ is the kernel of the canonical map 
\[\ag \lrar \ag/N_1\times \ag/N_2.\]
 Now consider the following inverse system of groups:
\[\mathcal{N} = \{\ag/N : N\in \mathcal{N}_0\}.\]
where it is viewed with canonical maps $\phi_{N_1 N_2}: \ag/N_1\lrar \ag/N_2$ whenever $N_1\sub N_2$. Denote its inverse limit by 
\[\GG=\varprojlim_{N\in \mathcal{N}_0} \ag/N.\]
By the universal property \cref{uniinverse} of $\GG$, the canonical projections $\af_k: \ag \lrar \ag $ induce a map $\iota: \ag\lrar \GG$, which has dense image. 

We  think of $\ag$ as a  topological group that has $\mathcal{N}_0$ as a basis of neighbourhoods of $1\in \ag$. This way, we can view $\GG$ as the completion of $\ag$ with respect to this topology and the map $\iota $ is continuous. Similarly, the maps $\phi_k=\iota \circ \af_k: \GG_k\lrar G$ are continuous. We have that $\phi_k$ is continuous because, for each $N\in \mathcal{N_0}$, $\phi_k^{-1}(N)$ is open in $\GG_k$, by assumption. So $\phi_k=\iota \circ \af_k$ is the composition of continuous maps. 

We now want to establish the universal property (\ref{unicop}) for the group $\GG$ and the continuous homomorphisms $\phi_k: \GG_k\lrar G$.  Let $\KK$ be a pro-$p$ group and let  $\psi_k: \GG_k\lrar \KK$ be continuous homomorphisms. Due to the universal property of the inverse limit, condensed in the diagram (\ref{uniinverse}), we can assume that $\KK$ is a finite $p$-group in order to construct a continuous homomorphism $\psi: G\lrar \KK$ such that the diagram (\ref{unicop}) commutes.   

 Since $\KK$ is finite (and hence discrete), the continuity of $\psi_k$ translates into each  $\ker \psi_k$ being open. 
Working in the category of groups, by the universal property of $\ag=\GG_1*\cdots *\GG_n$, there exists a unique $\ap: \ag\lrar \KK$ such that the diagram 
\begin{equation} 
    \begin{tikzcd}
     \GG_k \ar[r, "\af_k"] \ar[dr, "\psi_k"] & \ag \ar[d, dashrightarrow, "\ap"]\\
      & \KK
    \end{tikzcd}
\end{equation}
is commutative.  Furthermore, it can be said that $\ap$ is continuous. To prove this claim, it suffices to check that $\ker \ap$ is open. In fact, we see that $\ker \ap\in \mathcal{N}_0$, since \[\phi_k^{-1}(\ker \ap)=\ker \psi_k\] is open in $\GG_k$ for all $k$. 

So $\ap$ is continuous and can be extended to a continuous map $\psi: \GG\lrar \KK$ in such a way that the diagram
\begin{equation}
    \begin{tikzcd}
     \GG_k \ar[rrrr, red, bend left,  "\phi_k"] \ar[rr, "\af_k"] \ar[ddrrrr, red,  bend right,  "\psi_k"] & & \ag \ar[ddrr, bend right, "\ap"] \ar[rr, "\iota"] && G \ar[dd, red, dashrightarrow, "\psi"] \\
     & & & &\\
     & & & & \KK 
    \end{tikzcd}
\end{equation}
is commutative. The commutative subdiagram that is highlighted red shows that the resulting $\psi: G\lrar \KK$ has the desired properties. 

It rests to check that this  $\psi$ is unique. Suppose that $\psi_0$ is a any other such map. We can define $\ap_0: \ag \lrar \KK$ by $\ap_0=\psi_0\circ \iota$ and observe that there are, again,  commutative diagrams of the form  
\begin{equation}
    \begin{tikzcd}
     \GG_k \ar[rrrr, bend left,  "\phi_k"] \ar[rr, blue,  "\af_k"] \ar[ddrrrr, blue,  bend right,  "\psi_k"] & & \ag \ar[ddrr, blue, bend right, "\ap_0"] \ar[rr, "\iota"] && G \ar[dd, dashrightarrow, "\psi_0"] \\
     & & & &\\
     & & & & \KK.
    \end{tikzcd}
\end{equation}
Analysing the subdiagram that is highlighted in blue; we prove that $\ap_0\equiv\ap$, due to the universal property of $\ag$ in the category of groups. Lastly, it is also clear that $\psi_0\equiv\psi$, since they are both continuous and both coincide every $\iota(\af_k(\GG_k))$, which generate the dense subgroup $\iota (\ag)$ of $\GG$. 

This completes the proof of \cref{coprod}.\end{proof}

As we remarked before, it is convenient to think of $\GG_1\coprod \cdots \coprod G_n$ as the completion of the topology  of $\ag=G_1*\cdots *G_n$ generated by the finite-index normal subgroups $N$ such that $\ag/N$ are finite $p$-groups and $N\bigcap \GG_k$ is open in $\GG_k$. 

In any case, the way we are going to work with this object is by means of its universal property (\ref{unicop}). Observe that $\FF_n$ is isomorphic to the $n$-fold coproduct $\Z_p\coprod \cdots \coprod \Z_p$. We give some of the expected properties of the coproduct. 

\begin{cor} \label{cop1} Let $\GG_1, \dots, \GG_n$ be pro-$p$ groups and let $\phi_k: \GG_k\lrar \GG_1\coprod \cdots \coprod \GG_n$ be the canonical maps. Then 
\begin{itemize}
    \item each $\phi_k$ is injective; and 
    \item $\GG$ is topologically generated by $\bigcup_{k} \phi_k(\GG_k)$.
\end{itemize}
\end{cor}
\begin{proof}
Notice that each $\GG_k$ is a retract of $\GG=\GG_1\coprod \cdots \coprod \GG_n$. In fact, there exists continuous homomorphisms $\psi_k: \GG\lrar \GG_k$ such that $\psi_k\circ \phi_k=\id_{\GG_k}$. From this, the first claim follows. The second one is a consequence of the fact that, as we saw in the proof of \cref{coprod}, $\iota(\GG_1*\cdots *\GG_n)$ is dense in $\GG$ and $\GG_1*\cdots *\GG_n$ is generated by $\cup_k \GG_k$.
\end{proof}
Since the canonical $\phi_k: \GG_k\lrar \GG_1\coprod \cdots \coprod \GG_n$ are injective, we will sometimes identify $\GG_k$ with its isomorphic copy in 
 $\GG_1\coprod \cdots \coprod \GG_n$. More than this can be said. The following is \cite[Proposition 9.1.8]{Rib00}. 
 
 \begin{prop}\label{injectiveabs} Let $\GG_1, \dots, \GG_n$ be pro-$p$ groups and consider its free product $\ag=\GG_1*\cdots *\GG_n$ as abstract groups. Then the canonical map 
 \[\iota: \ag\lrar \GG_1\coprod \cdots \coprod \GG_n\]
 is injective. 
 \end{prop}

The Grushko-Neumann theorem in the category of abstract groups is a deep result. Interestingly, in the category of profinite groups, its analogue turns out to be false. However, as we are about to see, in the category of pro-$p$ groups this theorem can be reduced to elementary linear algebra. 

\begin{prop}\label{additivecop} Let $\GG_1, \GG_2$ be two  pro-$p$ groups of finite rank. Then $d\left(\GG_1\coprod \GG_2\right)=d(\GG_1)+d(\GG_2)$. 
\end{prop}
\begin{proof}
Since $\GG_1$ and $\GG_2$ generate topologically $\GG=\GG_1\coprod \GG_2$, it is clear that 
\begin{equation} \label{fr1} d(\GG)\leq d(\GG_1)+d(\GG_2).\end{equation}
Furthermore, by the universal property of $\GG$, there is a continuous homomorphism $f: \GG\lrar \GG_1/\Phi(\GG_1)\times \GG_2/\Phi(\GG_2)$
such that the restrictions to each $\GG_i$ are the canonical projections $\GG_i\lrar \GG_i/\Phi(\GG_i)$. In particular, $f$ is surjective. So 
\begin{equation} \label{fr2} d(\GG)\geq d\left( \GG_1/\Phi(\GG_1)\times \GG_2/\Phi(\GG_2)\right).\end{equation}
Recall that each $\GG_i/\Phi(\GG_i)$ is a $p$-elementary abelian group by \cref{Fra2} and that, by \cref{Fra1}, $d(\GG_i)=d(\GG_i/\Phi(\GG_i))$. Thus
\[d\left( \GG_1/\Phi(G_1)\times \GG_2/\Phi(\GG_2)\right)=d(\GG_1/\Phi(\GG_1)) +d(\GG_2/\Phi(\GG_2))=d(\GG_1)+d(\GG_2),\]
which, in addition to (\ref{fr1}) and (\ref{fr2}), implies that 
$d(\GG)=d(\GG_1)+d(\GG_2)$. 
\end{proof}

\begin{cor} \label{Fra3} Let $\GG_1$ and $\GG_2$ be two  pro-$p$ groups of finite rank and let $\GG=\GG_1\coprod \GG_2$. Then there is a canonical isomorphism 
\[\frac{\GG}{\Phi(\GG)}\lrar \frac{\GG_1}{\Phi(\GG_1)}\times \frac{\GG_2}{\Phi(\GG_2)}.  \]
\end{cor}
\begin{proof}
In the proof of \cref{additivecop}, we defined the canonical surjective map 
\[f: \frac{\GG}{\Phi(\GG)}\lrar \frac{\GG_1}{\Phi(\GG_1)}\times \frac{\GG_2}{\Phi(\GG_2)}. \]
We now know that both finite $p$-elementary abelian groups have the same size by the same proposition. So $f$ is, in fact, an isomorphism. Alternatively, we can explicitly construct an inverse of $f$. 
\end{proof}

\section{Presentations of pro-$p$ completions} \label{presentationsection}

From this point on, we shall discuss how the pro-$p$ completions of an amalgamated products or of an HNN extensions can be described in terms of the pro-$p$ completions of the involved factors. The most natural way to compare pro-$p$ completions  is by working with a defining universal property. Before this, we give an example of what we mean by presenting a pro-$p$ completion. 

{\bf Abstract group presentation:}
Let $X$ be a set and let  $R$ be a subset of the free group $F(X)$ on  $X$. We denote the abstract group 
\[\lan X | R\ran ={F(X)}\big/{\lan \lan R\ran \ran}.\]

{\bf Pro-p group presentation:} Let $X$ be a set and let $\FF(X)$ be the free pro-$p$ group of free generators indexed by the $X$. We denote the pro-$p$ group
\[\lan X | R\ran_p ={\FF(X)}\big/{\overline{\lan \lan R\ran \ran}}.\]

A particular property that we will establish is the following. 
\begin{prop} \label{proppres} Let $G$ be an abstract group of presentation $\lan X | R\ran$. Then there is a canonical isomorphism
\[G_{\hat{p}}\cong \lan X |R\ran_p.\]
\end{prop}

Our starting point, as we mentioned above, is re-defining the canonical map $G\lrar G_{\hat{p}}$ in terms of a universal property. 

\begin{prop} \label{uniprop} Let $\Ga$ be an abstract group. Let $\iota_{\hat{p}}: \Ga\lrar \Ga_{\hat{p}}$ be the canonical map. The data $(\Ga_{\hp}, \iota_{\hp})$ is characterised by the following properties. 
\begin{enumerate}
    \item $\Ga_{\hat{p}}$ is a pro-$p$ group. 
    \item The map $\iota_{\hat{p}}$ has dense image. 
    \item For any pro-$p$ group $\KK$ and every group homomorphism $f: \Ga\lrar \KK$ with dense image, there exists a continuous homomorphism $f_{\hat{p}}: \Ga_{\hat{p}}\lrar \KK$ such that the diagram 
    \begin{equation}
        \begin{tikzcd}
         \Ga \ar[dr, "f"] \ar[r, "\iota_{\hp}"] & \Ga_{\hat{p}}\ar[d, dashrightarrow, "f_{\hat{p}}"]\\
           & \KK
        \end{tikzcd}
    \end{equation}
    is commutative. 
\end{enumerate}
When we say $(\Ga_{\hp}, \iota_{\hp})$ are unique, we mean if $(\HH, \iota)$, where $\HH$ is a pro-$p$ group and $\iota: \Ga\lrar \HH$ is a group homomorphism, is a pair that verifies the three above properties; then there exists an isomorphism of pro-$p$ groups $\al: G_{\hp}\lrar \HH$ such that $\al\circ \iota_{\hp}=\iota$. 
Moreover, by the universal property of the inverse limit (\ref{uniinverse}), it suffices to check the third condition for finite $p$-groups $\KK$. 
\end{prop}

\begin{proof} We already know that $(G_{\hp}, \iota_{\hp})$ verifies the given properties. The proof of the uniqueness of the objects $(G_{\hp}, \iota_{\hp})$ is standard and works the same way as in the uniqueness of the free product of pro-$p$ groups. 
\end{proof}

From this point of view, we can compute pro-$p$ completions of quotients and of free products. Notice that the following result generalises \cref{proppres}.

\begin{prop}[{\bf Pro-$p$ completion of a quotient}] \label{propquo} Let $\Ga$ be an abstract group and let $N\n \Ga$ be a normal subgroup. There is a canonical isomorphism 
\[\left(\Ga/N\right)_{\hp} \cong \Ga_{\hp}\Big/ \ove{\iota_{\hp}(N)}. \]
\end{prop}
\begin{proof}
Consider the canonical group homomorphisms 
\begin{equation*}
\begin{tikzcd}
    \Ga \ar[r, "\iota_{\hp}"] & \Ga_{\hp} \ar[r] & \Ga_{\hp}\Big/ \ove{\iota_{\hp}(N)},
\end{tikzcd}
\end{equation*}
whose composition $\Ga\lrar \Ga_{\hp}\Big/ \ove{\iota_{\hp}(N)} $ is a map with dense image and with a kernel that contains $N$. This induces a group homomorphism $\iota: \Ga \lrar \Ga_{\hp}\Big/ \ove{\iota_{\hp}(N)}  $. We are going to prove the required properties of \cref{uniprop} for the pair $(\Ga_{\hp}, \iota)$. The first two properties are immediate. We study the universal property of the third point.

Let $\KK$ be a pro-$p$ group and let $f: \Ga/N\lrar \KK$ be a group homomorphism with dense image. We have a diagram of the form 
\begin{equation}\label{propquo0}
    \begin{tikzcd}
        \Ga \ar[d, blue,  "q"] \ar[r, blue, "\iota_{\hp}"] & \Ga_{\hp} \ar[d, blue, "q_{\hp}"] \ar[ddr, bend left, dashrightarrow, "g"] & \\
         \Ga/N \ar[r, blue,  "\iota"] \ar[rrd, bend right, "f"] &  \Ga_{\hp}\Big/ \ove{\iota_{\hp}(N)} \ar[dr, dashrightarrow, "h"] & \\
          & & \KK,
    \end{tikzcd}
\end{equation}
where the subdiagram that is highlighted in blue is commutative. 

There exists a continuous group homomorphism $g: \Ga_{\hp}\lrar \KK$ such that $g\circ \iota_{\hp}=f\circ q$. Since $f$ has dense image and $q$ is surjective, the map $g$ is surjective. We also notice that $N\leq \Ga_{\hp}$ is contained in the kernel of $g$, so $\ove{\iota_{\hp}(N)}\unlhd \ker g$ is closed. By  \cref{first}, there exists a continuous group homomorphism $h: \Ga_{\hp}\Big/ \ove{\iota_{\hp}(N)} \lrar \KK$ such that $h\circ q_{\hp}=g$. 

The overall commutativity of the diagram (\ref{propquo0}) follows directly and the proof is complete. 
\end{proof}

It might be enlightening to remark that this provides an alternative proof of \cref{free pro-p completion}, the simple criterion for the pro-$p$ completion of a one-relator group to be free.

\begin{proof}[Second proof of \cref{free pro-p completion}] We identify $F(X)$ with its image in $\FF(X)$. The element $\omega\in F$ does not belong to $\FF(X)$ because the continuous map $f: \FF(X)\lrar \Z/p$ that determined by $f(x_k)=1$ and $f(x_i)=0$ for $i\neq k$ verifies that $f(\omega)\cong \omega_{x_k}\mod p$ so $\omega\notin \ker f$. In particular, by \cref{Fra2}, this implies that $\omega\notin \Phi(\GG)$. Furthermore, by \cref{freebasis}, there exists a topological generating set $\{\omega, \omega_2, \dots, \omega_n\}$ of $\FF(X)$. We consider $Y=\{y_1, \dots, y_n\}$ and a pro-$p$ completion map $\iota_{\hp}: F(Y)\lrar \FF(X)$ given by $\iota_{\hp}(y_1)=\omega$ and $\iota_{\hp}(y_i)=\omega_i$ for $i\geq 2.$ Then, by \cref{proppres}, there are canonical isomorphisms
\[\lan X, \omega \ran_{\hp}\cong \FF(X)/\ove{\lan\lan \omega\ran\ran} \cong \FF(X)/\ove{\iota_{\hp} (y_1)}\cong \left(F(Y)/\lan \lan y_1\ran\ran\right)_{\hp}\cong (F_{n-1})_{\hp}\cong \FF_{n-1}.\qedhere \]
\end{proof}

A sufficient understanding of the pro-$p$ completion of one-relator groups gave us a criterion for ensuring this to be free. We want to extend this principle for more general classes of groups.

The next step is to understand pro-$p$ completions of free products. This  will exhibit the functorial behaviour of the assignment $G\mapsto G_{\hp}$, from the category of abstract groups to the category of pro-$p$ groups.

\begin{prop}[{\bf Pro-$p$ completion of a free product}] \label{propprod} Let $\Ga$ and $\La$ be two abstract groups. There is a canonical isomorphism 
\[\left(\Ga*\La\right)_{\hp}\cong \Ga_{\hp}\coprod \La_{\hp}.\]
\end{prop}
\begin{proof} 
We define 
$j: \Ga* \La\lrar  \Ga_{\hp}\coprod \La_{\hp}$ by $j=\iota\circ (\iota_{\hp}*\iota_{\hp})$, where  $\iota_{\hp}*\iota_{\hp}: \Ga*\La\lrar \Ga_{\hp}* \La_{\hp}$ is canonical and $\iota: \Ga_{\hp}*\La_{\hp}\lrar \Ga_{\hp}\coprod \La_{\hp} $ is as in the proof of \cref{coprod}. 

To prove this proposition, it suffices to ensure the conditions of \cref{uniprop} for the pair $(\Ga_{\hp}\coprod \La_{\hp}, j)$. 

On the one hand, we claim that  the map $j$ has dense image. To verify this, we will need to refer to the proof of \cref{coprod}. Here, the free product $\Ga_{\hp}*\La_{\hp}$ was endowed with a topology under which the map $\iota: \Ga_{\hp}*\La_{\hp}\lrar \Ga_{\hp}\coprod \La_{\hp} $ is continuous and has dense image. It remains to proof that the map 
\[\iota_{\hp}*\iota_{\hp}: \Ga*\La\lrar \Ga_{\hp}* \La_{\hp}\]
has dense image. For doing so, we use  the basis of neighbourhoods of the identity $\mathcal{N}_0$  that defines the topology of $\Ga_{\hp}* \La_{\hp}$. Since these open subsets are normal subgroups, it suffices to check that, for every $N\in \mathcal{N_0}$, the induced map 
\begin{equation}\label{propprod1} \al_N: \Ga*\La\lrar \left(\Ga_{\hp}* \La_{\hp}\right)/N\end{equation}
is surjective. Every $N\in \mathcal{N_0}$ verifies, by definition, that $N\bigcap \Ga_{\hp}$ (resp. $N\bigcap \La_{\hp}$) is an open subgroup of $\Ga_{\hp}$ (resp. $\La_{\hp})$). So the canonical maps 
\[\Ga\lrar \Ga_{\hp}/N\bigcap \Ga_{\hp}\, \cong \, \Ga_{\hp}N/N\, ; \, \, \, \, \, \, \La\lrar \La_{\hp}/N\bigcap \La_{\hp}\, \cong \, \La_{\hp}N/N\]
are surjective. Hence the image of $\al_N$ contains both $\Ga_{\hp} N/N $ and $\La_{\hp} N/N$. Since these subgroups generate the whole $\left(\Ga_{\hp}* \La_{\hp}\right)/N$, then the map $\al_N$ of (\ref{propprod1}) is surjective.

 It rests to verify the  third condition of \cref{uniprop}, which is a universal property. Let $\KK$ be a pro-$p$ groups and let $f: \Ga*\La\lrar \KK $ be a group homomorphism. We have a diagram of the form
\begin{equation} \label{propprod2}
    \begin{tikzcd}
        \Ga \ar[rrdddd, bend right=110,  "f|_{\Ga}"'] \ar[rrdd, blue, "j|_{\Ga}"] \ar[rr, blue, hookrightarrow]\ar[dd, blue, "\iota_{\hp}"] & &\Ga*\La \ar[dd, blue, "j"]& &  \ar[dd, blue, "\iota_{\hp}"]\ar[ll, blue,  hookrightarrow]\ar[lldd, blue, "j|_{\La}"'] \ar[lldddd, bend left=110,  "f|_{\La}"]  \La\\
        & & & & \\
        \Ga_{\hp} \ar[rr, blue,  hookrightarrow] \ar[rrdd,  dashrightarrow, "f_{\hp}"] & & \Ga_{\hp}\coprod \La_{\hp} \ar[dd, dashrightarrow, "g"] & & \La_{\hp} \ar[ll, blue,  hookrightarrow]  \ar[lldd,  dashrightarrow, "f_{\hp}"']\\
        & & & &\\
        & & \KK & & 
    \end{tikzcd}
\end{equation}
The subdiagram that is highlighted in blue is commutative. By the universal properties of the pro-$p$ completions $\Ga_{\hp}$ and $\La_{\hp}$, there exist continuous homomorphisms $f_{\hp}^{\Ga}: \Ga_{\hp}\lrar \KK$ and $f_{\hp}^{\Ga}: \La_{\hp}\lrar \KK$ such that $f|_{\Ga}=f_{\hp}^{\Ga}\circ \iota_{\hp}$ and $f|_{\La}=f_{\hp}^{\La}\circ \iota_{\hp}$. 

By the universal property of the coproduct, there exists a continuous homomorphism $g: \Ga_{\hp}\coprod \La_{\hp}\lrar \KK $ that extends the previous $f_{\hp}^{\Ga}$ and $f_{\hp}^{\La}$.

Finally, by the universal property of the coproduct (now in the category of abstract groups), we can uniquely extend $f|_{\Ga}$ and $f|_{\La}$ to a group homomorphism  $f: \Ga*\La\lrar \KK$. 

The resulting diagram (\ref{propprod2}) is commutative. We are interested in verifying that $f=g\circ j$. The groups $\Ga$ and $\La$ generate $\Ga*\La$, so it suffices to prove that $f|_{\Ga}=g\circ j|_{\Ga}$ and  $f|_{\La}=g\circ j|_{\La}$. These identities can be read in the diagram (\ref{propprod2}). 
\end{proof}

We can now discuss two applications of the two previous propositions \ref{propquo} and \ref{propprod}.  These will be used during \cref{pro-$p$ embeddings} and they describe completions of amalgamated products and HNN extensions using the constructions that have just been introduced. 

\begin{cor}[{\bf Pro-$p$ completions of amalgamated products}] \label{propamal} Let $\theta_1: A\lrar \Ga$ and $\theta_2: A\lrar \La$ be group homomorphisms. Consider the corresponding amalgamated product $\Ga\underset{A}{*}\La$ of abstract groups. Then there is a canonical isomorphism of pro-$p$ groups 
\[ \left(\Ga\underset{A}{*}\La\right)_{\hp}\, \, {\bf \cong} \, \, \Ga_{\hp}\coprod \La_{\hp}\, \Bigg/ \ove{\lan \lan \iota(\theta_1(a)\, \theta_2^{-1}(a)), \, a\in A\ran\ran}. \]
\end{cor}
\begin{proof} Since \[\Ga\underset{A}{*}\La\cong \Ga*\La \Big/  \lan \lan \iota(\theta_1(a)\, \theta_2^{-1}(a)), \, a\in A\ran\ran,\] the claimed isomorphism comes from the canonical isomorphisms described in \cref{propquo} and \cref{propprod}.
\end{proof}

\begin{cor}[{\bf Pro-$p$ completions of HNN extensions}]\label{prophnn}  Let $A$ be a subgroup of $\Ga$ and let $\theta: A\rar \Ga$ be a group monomorphism. Consider the corresponding HNN extension $\Ga \hn{A}$. Then there is a canonical isomorphism of pro-$p$ groups
\[ \left(\Ga\hn{A}\right)_{\hp}\, \, {\bf \cong} \, \, \Ga_{\hp}\coprod \Z_p\Bigg/ \ove{\lan \lan \iota(t\, a\, t^{-1} \theta^{-1}(a)), \, a\in A\ran\ran},   \]
where $t$ also corresponds to the topological generator of $\Z_p$. 
\end{cor}

We finish this section with a non-trivial result about free pro-$p$ groups. It is obvious, by construction, that these are residually-$p$. However, we are going to verify that they are residually -(torsion-free nilpotent).

\begin{prop}\label{restorfree} A free pro-$p$ group is residually-(torsion-free nilpotent). 
\end{prop}

\begin{proof} A free pro-$p$ group is residually-(free of finite rank), so it suffices to restrict to the case of finite rank. Let $\FF=\FF_n$ be the free pro-$p$ group of rank $n$. 
Since $\FF$ is residually-$p$,  $\bigcap_k \ga_{k}^{(p)}\FF=\{1\}$ and, in particular, $\bigcap_k \ga_{k}\FF =1$. 

It rests to prove that each $\FF/\ga_k \FF$ is torsion-free. 

Now let $\Ga$ be a  finitely generated torsion-free nilpotent group with lower central series $1=\ga_{c+1}\Ga\unlhd \ga_c\Ga\no \cdots \ga_2\Ga\unlhd \ga_1\Ga=\Ga.$ Then $\ga_{k}\Ga/\ga_{k+1}\Ga\cong \Z^{n_k}$ for all $1\leq k\leq c$ and some non-negative integers $n_1, \dots, n_c$. By the indications of \cite[Exercises 21 and 22]{Dix99}, it follows that the canonical maps $(\ga_k \Ga)_{\hp}\lrar \ga_k \Ga_{\hp}$ and $(\ga_{k+1}\Ga/\ga_k \Ga)_{\hp}\lrar \ga_{k+1}\Ga_{\hp}/\ga_k\Ga_{\hp}$ are isomorphisms of pro-$p$ groups for all $1\leq k\leq c$. So there is a lower central series of $\Ga_{\hp}$ by closed subgroups
\[1\unlhd \ga_c\Ga_{\hp}\unlhd \cdots \unlhd \ga_2\Ga_{\hp}\unlhd \Ga_{\hp},\]
for which each quotient $\ga_{k+1}\Ga_{\hp}/\ga_k \Ga_{\hp}\cong \Z_p^{n_k}$ is torsion-free. So $\Ga_{\hp}$ is torsion-free. 

Now let $F=F_n$ be the free group of rank $n$. We know that $\Ga=F/\ga_{c+1} F$ is  a finitely generated torsion-free nilpotent group (this was checked during the proof of \cref{freeres1}). Hence by the previous remark, $\Ga_{\hp}$ is torsion-free. In addition, the canonical map $\Ga_{\hp}\lrar \FF/\ga_k(\FF)$ is an isomorphism by \cref{propquo}.
\end{proof}

\section{The completed group algebra} \label{completedsection}

Given an abstract group $G$, we denote by $\F_p[G]$, or simply by $\F_p G$, the group algebra of $G$ with coefficients in the field $\F_p$ of $p$ elements. In the case of pro-$p$ groups $\GG$, there is an analogous ring which is more interesting than its group algebra $\F_p \GG$, namely the completed group algebra, since it captures the topology of $\GG$ as well. 

\begin{defi}Let $\GG$ be a pro-$p$ group. Its {\bf completed group algebra} $\F_p [[\GG]]$ is defined by the inverse limit of the $\F_p$-group algebras $\F_p[ \GG/\NN]$, where $\GG$ ranges over open normal subgroup of $\GG$. We can write 
\[\p [[\GG]]= \varprojlim_{\NN\no \GG} \p[\GG/\NN].\]
\end{defi}

If $\GG$ is a free pro-$p$ group, there is a more explicit description of the completed group algebra. 

\begin{thm} \label{completediso}
Let $\FF$ be a free pro-$p$ group freely generated by $f_1, \dots, f_d$, then the continuous homomorphism $\F_p\lan \lan x_1, \dots, x_n\ran\ran \lrar\F_p[[\FF]] $ that sends $x_i$ to $f_i-1$ is an isomorphism.
\end{thm}
For a proof of the previous result, the reader is referred to \cite[Theorem 7.3.3]{Wil98}.

\begin{remark}[$\pF$ as a local ring] \label{completedinv} Let $P$ be a finite $p$-group. Then $\p[P]$ is a local ring with maximal ideal $\mm =\p I_{P}$ and residual division ring $\p [P]/\mm\cong \p $. Observe that given an inverse limit of rings $R=\varprojlim R_i$, its set of units  $R^{*}$ can be identified with  $\varprojlim R_i^*$. Given a pro-$p$ group $\GG$, we can use the latter observation to prove that $\p[[\GG]]$ is also a local ring. The completed group algebra $\F_p[[\GG]]$ also has an augmentation map $\epsilon: \F_p[[\GG]]\lrar \F_p$, whose kernel $\mm$ is the only maximal ideal and equals the closure of $\p \, I_{\F}\sub \p[\FF]$ in $\p[[\FF]]$ (endowed with the profinite topology). In case $\GG$ is a free pro-$p$ group $\FF$, one can alternatively check that $\pF$ is local as follows as a corollary of the isomorphism of \cref{completediso}. In any case, we conclude that an element $x\in \p[\FF]$ is invertible in $\F_p [[\FF]]$ if and only if $x\notin \p I_{\FF}$. 
\end{remark}

\section{Cohomological characterisation of freeness} \label{propcohomology}

 Free abstract groups and free pro-$p$ groups share some similar properties. They are both, in their corresponding categories of groups, the free objects. From the point of view of this section, these are the groups of cohomological dimension at most one. In the category of abstract groups, this characterisation is due to deep arguments given by Stallings and Swan. 

We are going to review the proof of this characterisation in the pro-$p$ setting, assuming some bits of the cohomology theory of profinite groups. We will assume, for simplicity, finite topological generation  (finite rank). 


\begin{defi} A  pro-$p$ group $\GG$ is said to have {\bf $p$-cohomological dimension $d$} if $n=d$ is the smallest $n$ such that $H^m(\GG, \p)=0$ for every $m\geq n+1$, where $\p$ is considered with the discrete topology.  We denote  $\cdp(\GG)=d$.  
\end{defi}

The following is \cite[Theorem 7.3.1]{Rib00}.

\begin{thm}\label{cdH} Let $\GG$ be a pro-$p$ group and let $\HH$ be a closed subgroup. Then $\cdp (\HH)\leq \cdp (\GG)$.
\end{thm}

\begin{thm} \label{preSwan} Let $\GG$ and $\GG'$ be finite-rank pro-$p$ groups. Let $f: \GG\rar \GG'$ be a continuous homomorphism such that \[f^*: H^1(\GG', \F_p)\rar H^1(\GG, \F_p)\] is an isomorphism and \[f^*: H^2(\GG', \F_p)\rar H^2(\GG, \F_p)\] is an injection. Then $f$ is an isomorphism.
\end{thm}

The idea of the proof of the previous result is that a pro-$p$ groups $\GG$ is an inverse limits of $p$-groups that are related to $H^1(\GG, \F_p)$ and $H^2(\GG, \F_p)$ by short exact sequences. 
We should remark that the root of this principle is specially  the paper of Stallings \cite{Sta65}. We are giving an adaptation to the setting of pro-$p$ groups, taken from \cite{Wil21}, avoiding the technicalities of the cohomology theory of profinite groups.

This result has several important consequences, which we shall state later. Briefly, it says that abstract groups with trivial $H^2(\, -, \F_p)$ have free pro-$p$ completion. On other side, this characterisation of freeness is easily seen to be stable with respect to subgroups. In particular,  it allows to  derive the non-trivial fact that closed subgroups of free pro-$p$ groups are pro-$p$ free.

Before giving the proof, we state the Five Term Exact sequence lemma. Given an open subgroup $\HH$ of a pro-$p$ group $\GG$, this lemma relates the cohomology (with coefficients in $\F_p$) of the three groups $\HH, \FF, \FF/\HH$. 

\begin{lemma}[Five Term Exact sequence] Given an open subgroup $\HH$ of a pro-$p$ group $\GG$ and a $\GG$-module $M$, we have a natural short exact sequence of the form
\[0\rar H^1(\GG/\HH, M^{\HH})\rar H^1(\GG, M)\rar H^1(\HH, M)^{\GG}\rar H^2(\GG/\HH, M^{\HH})\rar H^2(\GG, M).\]
\end{lemma}

\begin{remark} \label{trivialaction} If $M$ is a trivial $\GG$-module, then $H^1(\GG, M)=\Hom (\GG, M)$. In this case, we can also describe more precisely $H^1(\HH, M)^{\GG}$, for an open normal subgroup $\HH$ of $\GG$. The underlying action of $\GG$ on $H^1(\HH, M)$ is the one induced by the conjugation action of $\FF$ on $\HH$. So 
\begin{align*} H^1(\HH, M)^{\GG} =& \{\phi\in \Hom (\HH, M): \phi (ghg^{-1})=\phi (h)\, \mbox{for all $g\in \GG, h\in \HH$}\}\\
=& \{\phi\in \Hom (\HH, M): \phi ([g, h])=0\, \mbox{for all $g\in \GG, h\in \HH$}\}\\
=& \Hom (\HH/\overline{[\GG, \HH]}, M).
\end{align*}
\end{remark}
\begin{proof} 
For each $n$, recall the notation $\GG_{(n)}=\gamma_n^{(p)}\GG$ and $\GG'_{(n)}= \gamma_n^{(p)}\GG'$. Since $\GG\cong  \varprojlim_{n} \GG/ \GG_{(n)} $, it suffices to check that, for each $n$, the induced map 
\begin{equation}\label{f_n} f_n : \GG/ \GG_{(n)}\rar \GG'/\GG'_{(n)}\end{equation}
is an isomorphism. We proceed by induction on $n$. The case $n=1$ is trivial. Let us suppose that $f_n$ is an isomorphism for some $n\geq 1$. 
\begin{equation*}
\begin{CD}
 H^1 (\GG'/\GG_{(n)}') @>>> H^1(\GG')   @>>> H^1(\GG_{(n)}')^{\GG'} @>>>  H^2(\GG'/\GG_{(n)}') @>>> H^2(\GG')\\
@V{g_1}VV @V{g_2}VV @V{g_3}VV @V{g_4}VV @V{g_5}VV \\
H^1 (\GG/\GG_{(n)}) @>>> H^1(\GG)   @>>> H^1(\GG_{(n)})^{\GG} @>>>  H^2(\GG/\GG_{(n)}) @>>> H^2(\GG)\\
\end{CD}
\end{equation*}
By assumption, $g_2$ is an isomorphism and $g_5$ is an injection, and by the induction hypothesis, $g_1$ and $g_4$ are isomorphisms. So applying the five lemma to the previous diagram we deduce that $g_3$ is an isomorphism\footnotemark. \footnotetext{Recall that, in order to apply the five lemma, we only needed $g_1$ to be surjective.}  By \cref{trivialaction} \[H^1(\GG_{(n)})=\Hom (\GG_{(n)}/[\GG, \GG_{(n)}], \F_p),\] and since $\F_p$ has exponent $p$, the previous equals  \[\Hom \left(\GG_{(n)}\Big/\GG_{(n)}^p\, [\GG, \GG_{(n)}], \F_p\right),\] which is the dual of $\GG_{(n)}/\GG_{(n+1)}$ as $\F_p$-vector space. So the dual map $h_n$ of $g_3$ is the map
\[\GG_{(n)}/\GG_{(n+1)}\xrightarrow{h_n} \GG'_{(n)}/\GG'_{(n+1)}\] induced by $f$ and this is an isomorphism.  
Lastly, we consider the commuting diagram:
\begin{equation*}
\begin{CD}
1 @>>> \GG_{(n)}/\GG_{(n+1)}   @>>> \GG/\GG_{(n+1)} @>>>  \GG/\GG_{(n)} @>>> 1\\
@VVV @V{h_n}VV @V{f_{n+1}}VV @V{f_n}VV @VVV \\
1 @>>> \GG'_{(n)}/\GG'_{(n+1)}   @>>> \GG'/\GG'_{(n+1)} @>>>  \GG'/\GG'_{(n)} @>>> 1\\
\end{CD}
\end{equation*}
Since $h_n$ and $f_n$ are isomorphisms, another application of the five lemma yields that $f_{n+1}$ is an isomorphism. This completes the induction. 
\end{proof}

\begin{cor} \label{cd1}
A pro-$p$ group is free if and only if its cohomological dimension is at most  1. 
\end{cor}
\begin{proof}
We only prove the hard direction. It is easy to check that free pro-$p$ groups $\FF$  have cohomological dimension one by using the interpretation of $H^2(\FF, \p)$ in terms of the splitting of short exact sequences $0\lrar \p \lrar \GG\lrar \FF \lrar 1$.  For a complete proof, we refer the reader to \cite[Theorem 7.7.4]{Rib00}. Let $\GG$ be a pro-$p$ group of cohomological dimension 1. Let $\FF$ be a free pro-$p$ group such that $d(\GG)=d(\FF)$. In particular, we have 
\[H^1(\GG, \F_p)\cong \GG/\Phi (\GG)\cong \FF/\Phi (\FF)\cong H^1(\FF, \F_p),\] so any surjective map $f: \FF\rar \GG$ induces an surjection between finite groups $f^* : H^1(\GG, \F_p)\rar H^1(\FF, \F_p)$. Hence the latter $f^*$ is an isomorphism and the corresponding $f^* : H^2(\GG, \F_p)\rar H^2(\FF, \F_p)$ would be an injection because $H^2(\GG, \F_p)=0$. Thus any such $f$, which can be constructed because $\FF$ is free, would be an isomorphism by \cref{preSwan}.  
\end{proof}

Essentially the same proof of \cref{preSwan} gives a way of ensuring that an abstract group with $H^2(G, \p)=0$ has free pro-$p$ completion. This observation, in combination with Hopf's formula on $H^2$, can give another proof of \cref{free pro-p completion}. A more general consequence of the argument of \cref{preSwan} is the following.

\begin{cor} \label{postSwan} Let $G$ and $G'$ be finitely generated (abstract) groups and let $f: G\rar G'$ be a group homomorphism such that \[f^*: H^1(G', \F_p)\rar H^1(G, \F_p)\] is an isomorphism and
\[f^*: H^2(G', \F_p)\rar H^2(G, \F_p)\] is an injection. Then the induced map  $f_{\hp}: G_{\hat{p}}\lrar G'_{\hat{p}}$ is an isomorphism.
\end{cor}

We also leave stated an important result which we shall use in the arguments of \cref{pro-$p$ embeddings}.

\begin{cor}\label{subfree} A (finite rank) closed subgroup $\HH$ of a (finite rank) free pro-$p$ group $\FF$ is free. 
\end{cor} 
\begin{proof} By \cref{cdH}, $\cdp(\HH)\leq \cdp(\FF)\leq 1$. By \cref{cd1}, we conclude that $\HH $ is free. 
\end{proof}

\chapter{Parafree groups} \label{parafreesection}

\begin{defi}
We say that two groups $G_1$ and $G_2$ have the same {\bf lower central quotients} if, for every $k\geq 1$, we have 
\[\frac{G_1}{\gamma_k G_1} \cong \frac{G_2}{\gamma_k G_2}.\]
\end{defi}

\begin{defi}
A finitely generated abstract group $G$ is termed {\bf parafree} if it has the following two properties: 
\vspace{-1cm}
\begin{item}
\item[(i)] It has the same lower central sequence as some free group $F$.
\item[(ii)] It is residually nilpotent. 
\end{item}
\end{defi}


\section{Classical examples and counterexamples}

In this section we will collect some known families of parafree groups. We first remark how they relate to other groups that we have already mentioned. First, notice that Baumslag-Solitar groups are not residually nilpotent. In fact, by \cref{BShopf}, $B(n, m)$ is not even residually finite if $n, m\geq 2$. If $n=1$ and $m=2$, we have $\Ga=B(1, 2)=\lan a, b|bab^{-1}=a^2\ran$. Observe that $a=[a, b^{-1}]\in \ga_2 \Ga$. Recursively, we deduce that $a$, which is a non-trivial element of $\Ga$, belongs to $\bigcap \ga_k \Ga$. So $B(1, 2)$ is not residually nilpotent. 
It is clear that the remaining group $B(1, 1)=\Z^2$ is not parafree.  
A different way to verify Baumslag-Solitar groups are not parafree goes as follows. It is easy to see that exactly when $n, m$ are coprime, $B(n, m)$ has abeliank rank equal to 1, though they are not isomorphic to $\Z$ (as they should, by \cref{Z}). In other cases, it is not even true that the abelianisation of $B(n, m)$ is torsion-free.  

Another counterexample to parafree groups are non-free surface groups. For example, the fundamental group of non-orientable closed surfaces $\pi_1(N_g)$ do not have torsion-free abelianisation. On the other side, the groups $\pi_1(\Sigma_g)$ have abelianisation isomorphic to $\Z^{2g}$. If these were parafree, then they would be non-free parafree groups $\Ga$ with $d(\Ga)=d(\Ga_{\ab})$. We will see in \cref{parastrongHopf} that this cannot happen.

We now turn to review some examples. Recall that our aim is to study parafree groups and that these groups are characterised by two properties. The first property, having the same nilpotent genus as a free group, was studied and characterised in different ways during the previous section. The second property enjoyed by groups that are named parafree is their residual nilpotence. This condition is considerably harder to verify. 

Bausmlag also produced many examples parafree one-relator groups and as amalgamated products of a parafree group and $\Z$ in \cite{Bau67} and \cite{Bau69}.

We can described some of these families, which are also surveyed in \cite{Bau05}.

\begin{defi} \label{parafamilies}
We introduce the following families of groups. 
\begin{align*}
    G_{i,j} &=\lan a, b, c\, |\,  a=[c^i, a][c^j, b] \ran\,\, \,  \, \mbox{for any positive integers $i, j$.} \\
     H_{i,j} &=\lan a, b, c\, |\,  a=[a^i, t^j][s, t] \ran\, \, \, \, \mbox{for any positive integers $i, j$.} \\
      K_{i,j} &=\lan a, b, c\, |\,  a^i[s, a]=t^j \ran\, \, \, \, \mbox{for coprime integers $i, j$.} \\
      N_{p, q, r} &=\lan a, b, c\, | \, a^p b^qc^r\ran\, \, \, \, \mbox{for non-zero integers with $\gcd (p, q, r)=1$.} 
\end{align*}
\end{defi}
These families of groups split over $\Z$. The  first  are cyclic HNN extensions of $F_2$. On the other side, the two last families are amalgamated products of $F_2$ and $\Z$ with cyclic amalgams.

For the question of whether the groups of \cref{parafamilies} are pairwise non-isomorphic or not, see \cite[Section 9]{Bau05} and the references therein.

 The following source of examples is due to Baumslag and Cleary. We first need to introduce a definition to describe the structure of the defining relations of this family.

\begin{defi} Given free generators $s_1, \dots, s_n, t$ of $F_{n+1}$ and $w\in [F_{n+1}, F_{n+1}]$, we set $s_{i, j}=t^{j} s_i t^{-j}$ for each $1\leq i\leq n$ and $j\in \Z$. We can express $w$ uniquely  as a word on the previous $s_{i, j}$. Given $i$, we define $\mu (i)$ (resp. $\nu (i)$) to be the minimum (resp. maximum) of those $j$ such that $s_{i, j}$ appears in $w$. We say that $w$ satisfies \textit{the redundancy condition} on $s_{i'}$ if $\mu (i')$ and $\nu(i')$ are distinct and both $s_{i', \mu(i')}$ and $s_{i', \nu(i')}$ appear once and only once in $w$ (possibly as inverses).
\end{defi}

We can now formulate the following result from \cite[Theorem 3]{Bau06}. 
\begin{cor} \label{parafamilies2} Let $F$ be the free group on $a_1, \dots, a_p, s_1, \dots, s_n, t$, where $p\geq 1$ and $n\geq 1$. Furthermore, let $E$ be the subgroup of $F$ generated by $s_1, \dots, s_n, t$ and let $w$ be a cyclically reduced word in $[E, E]$ that satisfies the redundancy condition on $s_{i'}$. Finally, let $v$ belong to $[F, F]$ and suppose that $v$ does not involve $s_{i'}$. Then the one-relator group 
 \[G =\lan a_1, \dots, a_p, s_1, \dots, s_n, t\, | \, a_1=v w \ran\]
is parafree.  
\end{cor}



\section{Free nilpotent genus}

\begin{defi}
We say that two groups have the same {\bf nilpotent genus} if they have the same isomorphism types of nilpotent quotients. Similarly, given a prime $p$, we say that two groups have the same {\bf $p$-genus} if they have the same isomorphism types of finite $p$-groups quotients. We say that a group has {\bf free nilpotent genus} if it has the same the same nilpotent genus as some free group.
\end{defi}

 In these terms, studying parafree groups is studying residually nilpotent groups with free nilpotent genus. 
Let us look at a few alternative ways to look at this family of groups.

\begin{prop}[Characterisations of the property of having free nilpotent genus] \label{characterisations} Let $G$ be a finitely generated group. The following conditions are equivalent. 
\begin{itemize}
    \item[(i)] There exists a free group $F$ that has the same nilpotent genus as $G$.
    \item[(i')] There exists a free group $F$ such that, for every $k\geq 1$, $ F/\gamma_kF\cong G/\gamma_kG$.
    \item[(i'')] There exists a free group $F$ and an injection $\phi: F\lrar G$ such that, for every $k\geq 1$, the induced map $\phi_k : F/\gamma_kF\lrar G/\gamma_kG$ is an isomorphism. 
    \item[(ii)] There exists a free group $F$ with the same $p$-genus as $G$ for every prime $p$.
    \item[(ii')] There exists a free group $F$ such that, for every $p$ and every $k\geq 1$, $ F/\gamma_{k, p}F\cong G/\gamma_{k,p} G$.
    \item[(ii'')] There exists a free group $F$ and an injection $\phi: F\lrar G$ such that, for every prime $p$, the induced maps $\phi_{k,p} : F/\gamma_{k,p} F\lrar G/\gamma_{k, p}G$ are isomorphism. 
    \item[(iii)] There exists a free group $F$ such that, for every prime $p$, $G_{\hat{p}}\cong F_{\hp}$. 
    \item[(iii')] There exists a free group $F$ and an injection $\phi: F\lrar G$ such that, for every prime $p$, the induced map $\phi_{\hat{p}}: F_{\hat{p}}\lrar G_{\hat{p}}$ is an isomorphism.
\end{itemize}
\end{prop}
These equivalences are well-known. We include a proof for the convenience of the reader.

\begin{proof} We discuss each implication separately. We simply observe in advance that, under all circumstances, it is clear that the free group $F$ must have finite rank $n=d(F)$.  
\begin{enumerate}
    \item $(i)\implies (i')$. Let $k\geq 1$. The group $G/\ga_k G$ is a nilpotent quotient of $G$. So there is an epimorphism $f: F\lrar G/\ga_k G$, which induces an epimorphism $f_k: F/\ga_k F\lrar G/\ga_k G$. Similarly, there is an epimorphism $g_k: G/\ga_k G\rar F/\ga_k F$. The group $F/\ga_k F$ is finitely generated nilpotent. By the Hopf property, this implies that $f_k$ and $g_k$ are inverses of each other. 
    \item $(i')\implies (i)$. This is trivial because the collection of isomorphism types of nilpotent quotients of a group $\Ga$ is the union of collections of isomorphisms types of quotients of each $\Ga/\Ga_k$. 
    \item $(i'')\implies (i')$ is trivial.
    \item $(i')\implies (i'')$.  Take  $g_1, \dots, g_n\in G$  lifts of  generators of $G/\gamma_2 G$. We name by $f_1, \dots, f_n$ free generators of $F$ and consider the homomorphism $\phi: F\rar G$ such that $\phi (f_i)=g_i$. We will prove that the induced maps $\phi_k : F/\gamma_k F\rar G/\gamma_k G$ are isomorphisms. By construction, $\phi_2$ is an isomorphism. In particular, $\phi (F)\, [G, G]=G$. This implies that, naming $G_{(k)}=G/\gamma_k G$,  $\phi_k (F_{(k))}\, [G_{(k)}, G_{(k)}]=G_{(k)}$. Since each $G_{(k)}$ is nilpotent, we deduce from  \cref{Basisnilp} that $\phi_k(F_{(k)})=G_{(k)}$; hence  $\phi_k$ is surjective. We also know that each $F_{(k)}$ is Hopfian (\cref{Hopfnilp}). Since $F_{(k)}\cong G_{(k)}$, the latter implies that each $\phi_k$ is injective (see \cref{Hopf1}). So every $\phi_k$ is an isomorphism. Lastly, we can claim that $\phi$ is injective because it factors through every isomorphism $\phi_k$;  and $F$ is residually nilpotent, which means that $\bigcap \ga_k F=1$.  
    \item The equivalence $(ii)\iff (ii')\iff (ii'')$ is entirely analogous to the already proven equivalence $(i)\iff (i')\iff (i'')$. 
    \item $(i)\implies (ii)$ is trivial. 
    \item $(ii'')\implies (i'')$. This is simply due to the fact that $\bigcap_p \ga_{ k, p}F=\ga_k F.$
    \item $(iii)\implies (ii)$. This is a direct consequence of the correspondence (\ref{correspondence}) between the finite $p$-quotients of an abstract group $\Ga$ and those of its pro-$p$ completion $\Ga_{\hp}$. 
    \item $(ii')\implies (iii)$. Notice that this is a particular case of \cref{proC}. In this particular case, we give a more direct proof. The given $\phi$ induces isomorphisms 
    \[F_{\hp}\cong \varprojlim_k F/\gamma_{k, p}F\cong \varprojlim_k G/\gamma_{k,p} G\cong G_{\hp}. \]
    \item $(iii')\implies (iii)$ is trivial.
    \item $(iii)\implies (iii')$.  Since $G_{\hp}\cong F_{\hp}$, then $G_{\hp}$ is pro-$p$ free. In addition, since $(iii)\implies (i)$, $G$ and $F$ have the same abelian quotients. From this, we  deduce that $G_{\ab}\cong F_{\ab}$. Take $g_1, \dots, g_n$ such that their reductions modulo $[G, G]$ form a $\Z$-basis of $G_{\ab}\cong \Z^n$.
    Consider the map $\phi: F\lrar G$ defined by sending $\phi(f_i)=g_i$, for some generating set $\{f_1, \dots, f_n\}$ of $F$. For every prime $p$, there is a commutative diagram 
    \begin{equation*}
        \begin{tikzcd}
            F\ar[d, hookrightarrow, "\iota_{\hp}"] \ar[r, "\phi"] & G \ar[d]\\
            F_{\hp}\ar[r, "\phi_{\hp}"] & G_{\hp}.
            \end{tikzcd}
    \end{equation*}
    Notice that $\phi_{\hp}$ is surjective, since the induced $\phi: F/F^p[F, F]\lrar G/G^p[G, G]$ is surjective, by construction. 
    So $\phi_{\hp}$ is a continuous epimorphism between two free pro-$p$ groups of the same rank. By \cref{Hopf2}, $\phi_{\hp}$ is an isomorphism. On the other side, 
     $F$ is residually-$p$, so $\iota_{\hp}$ is injective, too. This implies that $\phi$ is injective.  
\end{enumerate}
With these equivalences, the proof is complete.
\end{proof}


From this point of view, we can introduce another criterion for a finitely generated $G$ to have free nilpotent genus; which is a direct consequence of \cref{postSwan} and \cref{characterisations}.

\begin{prop} \label{postSwan1} Let $G$ be a finitely generated group such that $G/[G, G]$ is free abelian and $H^2(G, \F_p)=0$ for every prime $p$. Then $G$ has the same lower central series as a free group of the same rank as $G/[G, G]$. 
\end{prop}

\section{Similarities and differences with free groups}

We know that finitely generated residually nilpotent groups are  Hopfian. In particular, parafree groups are Hopfian. Moreover, we have an strong Hopf property for parafree groups, which says that parafree proper quotients of parafree groups have strictly lower rank. 

\begin{prop} \label{parastrongHopf} Let $N$ be a normal subgroup of a group $P$. Suppose that both $P$ and $P/N$ are parafree groups of the same rank, then $N=1$. 
\end{prop}
\begin{proof} We consider the natural projection $p: P\rar P/N$ and we will show that its kernel $N$ is trivial. This projection induces a surjective map 
\[p_n: \,  \frac{P}{\gamma_n P} \longrightarrow \frac{P/N}{\gamma_n (P/N)}\cong \frac{P}{(\gamma_n P)N}, \]
for every $n$. Since both of the previous groups are isomorphic to $F/\gamma_n F$, for the same free group $F$, $p_n$ is a surjective endomorphism of the Hopfian group $F/\gamma_n F$, so it is injective. This implies that $N\subset \gamma_n P$ for every $n\geq 1$, so $N=1$. 
\end{proof}
This is due to the fact that finitely generated nilpotent groups are Hopfian. 
We have two immediate consequences from the last proposition.
\begin{cor} \label{Z} Parafree groups of abelian rank 1 are isomorphic to $\Z$.
\end{cor}
\begin{proof}
Applying \cref{parastrongHopf} to $N=[G, G]$ it follows that $G\cong G/N\cong \Z$. 
\end{proof}

By similar reasons, we can also prove that the only abelian (equivalently, nilpotent) parafree groups are the trivial group and $\Z$. 

\begin{prop} The free product of parafree groups is parafree. 
\end{prop}

\begin{prop}
Subgroups of parafree are not necessarily parafree. 
\end{prop}
\begin{proof}
The group $G=\lan a, b, c| a^2b^2c^3\ran$ is parafree. The subgroup of $G$ generated by $a, b, cac^{-1}, cbc^{-1}$ has presentation $\lan x, y, z, t| x^2y^2z^2t^2\ran\cong \pi_1(S_3)$. So $G$ contains the surface groups 
\[\pi_1(\Sigma_g)\inc \pi_1(S_3)\inc G\]
for $g\geq 2$; which are not parafree.
\end{proof}
The following proposition is harder. Its proof can be found in \cite{Bau69}.

\begin{prop} \label{twogen}
The two-generator subgroups of parafree groups are free. 
\end{prop}

\begin{prop} For every integer $n\geq 2$, there exists a parafree $G$ with $\abr(G)=n$ and with a non-free three-generated subgroup.
\end{prop}

\begin{proof} Let $n\geq 2$ and consider the 3-generated group $G_3=\lan a, b, c | a^2b^2c^3\ran$, which is a non-free parafree group. By taking $G_n=G_3* \Z *\cdots \Z$, the free product of $G_3$ with $n-2$ copies of $\Z$, we get a group $G_n$ with $\abr(G_n)=n$ that contains a subgroup isomorphic to $G_3$.
\end{proof}
 
 The previous group $G_3=\lan a, b, c | a^2b^2c^3\ran$ is quite interesting. It is not only a non-free parafree group, it is  residually free (\cite{Bau62} and non-limit (\cite{Lyn62}).

\begin{cor} The centre of a parafree group is either trivial or the whole group. This last case only occurs  for $\Z$. More generally, any abelian subgroup of a parafree group is free. 
\end{cor}

The following is \cite[Theorem 7.1]{Bri14}. It uses the $L^2$-methods of \cref{L2section}.
\begin{thm} Let $G$ be a finitely generated parafree group and let $N$ be a finitely generated normal subgroup of infinite index. Then $N=1$. 
\end{thm}
\begin{proof} The group $G$ has rank $n$. If $n=1$ then, by  \cref{Z}, $G\cong \Z$ and the statement is trivial. Suppose that $n\geq 2$ and let $F$ be the free group of rank $n$. Since $G$ is parafree, $G_{\hp}$ is pro-$p$ free by \cref{characterisations} and  the canonical map $G\inc G_{\hp}\cong F_{\hat{p}}$ is injective and has dense image. By \cref{Bettiest1}, $\b(G)\geq \b(F)=n-1>0$. Let us suppose that $N$ is not trivial. Then, since $G$ is torsion-free, $N$ is infinite. We have a short exact sequence
\[1\rar N\rar G\rar G/N\rar 1,\]
where both $N$ and $G/N$ are infinite, $\b(G)>0$ and $\b(N)<\infty$ (since $N$ is finitely generated). This contradicts \cref{shortBetti}. Thus $N=1$.
\end{proof}

\chapter{The theory of group rings} \label{groupringsection}
The main idea that underlies this chapter can be roughly outlined as follows. Given a group homomorphism $G\lrar H$, we want to ``linearise'' it by studying instead the corresponding map $I_G\lrar I_H$. In the last setting, cohomological and ring-theoretical arguments are applicable. From this, one wants to derive information about the initial group homomorphism. The most important result of this chapter is \cref{kernelmodp}.

\section{The bar resolution}
Let $G$ be a group. We review the construction of the bar resolution $\bC^*(G, k)\rar k\rar 0$, a resolution of $k$ over $kG$. It is also named \textit{standard resolution} and it allows us to explicitly compute homology groups of low degree.  Let $\bC^n(G, k)$ be a $kG$-module with free basis given by the symbols $\{|g_1|\cdots |g_n|: g_i\in G\}$, a collection that is naturally in bijection with $G^n$. Observe that $\bC^0(G, k)=kG$. We define $d_n: \bC^{n}(G, k)\rar \bC^{n-1}(G, k) $ by 
\begin{align*} d_n(|g_1|\cdots|g_n|) = & \, g_1\, |g_2|\cdots|g_n|\\
&-|g_1g_2|g_3|\cdots |g_n|\\
&+\dots \\
&+(-1)^{n-1}\, |g_1|\dots |g_{n-1}g_n|\\
&+(-1)^n\, |g_1|\cdots |g_{n-1}|.\end{align*}
We claim that the following is a resolution of projective (in fact, free) $kG$-modules
\[\xrightarrow{d_{n+2}}\bC^{n+1}(G, k)\xrightarrow{d_{n+1}} \bC^{n}(G, k)\xrightarrow{d_{n}} \cdots \xrightarrow{d_{3}} \bC^{2}(G, k)\xrightarrow{d_{2}} \bC^1(G, k)\xrightarrow{d_{1}} kG \xrightarrow{d_{0}} k\rar 0. \]
This particular resolution allows us to compute $H^n(G, k)$ as the homology of the chain complex obtained after applying the  functor $\tens{kG} k$ to the previous resolution. This sequence is natural from the point of view of topology. 

Homology groups also appear as correction terms in sequences that are not exact and homology groups of low degree are very meaningful. In fact, $H^0(G, k)\cong k$ and $H_1$ verifies the following. 

\begin{lemma} \label{H_1} There is a natural isomorphism $H_1(G, k)=k\tens{\Z} G_{\ab}$.
\end{lemma}

\begin{proof} We use the bar resolution of $k$ to compute $H_1(G, k)$. By applying to 
\[ \xrightarrow{d_{3}} \bC^{2}(G, k)\xrightarrow{d_{2}} \bC^1(G, k)\xrightarrow{d_{1}} kG \xrightarrow{d_{0}} k\rar 0\]
the functor $\tens{kG} k$, we obtain the sequence of $k$-modules
\[ \xrightarrow{{d}_{3}\otimes k} kG^2 \xrightarrow{d_2\otimes k} kG  \xrightarrow{0} k \xrightarrow{\id} k\rar 0,\]
where $\ti{d}_2$ sends each $(g, h)$ to $g+h-gh$. So 
\[H_1(G, k)\cong \frac{k\,  G^2}{\im \ti{d}_2}\cong k\otimes_{\Z} G_{\ab}.\qedhere\]
\end{proof}

More importantly, this particular resolution, the bar resolution, allows us to add explicitly an homology group to a non-exact resolutions of $k$ over $kG$, to obtain an exact sequence. Of course, here underlies the fact that the computation of the homology does not depend on the choice of projective resolution of $k$ over $kG$, since any pair of such resolutions is homotopically equivalent and the application of $\tens{kG}k$ preserves chain homotopies, as it is an additive functor. This fact will underlie the arguments of subsequent sections. 

\section{The augmentation ideal}\label{augsection}
During \cref{presentationsection}, we saw the importance of defining objects in terms of universal properties. We worked with pro-$p$ completions in terms of neat diagrams and, roughly, we could observe that it is functorial. The approach of this section is similar. We introduce the notion of derivation to describe the augmentation ideal $k I_G$ in terms of a universal property. 

\begin{defi}[Derivation] Let $M$ be a $G$-module. We say that a function $d: G\rar M$ is a derivation if $d(xy)=d(x)+xd(y)$ for every $x, y\in G$. We denote by $\Der (G, M)$ the abelian group of derivations $G\rar M$. 
\end{defi}

There is a {\bf universal derivation} $D$ of $G$. We define $D: G\rar kI_G$ to be the derivation such that $D(g)=g-1$. There is an isomorphism of abelian groups $\Der( G, M)\cong \Hom_G(kI_G, M)$ described as follows. Given $f\in \Hom_{kG}(kI_G, M)$, we define a derivation $d=f\circ D: G\rar M$. 
\begin{equation} \label{deriv}
\begin{tikzcd}
 G \ar[rd, dashrightarrow, "d"] \ar[r, "D"] & k I_G \ar[d, "f"] \\
  & M.
\end{tikzcd}
\end{equation}

Reciprocally, given $d\in \Der (G, M)$, we can define the associated $kG$-homomorphism $f: kI_G\rar M$ by $f(g-1)=d(g)$. These correspondences are well-defined and are inverses of each other.

\begin{lemma}  Let $F$ be a free group freely generated by  $S\sub F$. Then $\{s-1: s\in S\}$ is a free basis of the $kF$-module $kI_F$. 
\end{lemma}
\begin{proof} It is clear that $\Der (G, A)$ is in bijection with group homomorphisms $\ga : F\rar M\rtimes F$ that are splittings of the exact sequence 
\[1\rar M\rar M\rtimes F\rar F\rar 1.\]
In fact, those $\ga$ are exactly those maps of the form $\ga(f)=(d(f), f)$ for some derivation $d$. In other words, for any $\{a_s: s\in S\}\sub M$, there is exactly one $d\in \Der(F, M)$ such that $d(s)=a_s$ for all $s\in S$. Using the characterisation $\Der (F, M)\cong \Hom_{kF} (kI_F, M)$ of \cref{deriv}, we can also conclude that for any $kF$-module $M$ and any 
$\{a_s: s\in S\}\sub M$, there is exactly one $kF$-homomorphism $f: I_F\rar M$ such that $f(s-1)=a_s$ for all $s\in S$. The conclusion follows. 
\end{proof}

\begin{defi}[Fox derivatives]  Let $F$ be a free group freely generated by $S\sub F$.  We denote by $\pd{}{s}: F\rar k F$ the maps defined by 
\[f-1=\sum_{s\in S} \pd{f}{s}\, (s-1), \,  \, \, \mbox{for all $f\in F$.}\]
\end{defi}

Our aim is to prove the following result. 

\begin{thm} \label{Relmod} Let $F$ be a free group freely generated by $S\sub F$. 
Let $T\sub F$ and consider the normal closure $N=\lan \lan T\ran \ran$. Consider $G=F/N$. 
Given $f\in F$, we denote its image under the projection $F\rar G$ by $\overline{f}\in G$. There is an exact sequence of $kG$-modules of the form 
\begin{equation*}
\begin{tikzcd}
0\ar[r] & k\tens{\Z}R_{\ab}  \ar[r, "\partial_2"] & kG^{(S)} \ar[r, "\partial_1"] & kG \ar[r, "\epsilon"] & k\ar[r] & 0,
\end{tikzcd}
\end{equation*}
where the maps $\delta_1$ and $\delta_2$ are defined as follows:
\begin{align*} 
\partial_1\left(g_s\, e_s\right) =& g_s \, (s-1), \, \, \, \mbox{for all $s\in S$ and $g_s\in G$,}\\
\partial_2\left(c\otimes \overline{r}\, [G, G]\right) =& \sum_{s\in S} c\,  \overline{\pd{r}{s}}\, e_s,\, \, \, \mbox{for all $c\in k$, $r\in R$.}
\end{align*}
\end{thm}

\begin{defi} In the context of the previous theorem, the group $R_{\ab}$, or more broadly $k\otimes_{\Z} R_{\ab}$, is called the \textit{relation module} of $G$. Notice that it requires a choice of $F$ and a map $F\rar G$. It has a natural $kG$-structure described as follows. The conjugation action of $F$ on $R$, leads to a right $kF$-module structure on $k\otimes R_{\ab}$. However, the action of $R$ is trivial since $r^{t}\equiv r\mod [R, R]$ for all $r, t\in R$. So this induces a natural $kG$-modules structure on $k\otimes R_{\ab}$.
\end{defi}

\begin{proof} The sequence 
\begin{equation*}
\begin{tikzcd}
0\ar[r] &  kG^{(S)} \ar[r, "\partial_1"] & kG \ar[r, "\epsilon"] & k\ar[r] & 0, \\
\end{tikzcd}
\end{equation*}
is not exact at $kG^{(S)}$. It needs a correction of $H_1(R, k)$ and we compute this correction explicitly with the aid of the bar resolution. 
Consider the following commutative diagram of $kR$-modules, whose rows are resolutions of $k$ over $kR$,
\begin{equation*}
\begin{tikzcd}
\cdots \ar[r, "d_4"]  & \bC_3(R, k) \ar[r, "d_3"] \ar[d, "f_3"]] & \bC_2(R, k) \ar[d, "f_2"] \ar[r, "d_2"] & \bC_1(R, k) \ar[d, "f_1"] \ar[r, "d_1"] & kR \ar[d, "f_0"] \ar[r, "d_0"] & k \ar[d, "\id"] \ar[r] & 0 \\
\cdots \ar[r] & 0\ar[r] & 0 \ar[r] &  kF^{(S)} \ar[r, "\partial_1"] & kF \ar[r, "\partial_0"] & k\ar[r] & 0.
\end{tikzcd}
\end{equation*}
Here $f_0$ is the natural inclusion, $f_k=0$ for $k\geq 2$ and 
\[f_1\left(r|t|\right) = \sum_{s\in S} r\, \pd{t}{s}\, e_s, \, \, \, \mbox{for all $r, t\in R$.}\]
Since this the arrows are projective resolutions of $kR$-modules, the chain map $f_k$ is an homotopy equivalence. We apply the functor $k\tens{kR} $ and  observe that the resulting diagram is
\begin{equation*}
\begin{tikzcd}
\cdots \ar[r, "k\otimes d_4"]  & kR^3  \ar[r, "k\otimes d_3"] \ar[d, "0"] & kR^2\ar[d, "0"] \ar[r, "k\otimes d_2"] & kR \ar[d, "k\otimes d_1"] \ar[r, "0"] & k \ar[d, "\iota"] \ar[r, "\id"] & k \ar[d, "\id"] \ar[r] & 0, \\
\cdots \ar[r] & 0\ar[r] & 0 \ar[r] &  kG^{(S)} \ar[r, "\partial_1"] & kG \ar[r, "\partial_0"] & k\ar[r] & 0,
\end{tikzcd}
\end{equation*}
from where we can adapt the following exact sequence of abelian groups
\begin{equation*}
\begin{tikzcd}
0\ar[r] &  \frac{kR}{\im d_2\otimes k} \ar[r, "k\otimes f_1 "] &  kG^{(S)} \ar[r, "\partial_1"] & kG \ar[r, "\partial_0"] & k\ar[r] & 0.
\end{tikzcd}
\end{equation*}
By substituting the first non-zero factor by $k\otimes_{\Z} R_{\ab}$, which is isomorphic to $H_1(R, k)$ according to the isomorphism that is implicit in the proof of \cref{H_1}, then the corresponding $k\otimes f_1$ turns into $\partial_2$ and the result is the following exact sequence of abelian groups 
\begin{equation*}
\begin{tikzcd}
0\ar[r] &  k\otimes_{\Z} R_{\ab} \ar[r, "\partial_2"] &  kG^{(S)} \ar[r, "\partial_1"] & kG \ar[r, "\partial_0"] & k\ar[r] & 0,
\end{tikzcd}
\end{equation*}
whose maps are also seen to be $kG$-homomorphisms. The only map for which this is not obvious is the map $\partial_2$.
Let $r\in R$ and let $f\in F$. Notice that 
\begin{align*} \partial_2\left(frf^{-1}\, [R, R]\right) =& \sum_{s\in S}   \overline{\pd{frf^{-1}}{s}}\, e_s \\
=& \sum_{s\in S}   \overline{\pd{f}{s}+ f\pd{r}{s}-frf^{-1}\pd{f}{s} }\, e_s \\
=& \sum_{s\in S}   \overline{f\pd{r}{s}}\, e_s \\
=& \overline{f} \,  \partial_2\left(r\, [R, R]\right).\qedhere
\end{align*}

\end{proof}

One can view an algebraic proof of \cref{Relmod} that uses no homological machinery in \cite[Chapter 11]{Joh97}. A purely topological argument is given in \cite[Chapter II, Section 5]{Bro82}, by means of the Hopf's formula. Our argument is based on the indications given in \cite[Exercise 4, Chapter IV, Section 2]{Bro82}.

We have given a description of the augmentation ideal of a group in terms of a presentation. We are also interested in giving a description of the augmentation ideal of an an amalgamated product and of an HNN extension. However, we should first study free products. 

\begin{notation} Given two groups $\La\leq \Ga$, we denote by $kI_\La^{\Ga}$ the left ideal of $kG$ generated by $I_{\La}$. Since $k\Ga$ is a natural free right $k\La$-module, it follows that the canonical map of $k\Ga$-modules 
\[k \Ga\tens{k \La} k I_{\La}\lrar kI_{\La}^{\Ga}\]
is an isomorphism.
\end{notation}

\begin{lemma} \label{augfree} Let $\Gamma$ be the free product of two groups $A$ and $B$. Then the $k\Gamma$-homomorphism
\[kI_{A}^{\Gamma}\oplus kI_{B}^{\Gamma} \lrar kI_{\Gamma}, \]
defined by $(x, y)\mapsto x+y$ for $x\in kI_A$ and $y\in kI_B$, is an isomorphism.
\end{lemma}
\begin{proof} We name $\psi$ the given map. This statement can be verified by brute force using Britton's lemma. However, it is cleaner to refer to the correspondence $\Der (\Ga, M)\cong \Hom_{\Z \Ga}(I_{\Ga}, M)$ of \cref{deriv}. 
\end{proof}

\begin{thm}[Swan] \label{augamal} Let $\theta_1: C\rar A$ and $\theta_2: C\rar B$ be group homomorphisms and consider the induced amalgamated product   $\Gamma=A\underset{C}{*} B$. There is an exact sequence of $k\Gamma$-modules 
\begin{equation*}
\begin{tikzcd}
0\ar[r] & kI_C^{\Gamma} \ar[r, "\al"] &  kI_{A}^{\Gamma}\oplus kI_{B}^{\Gamma} \ar[r, "\be"] & kI_{\Gamma} \ar[r] & 0,
\end{tikzcd}
\end{equation*}
where $\al(z)=(\theta_1(z), -\theta_2(z))$ for all $z\in kI_C$ and $\be(x, y)=x+y$ for all $x\in I_A$, $y\in I_B$.
\end{thm}
\begin{proof} It is clear that $\beta\circ \al=0$, that $\al$ is injective and that $\be$ is surjective. We simply have to check that $\im \al=\ker \be$. 

Let $\Ga_0=A*B$. We consider the normal subgroup $R$ of $\Ga_0$ generated the elements of the set $\{\theta_1(c)^{-1}\, \theta_2(c): c\in C\}$. Then $\Ga\cong \Ga_0/R$. The idea is to consider an exact sequence using \cref{augfree} for the free product $\Gamma_0$, to take $R$-coinvariants, and to study the correction term that arises in the new sequence for $\Gamma$ with the aid of the bar resolution of $R$. 

Consider a diagram of the following form

\begin{equation*}
\begin{tikzcd}
\cdots \ar[r, "d_4"]  & \bC_3(R, k) \ar[r, "d_3"] \ar[d, "f_3"] & \bC_2(R, k) \ar[d, "f_2"] \ar[r, "d_2"] & \bC_1(R, k) \ar[d, "f_1"] \ar[r, "d_1"] & kR \ar[d, "f_0"] \ar[r, "d_0"] & k \ar[d, "\id"] \ar[r] & 0 \\
\cdots \ar[r] & 0\ar[r] & 0 \ar[r] & kI_{A}^{\Gamma_0}\oplus kI_{B}^{\Gamma_0} \ar[r, "\psi"] & k\Gamma_0 \ar[r, "\partial_0"] & k\ar[r] & 0.
\end{tikzcd}
\end{equation*}
The map $\psi$ is the isomorphism $ kI_{A}^{\Gamma_0}\oplus kI_{B}^{\Gamma_0} \rar I_{\Gamma_0}$ of \cref{augfree} while the maps $f_k$ are zero for $k\geq 2$. Lastly, 
the map $f_1$ is the $kR$-linear extension of  $f_1(|r|)=\psi^{-1}(r-1)$.  It is clear that $f_0\circ d_1=\psi\circ f_1$. Moreover,  since the map $\psi^{-1}(r-1): R\rar  kI_{A}^{\Gamma_0}\oplus kI_{B}^{\Gamma_0}$ is a $R$-derivation, then $f_1$ verifies that $f_1\circ d_2=0$.
So the previous diagram is commutative. The arrows are projective resolutions of $kR$-modules and the maps $f_i$ define a chain homotopy. Applying the functor $k\tens{kR}$ we obtain
\begin{equation*}
\begin{tikzcd}
\cdots \ar[r, "k\otimes d_4"]  & kR^3 \ar[r, "k\otimes d_3"] \ar[d, "f_3"] & kR^2 \ar[d, "f_2"] \ar[r, "k\otimes d_2"] & kR \ar[d, "f_1"] \ar[r, "0"] & k \ar[d, "f_0"] \ar[r, "\id"] & k \ar[d, "\id"] \ar[r] & 0 \\
\cdots \ar[r] & 0\ar[r] & 0 \ar[r] & kI_{A}^{\Gamma}\oplus kI_{B}^{\Gamma} \ar[r, "\beta"] & k\Gamma \ar[r, "\partial"] & k\ar[r] & 0.
\end{tikzcd}
\end{equation*}
Arguing as in \cref{Relmod}, we have an exact sequence of the form 
\begin{equation*}
\begin{tikzcd}
 0\ar[r] & k\otimes R_{\ab} \ar[r, "f_1"] & kI_{A}^{\Gamma}\oplus kI_{B}^{\Gamma} \ar[r, "\beta"] & k\Gamma \ar[r, "\partial"] & k\ar[r] & 0,
\end{tikzcd}
\end{equation*}
where the image of $r+[R, R]\in R_{\ab}$ under $f_1$ is $\phi^{-1}(r-1)$ for any $r\in R$. Since  $R_{\ab}$ is generated by the elements of the set $\{\theta_1^{-1}(c)\, \theta_2(c): c\in C\}$, the image of $f_1$ is generated by their corresponding images. Let $c\in C$. 

Notice that 
\[\theta_1^{-1}(c)\, \theta_2(c)-1=\theta_1(c)^{-1}\, \left((\theta_2(c)-1)-(\theta_1(c)-1)  \right),\]
and so
\[f_1(\theta_1^{-1}(c)\, \theta_2(c))=\theta_1(c)^{-1}\, (1-\theta_1(c), \theta_2(c)-1).\]
The $k\Ga$-span of these elements, for $c\in C$, generate the image of $f_1$, which, by exactness, is the kernel of $\be$. Crucially, this is  the image of $\al$.
\end{proof}

\begin{thm} Let $A$ be a subgroup of $H$ and let $\theta: A\rar H$ be a group monomorphism. Consider the corresponding HNN extension $\Ga=H\hn{A}$.  The $k\Ga$-homomorphism
\begin{equation*}
\begin{tikzcd}
  kI_H^{\Ga}\oplus k\Ga \ar[r, "\be"] & kI_{\Ga},
\end{tikzcd}
\end{equation*}
where $\be(x, y)=x+y(t-1)$ for all $x\in I_H$ and $y\in k\Ga$, is surjective and the kernel is $k\Ga$-generated by the set $\{(\theta(a)-1-t(a-1), \theta(a)-1): a\in A\}$.
\end{thm}
Our argument is essentially the same as in \cref{augamal}. 

\begin{proof}
Let $\Ga_0=H*(t)$ and let $R$ be the normal subgroup of $\Ga_0$ generated by the elements $\{a^{-1}t^{-1}\theta(a)t: a\in A\}.$ Then $\Ga\cong \Ga_0/R$.
Consider a diagram of the form

\begin{equation*}
\begin{tikzcd}
\cdots \ar[r, "d_4"]  & \bC_3(R, k) \ar[r, "d_3"] \ar[d, "f_3"] & \bC_2(R, k) \ar[d, "f_2"] \ar[r, "d_2"] & \bC_1(R, k) \ar[d, "f_1"] \ar[r, "d_1"] & kR \ar[d, "f_0"] \ar[r, "d_0"] & k \ar[d, "\id"] \ar[r] & 0 \\
\cdots \ar[r] & 0\ar[r] & 0 \ar[r] &  kI_H^{\Ga_0}\oplus k\Ga_0 \ar[r, "\psi"] & k\Gamma_0 \ar[r, "\partial_0"] & k\ar[r] & 0.
\end{tikzcd}
\end{equation*}
The map $\psi: kI_H^{\Ga_0}\oplus k\Ga_0\rar kI_{\Ga_0}$ is the isomorphism of \cref{augfree}, defined by $(x, y)\mapsto x+y(t-1)$ for all $x\in kI_H$ and $y\in k\Ga_0$. Here we have implicitly identified $ k\Ga_0\cong I_{\lan t \ran}^{\Ga_0}$ by means of the isomorphism 
\[z\mapsto z(t-1), \, \, \, \mbox{for all $z\in k\Ga_0$.} \]
On the other side, the map $f_0$ is the natural inclusion, the maps $f_k$ are zero for $k\geq 2$ and $f_1$ is the $k\Ga_0$-extension of $f_1(|r|)=\psi^{-1}(r-1),$ which defines, when restricted to $G$, a derivation $G\rar kI_H^{\Ga_0}\oplus k\Ga_0$. So $f_1\circ d_2=0$. It is also direct to see that $f_0\circ d_1=\psi\circ f_1$. So the previous diagram represents a chain homotopy $f_i$ between projective resolutions of $k$ over $kR$. Applying $k\tens{kR}$ we get the commutative diagram
\begin{equation*}
\begin{tikzcd}
\cdots \ar[r, "k\otimes d_4"]  & kR^3 \ar[r, "k\otimes d_3"] \ar[d, "f_3"] & kR^2 \ar[d, "f_2"] \ar[r, "k\otimes d_2"] & kR \ar[d, "f_1"] \ar[r, "0"] & k \ar[d, "f_0"] \ar[r, "\id"] & k \ar[d, "\id"] \ar[r] & 0 \\
\cdots \ar[r] & 0\ar[r] & 0 \ar[r] &  kI_H^{\Ga}\oplus k\Ga \ar[r, "\beta"] & k\Gamma \ar[r, "\partial"] & k\ar[r] & 0,
\end{tikzcd}
\end{equation*}
From this, we can obtain, arguing as in \cref{Relmod}, an exact sequence of $k\Ga$-modules
\begin{equation*}
\begin{tikzcd}
 0\ar[r] & k\otimes R_{\ab} \ar[r, "f_1"] & kI_{A}^{\Gamma}\oplus kI_{B}^{\Gamma} \ar[r, "\beta"] & k\Gamma \ar[r, "\partial"] & k\ar[r] & 0,
\end{tikzcd}
\end{equation*}
where the image of $r+[R, R]\in R_{\ab}$ under $f_1$ is $\phi^{-1}(r-1)$ for any $r\in R$. By exactness, the image of $f_1$ must be the kernel of $\beta$, and we want to check that this can be $k\Ga$-generated by the elements of $\{(\theta(a)-1-t(a-1), \theta(a)-1): a\in A\}$. The image of $f_1$ will be $k\Ga$-generated by the images under $f_1$ of  $\{a^{-1}t^{-1}\theta(a)t : a\in A\}$, since this set generates the abelian group $R_{\ab}$. Lastly, observe that, for any $a\in A$, 
\[a^{-1}t^{-1}\theta(a)t-1=a^{-1}t^{-1}\, \left( \theta(a)-1-t(a-1)+(\theta(a)-1)(t-1) \right).\]
Hence 
\[f_1(a^{-1}t^{-1}\theta(a)t+[R, R])= (\theta(a)-1-t(a-1),  \theta(a)-1).  \]and the conclusion follows.
\end{proof}

\begin{cor}\label{aughnn1} Let $A=\lan a\ran $ be a cyclic subgroup of $H$ and let $\theta: A\rar H$ be a group monomorphism. Consider the corresponding HNN extension $\Ga=H\hn{A}$. Then 
\[\frac{kI_H^{\Ga}\oplus k\Ga}{k\Ga\, \left( \theta(a)-1-t(a-1), \theta(a)-1\right)} \cong  I_{\Ga}, \]
where this isomorphism of $k\Ga$-modules is given by 
\[[(x, y)]\mapsto x+y(t-1),\, \, \, \mbox{for all $x\in I_H$ and $y\in k\Ga$.}\]
\end{cor}

\section{Back to group homomorphisms}

We wanted to study the maps that group homomorphisms induce between their augmentation ideals. Here will give two instances in which one can derive information from these induced maps back to the initial group homomorphisms. 

\begin{lem} \label{kernelmodp} Let $\phi : \ti{G}\rar G$ be a surjective group homomorphism of kernel $K$. Suppose that the natural  homomorphism of $k G$-modules \[kI_{\ti{G}}\tens{k \ti{G}} k G\lrar k I_G,\]
defined by the $k$-linear extension of $a \otimes b\mapsto \phi (a)b$, is an isomorphism. Then 
\[k\otimes_{\Z} K_{\ab}=0.\]
\end{lem} 
\begin{proof}
We can view any group as a quotient of a free group. We take a free group $F$ freely generated by some $S\sub F$ and a surjective map $F\rar \ti{G}$ of kernel $\ti{N}\unlhd F$. By composing this projection with $\phi$, we have a surjective $F\rar G$, of (possibly) bigger kernel $N\unlhd F$. By using the functoriality of the exact sequence of \cref{Relmod}, we obtain a commutative diagram of $k$-modules, with exact rows, of the form
\begin{equation*}
\begin{tikzcd}
 k\otimes \ti{N}_{\ab} \ar[d, "\iota_{\ab}"] \ar[r, "\partial_2"] & k\ti{G}^{(S)} \ar[d, "\phi"] \ar[r, "\partial_1"] & k I_{\ti{G}} \ar[d, "\phi"] \ar[r] & 0\ar[d] \\
  k\otimes N_{\ab}\ar[r, "\partial_2"] & kG^{(S)} \ar[r, "\partial_1"] & k I_G \ar[r] & 0,
\end{tikzcd}
\end{equation*}
where the unspecified arrows are the natural ones, induced either by $\phi: \ti{G}\rar G$ or the natural inclusion $\iota: \ti{N}\inc N$.

Applying the right-exact functor $\otimes_{k\ti{G}} kG$, where $kG$ has the natural $k\ti{G}$ structure induced by $\phi$, we get a commutative diagram of $kG$-modules with exact arrows
\begin{equation*}
\begin{tikzcd}
 \left(k\otimes \ti{N}_{\ab}\right)\otimes_{k\ti{G}} kG \ar[d, "\iota_{\ab}\otimes 1"] \ar[r, "\partial_2\otimes 1"] & kG^{(S)} \ar[d, "\id"] \ar[r, "\partial_1"] & k I_G \ar[d, "\id"] \ar[r] & 0\ar[d] \\
  k\otimes N_{\ab}\ar[r, "\partial_2"] & kG^{(S)} \ar[r, "\partial_1"] & k I_G \ar[r] & 0.
\end{tikzcd}
\end{equation*}
We look at the arrows that point downwards. Since the second arrow is surjective and the third one is injective, then $\iota_{\ab}\otimes 1$, the first arrow, is surjective. This is one of the so called ``four lemma'' on diagram-chasing. 
We now rewrite the map $\iota_{\ab}\otimes 1$. There is a commutative diagram of the form 
\begin{equation} \label{square1}
\begin{tikzcd}
 \left(k\otimes \ti{N}_{\ab}\right)\otimes_{k\ti{G}} kG \ar[d, "\iota_{\ab}\otimes 1"] \ar[r, "\psi"] & k\otimes \frac{\ti{N}}{[\ti{N}, N]} \ar[d, "\al"] \\
  k\otimes N_{\ab} \ar[r, "\id"] & k\otimes \frac{N}{[N, N]},
\end{tikzcd}
\end{equation}
where $\psi$ is an isomorphism defined by 
\[\psi\left( (1\otimes \ti{n}+[\ti{N}, \ti{N}])\otimes g\right)= 1 \otimes gng^{-1}+[\ti{N}, N]\]
and $\al$ is defined by 
\[\al (1\otimes \ti{n}+[\ti{N}, N])=1\otimes \iota (\ti{n})+[N, N].\]
The inverse of $\psi$ would be given by 
\[\psi^{-1}\left( \ti{n}+ [\ti{N}, N]\right)=(1\otimes  \ti{n}+[\ti{N}, \ti{N}])\otimes 1.\]
For the previous map to be well-defined, we observe the following in the module $\left(k\otimes \ti{N}_{\ab}\right)\otimes_{k\ti{G}} kG$.  For any $n\in N$, $\ti{n}\in \ti{N}$, 
\begin{align*} (1 \otimes n\ti{n}n^{-1}+[\ti{N}, \ti{N}])\otimes 1 =& (1\otimes (\ti{n} +[\ti{N}, \ti{N}])\, n)\otimes 1 \\
=& ( 1\otimes \ti{n}+[\ti{N}, \ti{N}])\otimes \phi(n)\cdot 1 \\
=& ( 1\otimes \ti{n}+[\ti{N}, \ti{N}])\otimes 1. \end{align*}
Hence, for any $n\in N$ and $\ti{n}\in \ti{N}$, 
\[(1\otimes [\ti{n}, n] +[\ti{N}, \ti{N}])\otimes 1 =0.\]
This verification says that $\psi^{-1}$ is also well-defined. It is easy to see that $\psi$ and $\psi^{-1}$ are inverses of each other. So $\psi$ is an isomorphism of $kG$-modules. 

We proved that $\iota_{\ab}\otimes 1$ is surjective. So the map $\al$ of the commutative square (\ref{square1}) must be surjective, too. 

From the natural exact sequence of $\Z G$-modules 
\begin{equation*}
    \begin{tikzcd}
     \frac{\ti{N}}{[\ti{N}, N]}\ar[r, "\al'"] & \frac{N}{[N, N]}\ar[r, "\be'"] & \frac{N}{\ti{N}[N, N]}\ar[r] & 0.
    \end{tikzcd}
\end{equation*} 
we obtain, by applying the right exact functor $k \tens{\Z}$, an exact sequence of $kG$-modules 
\begin{equation*}
    \begin{tikzcd}
     k\otimes  \frac{\ti{N}}{[\ti{N}, N]}\ar[r, "\al"] & k\otimes \frac{N}{[N, N]}\ar[r, "\be"] & k\otimes \frac{N}{\ti{N}[N, N]}\ar[r] & 0.
    \end{tikzcd}
\end{equation*} 
We have already deduced that the map $\al$ is surjective, so the image of $\beta$ is zero and then the third module must be zero. Observe that $K\cong N/\ti{N}\leq F/\ti{N}\cong \ti{G}$, so  
\[k\otimes \frac{N}{\ti{N}[N, N]} \cong k\otimes K_{\ab}\cong 0.\qedhere\]
\end{proof}

A particular realisation of the previous lemma is a result due to Lewin.

\begin{cor}[Lewin] \label{Lewin} Let $G$ be a group generated by two subgroups $A$ and $B$ with $C=A\bigcap B$. Suppose that $kI_C^{G}=kI_C^G\bigcap kI_B^G$. Then the natural map  $A\underset{C}{*} B\lrar G$ is an isomorphism. 
\end{cor}
\begin{proof}
Denote $\ti{G}=A\underset{C}{*} B$ and name $\phi: \ti{G}\rar G$ the natural map, which is surjective.

Consider a commutative diagram of $k$-modules of the form
\begin{equation*}
    \begin{tikzcd}
    kI_C^{\ti{G}} \ar[r, "\al"] \ar[d, "\phi"] & kI_A^{\ti{G}}\oplus kI_{kB}^{\ti{G}} \ar[d, "\phi"] \ar[r, "\be"] & kI_{\ti{G}} \ar[d, "\phi"] \ar[r]& 0\ar[d] \\
    kI_C^{G} \ar[r, "\al"] & kI_A^{G}\oplus kI_{B}^{G} \ar[r, "\be"] & kI_G \ar[r]& 0,
    \end{tikzcd}
\end{equation*} 
where the arrows are the obvious ones. We define  $\al(z)=(z, -z)$ for all $z\in I_C$ and $\be(x, y)=x+y$ for all $x\in I_A$ and $y\in I_B$. By \cref{augamal}, the first row is an exact sequence of $k\ti{G}$-modules. Furthermore, by assumption, the second row is an exact sequence of $kG$-modules. 

We apply the right-exact functor $kG\tens{k\ti{G}}$ to the first row and we obtain a commutative diagram of $kG$-modules
\begin{equation*}
    \begin{tikzcd}
    kI_C^{G} \ar[r, "\al"] \ar[d, "\phi"] & kI_A^{G}\oplus kI_{kB}^{G} \ar[d, "\phi"] \ar[r, "\be"] & k G\otimes_{k\ti{G}} kI_{\ti{G}} \ar[d, "\phi"] \ar[r]& 0\ar[d] \\
    kI_C^{G} \ar[r, "\al"] & kI_A^{G}\oplus kI_{B}^{G} \ar[r, "\be"] & kI_G \ar[r]& 0,
    \end{tikzcd}
\end{equation*}
with exact rows. 

Now we look at the arrows that point downwards. The first is an epimorphism, and the second and fourth are monomorphism. A version of the four lemma ensures that the third arrow  $1\otimes \phi: k G\otimes_{k\ti{G}} kI_{\ti{G}}\rar kI_G$ is a monomorphism. Moreover, this map is surjective, since $\phi$ itself is surjective. By \cref{kernelmodp}, the map $\phi$ is an isomorphism. 
\end{proof}

\chapter{Universal division rings of fractions} \label{universalsection}

We want to emphasise our aim in talking about the concept of ``universal division rings of fractions'' before going into details. 

Let $\ti{G}$ be a candidate to being a parafree group. Suppose we already know that its pro-$p$ completions $\ti{G}_{\hat{p}}$ are free. Suppose, in addition, we also know that, for a suitable prime $p$, the kernel of the canonical map $\phi: \ti{G}\lrar \ti{G}_{\hat{p}}$ is free. Let $G$ be the image of $\ti{G}$ under the previous map. 
We consider the canonical $\p G$-module 
\[\p G\tens{\p \ti{G}} \p I_{\ti{G}}\]
and the canonical map of $\p G$-modules 
\begin{equation} \label{keymap} \p G\tens{\p \ti{G}} \p I_{\ti{G}}\lrar \p I_G.\end{equation}
If we check that the previous map is an isomorphism, we would deduce, by \cref{kernelmodp}, that
\[\F_p\tens{\Z} \ker \phi=0,\]
implying that $\ker \phi=1$, since it is free. 
In order to check that the map of (\ref{keymap}) is an isomorphism, we develop a ``dimension function'' on $\p G$-modules. This map is clearly surjective and, assuming reasonable properties about this ``dimension function'', if we prove that the dimension of both involved modules  is the same, then the surjective map of (\ref{keymap}) would be an isomorphism, as we want. 

This dimension function would be defined as follows. Let $\p G\inc \D$ be an embedding of  $\p G$ into a division ring $\D$. A tentative  
``dimension function'' for the collection of $\p G$-modules $M$ goes as follows. We can define 
\[\dim M=\dim_{\D} \D\tens{\p G} M.\]
However, for this dimension function to enjoy desirable properties, the ring $\D$ cannot be any division ring. It has to be ``minimal'' and ``universal'' in some sense. We will now give a brief account of how these notions are precised for some classes of rings, before introducing the universal division ring of fractions of $\p G$. 

\section{Embeddability of domains into division rings}

For a ring to be embeddable into a division ring, this ring must be a domain. In the commutative setting, this condition is not only necessary but also sufficient for the existence of such embedding. In fact, given a commutative domain  $A$, there is the standard construction of the field of fractions $\Frac(A)$ of $A$ with a canonical embedding $A\inc \Frac(A)$. It satisfies that, for any ring homomorphism $f: A\rar B$ such that $f(A\setminus\{0\} )\subset B^{\times}$, there exists a unique $\hat{f}: \Frac(A)\rar B$ such that $\hat{f}\circ \iota=f$; that is, such that the diagram
\begin{equation} \label{uni1}
\begin{tikzcd}
A\arrow[r,hook, "\iota"] \arrow[dr, rightarrow, "f"] & \Frac(A) \arrow[d, dashrightarrow, "\hat{f}"]\\& B
\end{tikzcd}
\end{equation}
is commutative. 

Underlying the construction of $\Frac(A)$, there is the more general notion of \textit{localization}, of fundamental importance in algebra and geometry. There is also a non-commutative analogue of this notion, namely the \textit{localization in the sense of Ore}. 

Given a ring $R$ and a (multiplicatively closed) subset $S\sub R$, we ask for the existence of a ring $R_S$ equipped with an injective map $\phi_S: R\lrar R_S$ such that
\begin{enumerate}
    \item the map $\phi_S$ is injective; 
    \item any element of $R_S$ can be written in the form $\phi_S(r)\phi_S(s)^{-1}$ for some $r\in R$ and $s\in S$; and
    \item  for any ring homomorphism $\psi: R\rar B $ with $\psi(S)\subset B^{\times} $, there exists a unique $\hat{\psi_S}: R_S\rar B$ such that $\psi=\hat{\psi}\circ \phi_S$. In other words, the diagram
\begin{equation} \label{uni2}
\begin{tikzcd}
A\arrow[r, rightarrow, "\phi_S"] \arrow[dr, rightarrow, "\psi"] & S^{-1}A \arrow[d, dashrightarrow, "\hat{\psi}"]\\& B,
\end{tikzcd}
\end{equation}
would be commutative.
\end{enumerate}

The natural candidate for $R_S$, and, in fact, the only one, is constructed from $R$ by attaching inverses $t_s$ to each element $s\in S$ as follows 
\begin{equation}\label{RS} R_S=\frac{R\,  \lan t_s; s\in S  \ran }{\lan t_s s-1, st_s-1; s\in S \ran}.\end{equation}

It is easy to check that it the universal property of  the diagram (\ref{uni2}). It would remain to ask under which conditions can we ensure that the canonical map $\phi_S: R\rar R_S $ is injective and that the elements of $R_S$ can be written as $a\, s^{-1}$. In this case, we name $R_S$ the \textit{right Ore localisation of $R$ with respect to $S$}. We give sufficient conditions due to 
Asano and Gabriel, inspired by  the work of Ore \cite{Ore31}.

\begin{defi} We say that  $S\subset R$ is a right Ore set of regular elements if 
\begin{enumerate}
    \item $S$ is closed under multiplication.
    \item For any $a\in R$ and $s\in S$, $as=0$ implies that $a=0$.
    \item For any $a\in R$ and $s\in S$, $aS\bigcap sR\neq \emptyset$. 
\end{enumerate}
\end{defi}

The following statement is taken from the book of Cohn \cite[Theorem 0.7.1]{Coh06}. 
\begin{thm}[Criterion for faithful localization] \label{Oree} Let $S\subset R$ be closed under multiplication. Then the canonical map $R\rar R_S$, defined in (\ref{RS}), is an embedding and every element of $R_S$ can be written in the form $r\, s^{-1}$, for $r\in R$, $s\in S$, if and only if $S$ is a right Ore set of regular elements.  
\end{thm}

\begin{defi} A ring $R$ is said to be a right Ore domain if $S=R\setminus \{0\}$ is a right Ore set of regular elements. In this case, $R$ can be embedded into $R_S$, which is a division ring generated by $R$. We name  $\Qo(R)=R_S$ the  {\bf (right) Ore division ring of fractions}. 
\end{defi}

It is not very difficult to find a big family of non-commutative Ore domains.

\begin{prop}[``Little Goldie's theorem'', 1957] \label{Goldie} Let $R$ be a right Noetherian domain, then $S=R\setminus \{0\}$ is a right Ore set of regular elements. In particular, $R$ can be embedded into its right Ore  ring of fractions, which is a division ring. 
\end{prop}
\begin{proof}
To check that $S=R\setminus \{0\}$ is a right Ore set of regular elements, we only have to verify that for any $a\in R$ and $s\in S$, $aS\bigcap sR\neq \emptyset$. If $a=0$, the previous intersection contains zero. Suppose that $a\neq 0$ and consider the chain of right ideals $I_m=(s, as, a^2s, \dots, a^{m}s)$, which is increasing, so it has to stabilise. In particular, there exists $n\geq 1$ so 
\[a^ns=\sum_{k=0}^{n-1} a^k\, s\, r_k,\, \, \, \mbox{for some $r_k\in R$.}\]
Consider the minimal $k_0$ so $r_{k_0}\neq 0$, then
\[a\left(a^{n-k_0-1}s-\sum_{k=k_0}^{n-2} a^{k-k_0}\, s\, r_k\right)=sr_{k_0}\in aS\bigcap sR. \qedhere\]
\end{proof}

This produces many examples of Ore domains, such as polynomial algebras.
We review the situation by saying that an {\bf Ore domain $R$} can be embedded into its Ore division ring of fractions $\iota_{\ore} :R\inc \Qo (R)$ and that this embedding verifies the following {\bf universal property}. For any ring homomorphism $\psi: R\lrar B$ such that $\psi(R\setminus \{0\})\sub B^*$, there exists a unique $\hat{\psi}$ such that the following diagram 
\begin{equation} \label{uni3}
\begin{tikzcd}
R \arrow[r, rightarrow, "\iota_{\ore}"] \arrow[dr, rightarrow, "\psi"] & \Qo (R) \arrow[d, dashrightarrow, "\hat{\psi}"]\\& B
\end{tikzcd}
\end{equation}
is commutative. In addition, if  $\psi$ is injective, the resulting $\hat{\psi}$ is injective, too. Moreover, $\Qo(R)$ is a flat $R$-module. This means that the functor $\tens{R} \Qo (R)$, from the category of $R$-modules to the category of abelian groups; preserves exactness of short exact sequences. 

When discussing the question of embedability of domains into division rings, it is natural, for historical reasons, to introduce  the class of Ore domains. Furthermore,  these domains will be part of our arguments in \cref{torfreesection}.  However, the class of Ore domains is still a very restrictive one. In fact, a result of Bartholdi and Kielak \cite{Bar16} states that a group ring $KG$, with coefficients in a field $K$, is an Ore domain if and only if $KG$ is a domain and $G$ is amenable. For this reason, the previous methods will not be directly applicable to our groups of interest, namely abstract subgroups $G$  of free pro-$p$ groups $\FF$. The crucial aspect of these groups $G$ is that they are orderable. 

Malcev and Neumann discovered independently a way to embed the group ring of an orderable group  into a division rings of power series. For the group $\Z$, notice that $K\Z\cong K[X, X^{-1}]$. The embedding provided by Malcev and Neumann, where $\Z$ is viewed with the standard ordering, would be $K[X, X^{-1}]\inc K((X))$, where $K((X))$ denotes the ring of one-variable Laurent series with coefficients in the field $K$.

\begin{defi}[{\bf Malcev-Neumann construction of $D((G))$}] Given a division ring $D$ and an orderable group $G$, we define the formal power series of $G$ over $D$, denoted $D((G))$, as the set of formal sums
\[\left\{ \sum_{g\in G} a_g \, g : \mbox{the support $\{g: a_g\neq 0\}$ is well-ordered}\right\}.\]
The operations on $D((G))$ are defined analogously as in $K((X)).$

\begin{equation} \label{sum2}\left(\sum_g a_g\, g\right) +\left(\sum_g b_g\,  g\right)=\sum_g (a_g+b_g) \, g.\end{equation}
\begin{equation}\label{prod2} \left(\sum_g a_g \, g\right)\left(\sum_h b_h \, h\right)=\sum_k \left(\sum_{gh=k} a_g\,  b_h\right) k.\end{equation}
\end{defi}

The difficult technical lemma about the order on $G$ which allows one to prove that (\ref{prod2}) is well defined and that $D((G))$ is a division ring is the following. 

\begin{lemma} Let $G$ be an orderable group and let $S\subset G$ be a well-ordered subset such that $s>1$ for every $s\in S$. Consider \[S^{\omega}=\bigcup_{n\geq 1} S^{n}.\] Then $S^{\omega}$ is well-ordered and every  $s\in S^{\omega}$ appears only in finitely many sets $S^n$. 
\end{lemma}

Given an element $0\neq \sum_{g\in G} a_g \, g \in D((G))$, we use the previous lemma to explicitly write down its inverse. Since its support is well-ordered, it has a minimum $g_0\in G$, and we can write 
\[0\neq \sum_{g\in G} a_g \, g=a_{g_0} \, g_0\,  \left( 1+ \sum_{1<g\in G} b_g \, g \right). \]
This way, writing $c=\sum_{1<g\in G} b_g \, g$,  its inverse will be 
\begin{equation} \label{inverse} \left( \sum_{k=0}^{\infty} (-1)^{k} c^k \right) a_{g_0}^{-1}\, g_0^{-1}.\end{equation}

We refer to the paper of Neumann \cite[Part I]{Neu49} for details.

This provides embeddings into division rings for group rings of orderable groups. However, for these group rings, there may not exist a division ring $\D$ and a monomorphism $R\inc \D$ verifying a universal property in the sense of the commutative diagram (\ref{uni3}). 

There is still a suitable and weakened notion of universality for division rings, due to Cohn, which we shall discuss next.

\section{Definition and existence of the ring $\D_{R}$}

\begin{defi}[Division closure]
Given a ring extension $R_0\sub \D_0$, where $\D_0$ is a division ring, we say that $\D\sub \D_0$ is the division closure of $R_0$ in $\D_0$ if $R_0$ generates rationally $\D_0$. In other words, each element of $\D_0$ can be built up from the elements of $R_0$ in stages, using addition, subtraction, multiplication, and division by nonzero elements.
We denote the division closure of $R_0$ in $\D_0$ by $\D=\D(R_0, \D_0)$ or, simply, by $\D(R_0)$. 
\end{defi}

\begin{defi} A $R$-ring is a homomorphism $\phi : R\rar S$. Two $R$-rings $\phi_1: R\rar S_1$ and $\phi_2: R\rar S_2$ are said to be isomorphic if there exists a ring isomorphism $\alpha: S_1\rar S_2$ such that $\phi_2=\alpha\circ \phi_1$. 
\end{defi}

We now precise Cohn's notion of universal division ring of fractions. 
\begin{defi}
Consider a $R$-ring $\phi: R\lrar \D$. We also say that it is 
\begin{enumerate}
    \item a division $R$-ring if $\phi$ is injective and $\D$ is a division ring;
    \item a $R$-ring of fractions if $\D$ is generated by $\phi(R)$.
\end{enumerate}
If the previous conditions hold, then we say that  $\phi: R\inc \D$ is a {\bf division ring of fractions of $R$}.
\end{defi}

\begin{defi}
Let $R$ be a ring and let $\D_1$ and $\D_2$ be two division $R$-rings of fractions.
\begin{itemize}
    \item We say write $\dim_{\D_1}\leq \dim_{\D_2}$ if   any finitely presented module $M$, it follows that $\dim_{\D} \D\tens{R} M\leq \dim_{\D'} D'\tens{R} M$. 
    \item We say that there is a specialization from $\D_1$ to $\D_2$ if there exists a local subring $R\sub B\sub \D_1$ with residual division ring isomorphic to $\D_2$. 
\end{itemize} 
\end{defi}

\begin{defi} Let $\D_R$ be a division $R$-ring of fractions. We say that it is a {\bf universal division $R$-ring of fractions} if, for any other division $R$-ring of fractions $\D'$, $\dim_{\D}\leq \dim_{\D'}.$
\end{defi}
We can give a few families of examples. 

\begin{eg} The following rings $R$ have a universal division $R$-ring of fractions $\D_R$.
\begin{itemize}
    \item Commutative domains $A$, for which $\D_A\cong \Frac(A)$. 
    \item More generally, right (resp. left) Ore domains $R$, for which $\D_R$ is their right (resp. left) Ore division ring of fractions $\D_R\cong \Qo (R)$.
    \item Semifirs \cite{Coh74}: A ring $R$ is said to be a semifir if every finitely generated left ideal of $R$ is free of fixed rank. 
    \item In particular, if $K$ is a field, then $R=K\lan \lan X\ran \ran$ has universal $\D_R$ since $R$ is a semifir by \cite[Proposition 2.9.19]{Coh06}).
\end{itemize}
\end{eg}

Another significant example for us is the completed group algebra $\p [[\FF]]$. We know, from \cref{completediso}, that this ring is isomorphic to $\p\lan \lan X\ran \ran$, where $|X|=d(\FF)$. From the previous examples, we know the following. 

\begin{prop} 
The ring $\p [[\FF]]$ has a universal division ring of fractions $\D_{\p [[\FF]]}.$
\end{prop}

As we remarked before, there is no chance in having universal property for the embedding $R\inc \D_R$ as in the setting of Ore division rings of fractions (\ref{uni3}). In fact,  Herbera and Sánchez \cite{Her15} prove that, if $R=K\lan X\ran$ with $|X|>1$,  there are infinitely many non-isomorphic $R$-division rings of fractions.  However, for some favourable situations one still expects to be able to compare different $R$-division rings. For example, let $G$ be a subgroup of $\FF$. Then $\p G$ embeds into $\D_{\p [[\FF]]}$ and also into $\p((G))$, after giving $G$ some ordering. It is not easy to compare both embeddings nor to study whether one of them turns out to provide a universal division $\p G$-ring. Both questions are answered positively. Interestingly, Cohn and  Hughes developed some methods that allow us to compare different $R$-division rings and, in particular, to prove the previous claim about  $\p G-$division rings of fractions. First of all, let us mention a more elementary criterion that can be used to compare division $R$-rings for Ore domains $R$.

\begin{lemma} \label{Orecriterion} Suppose that $R$ is an (left or right) Ore domain and that $R\inc \D$ is a division $R$-ring of fractions. Then $\D$ and $\Qo (R)$ are isomorphic $R$-rings. 
\end{lemma}
\begin{proof} By the universal property (\ref{uni3}) of the Ore division  $R$-ring of fractions  $R\inc \Qo(R)$, there is a commutative diagram
\begin{equation} 
\begin{tikzcd}
R \arrow[r, hookrightarrow, "\iota_{\ore}"] \arrow[dr, hookrightarrow, "\psi"] & \Qo (R) \arrow[d, dashrightarrow, "\hat{\psi}"]\\& \D,
\end{tikzcd}
\end{equation}
with $\hat{\psi}$ injective. Furthermore, since $R\inc \D$ is a division $R$-ring of fractions, then $\hat{\psi}$ is surjective. So $\hat{\psi}$ defines an isomorphism between the $R$-rings $\D$ and $\Qo(R)$. 
\end{proof}

We will talk about Hughes criterion later, we start with two results of Cohn \cite[Theorem 4.4.1 and Subsection 4.1]{Coh95}.

\begin{prop} \label{Cohn1} Let $R$ be a ring and let $\D_1$ and $\D_2$ be two division $R$-rings of fractions. Then the following holds.
\begin{itemize}
    \item The rings $\D_1$ and $\D_2$ are isomorphic, as $R$-rings,  if and only if $\dim_{\D_1}=\dim_{\D_2}$. 
    \item There is a specialization from $\D_1$ to $\D_2$ if and only if $\dim_{\D_1}\leq \dim_{\D_2}$.
\end{itemize}
In particular, if the universal division ring of fractions $\D_R$ exists, then it is unique. 
\end{prop}

\begin{defi} Let $R\inc \D_R$ be a universal division $R$-ring of fractions. Given a $R$-module $M$, we denote 
\[\dim_R M=\dim_{\D}\D\tens{R} M. \]
By \cref{Cohn1}, this is well-defined.
\end{defi}

In \cref{DpFring}, we will use the criterion of \cref{Cohn1} to check that for any closed subgroup $\HH$ of $\FF$, the division closure of $\p [[\HH]]$ in $\D_{\p [[\FF]]}$ is isomorphic to $\D_{\p [[\HH]]}.$

Using this fact, we would later prove, in \cref{DpGring},   that $\D_{\p G}$ exists and that the division closure of $\p G$ in $\D_{\p [[\FF]]}$ is isomorphic to $\D_{\p G}$. However, in this case, we will not use \cref{Cohn1}, but a criterion of Hughes on group rings of locally indicable groups (\cref{locind}).

\section{Sylvester rank functions}

Before establishing properties about the rings $\D_{\p G}$ and $\D_{\p [[\FF]]}$, we discuss some properties about their induced dimension functions in terms of Sylvester rank functions. 

Let $M$ be a finitely presented $R$-module. This is, there is an exact sequence of $R$-modules 
\[R^m \lrar R^n \lrar M \lrar 0. \]
In other words, there is a matrix $A\in M_{m\times n}(R)$ such that 
\[M\cong \frac{R^{n}}{R^m A}.\]
The consideration of a reasonable ``dimension function'' $\dim $ on $R$-modules $M$ is related to a reasonable ``rank function'' on matrices $A\in M_{m\times n} (R)$. 

For example, let $\D$ be a division $R$-ring. Then 
\begin{equation} \label{dimrank} \dim_{\D} \D\tens{R} M= n-\text{rk}_{\D} (A),\end{equation}
where $\text{rk}_{\D}(A)$ is the rank of $A$ as a matrix with coefficients in the division ring $\D$. 

For more general rings, that are not division rings, we do not expect a notion of rank as well behaved as the one coming from linear algebra, though we are interested in preserving some of them. 

\begin{defi} \label{Sylmatrix} Let $R$ be a ring. A {\bf Sylvester matrix rank function} $\rk$ on $R$ is a function that assigns a non-negative real number to each matrix over $R$ and that satisfies the following conditions. 
\begin{enumerate}
    \item $\rk(A)=0$ if $A$ is any zero matrix and $\rk(1)=1$.
    \item For all matrices $A_1, A_2$ \[\rk \begin{pmatrix} A_1 & 0\\
    0 & A_2 \end{pmatrix}=\rk (A_1)+\rk (A_2).\]
    \item For all matrices $A_1, A_2, A_3$, of appropriate sizes, \[\rk \begin{pmatrix} A_1 & A_3\\
    0 & A_2 \end{pmatrix}\geq \rk (A_1)+\rk (A_2).\]
\end{enumerate}
\end{defi}

Given a Sylvester matrix rank function $\rk$ on $R$, we have a ``dimension'' function on finitely-presented $R$-modules by setting $\dim (R^n/R^m A)=n-\rk (A)$, as occurs in (\ref{dimrank}). This $\rk$ enjoys the following properties. 

\begin{defi} \label{Sylmod} A {\bf Sylvester module rank function} on $R$ is a function that assings a non-negative real number to each finitely presented $R$-module and that satisfies the following conditions. 
\begin{enumerate}
    \item $\dim (\{0\})=0$ and $\dim R=1$. 
    \item $\dim (M_1\oplus M_2)=\dim M_1+\dim M_2$, and
    \item if $M_1\lrar M_2\lrar M_3\lrar 0$ is exact, then 
    \[\dim M_1+\dim M_3\geq \dim M_2\geq \dim M_3.\]
\end{enumerate}
\end{defi}

By previous comments, we know there is a natural correspondence between Sylvester matrix rank functions and Sylvester matrix module functions. Our next aim is to discuss a few aspects about the convergence of these Sylvester rank functions. However, these notions uniquely provide a ``dimension function'' for finitely presented modules. Recall that we are also interested in finitely generated modules. Given a Sylvester matrix rank function, Li \cite{Li20} provides a way to extend the induced dimension function to finitely generated modules so it verifies some favourable properties, as those listed in \cref{Sylmod}. 

\begin{defi} Given a Sylvester module rank function $\dim$, we define its {\bf extension} to finitely generated modules $M$ as follows. 
\[\dim M=\inf\{\dim \ti{M}: \mbox{ $\ti{M}$ is finitely presented and $M$ is a quotient of $\ti{M}$}\}. \]
\end{defi}

The following lemma reflects that this is a natural extension for the type of module rank functions we are interested in. 

\begin{lemma} Let $R\inc \D$ be a division $R$-ring. We consider the Sylvester module rank function $\dim$ defined by 
\[\dim M=\dim_{\D}\D\tens{R}M.\]
Then the extension $\ti{\dim}$ to finitely generated modules $M$ verifies that 
\[\ti{\dim} M=\dim_{\D} \D\tens{R} M.\]
\end{lemma}

This exhibits that in order to study these dimension functions, it suffices to restrict ourselves to the study of finitely presented modules. However, for other Sylvester rank functions, their extensions are not as easy to depict. Before describing more precise examples about this phenomena, we introduce the useful notion of convergence in Sylvester rank functions.

\begin{defi} Let $\dim_1$ and $\dim_2$ be two Sylvester module rank functions on $R$. We write $\dim_1\geq \dim_2$ if, for every finitely presented $R$-module $M$, $\dim_1 M\geq \dim_2 M$. 
\end{defi}

By \cref{Cohn1}, we can say that the universal division $R$-ring of fractions $\D_R$ is the the division $R$-ring of fractions which induces the minimal Sylvester module rank function $\dim_R$. 
 
\begin{defi}[Convergence] Let $\dim, \dim_i$ be two Sylvester module rank functions on $R$. We write $\rk=\lim_{i\rar \infty} \rk_i$ if, for every finitely presented $R$-module $M$, $\dim M=\lim_{i\rar \infty} \dim_i M$. 
\end{defi}

Equivalently, one could define the previous notion of convergence by setting a right topology in the collection of Sylvester matrix rank functions. Precisely, we denote by $\P(R) $ the collection of Sylvester matrix rank functions on $R$ and we view $\P(R)$ as a subspace of the space of real-valued functions on matrices,  $\R^\Mat(R)$. This bigger space is endowed  with the topology of pointwise convergence. Then $\P(R)$ is a compact convex subset of $\R^{\Mat (R)}$. 

There are many interesting problems, for example,  in the surroundings of Lück approximation conjectures, which can be rephrased, and generalised, in terms of convergence principles of certain  Sylvester matrix rank functions.

\begin{lemma}\label{conver} Let $M$ be finitely generated $R$-module and suppose that $\dim=\lim_{i\rar \infty} \dim_i$. Then
\[\dim M\geq \limsup_{i\rar \infty} \dim M_i.\]
If, in addition, $\dim_i\geq \dim$ for all $i$, then 
\[\dim M= \lim_{i\rar \infty} \dim M_i, \, \, \, \, \mbox{for all finitely generated $R$-modules $M$.}\]
\end{lemma}
\begin{proof}
Let $W$ be the set of all finitely presented $R$-modules $\ti{M}$ such that $M$ is a quotient of $\ti{M}$. Then
\begin{align*}
    \dim M=& \inf_{\ti{M}\in W} \dim \ti{M}\\
    =& \inf_{M\in W} \limsup_{i\rar \infty} \dim_i \ti{M}\\
    \geq&  \limsup_{i\rar \infty} \inf_{M\in W} \dim_i \ti{M}\\
    =& \limsup_{i\rar \infty} \dim_i M,
\end{align*}
and the first part of the statement follows. If, in addition, 
 $\dim_i\geq \dim$, then it is clear that
 \[\liminf_{i\rar \infty} \dim_i M\geq \dim M,\]
 so it follows that $\lim_{i\rar \infty} \dim_i M$ exists and equals $\dim M$. 
\end{proof}

\section{The ring $\D_{\p [[\FF]]}$} \label{DpFring}
Our next aim is to study concrete convergence of Sylvester module rank functions of $\pF$, following the treatment of  \cite[Section 3.1]{And20}.
\begin{prop} \label{dimapp} Let $\FF=\NN_1>\NN_2>\cdots$ be a chain of open normal subgroups with trivial intersection. Let $M$ be a finitely generated $\p [\FF]]$-module. Then 
\[\dim_{\p [[\FF]]}M =\lim_{i\rar \infty} \frac{\dim_{\p} \left(\p [\FF/\NN_i]\tens{\p[[\FF]]} M \right)}{|\FF: \NN_i|}. \]
\end{prop}
\begin{proof}
Let $\NN$ be a normal open subgroup of $\FF$, we define the Sylvester module rank function $\dim_{\p [\FF/\NN]} $ on $\p[[\FF]]$ by
\[\dim_{\p [\FF/\NN]} M=  \frac{\dim_{\p} \left(\p [\FF/\NN]\tens{\p[[\FF]]} M \right)}{|\FF: \NN|}. \]

We claim that  \begin{equation} \label{des0} \dim_{\p [\FF/\NN]} \geq \dim_{\p [[\FF]]}.\end{equation} In fact, let $M$ be a finitely generated  $\p [[\FF]]$-module with 
\[\dim_{\p [[\FF]]}=\dim_{\D_{\p [[\FF]]}} \D_{\p[[\FF]]}\tens{\p [[\FF]]} M=k.\]
Then there exist $k$ elements $\{m_1, \dots, m_k\}\sub M$ such that 
$\{1\otimes m_1, \dots, 1\otimes m_k\}$ are $\D_{\p[[\FF]]}$-linearly independent in $ \D_{\p[[\FF]]}\tens{\p [[\FF]]} M$. In particular, the $\p[[\FF]]$-submodule $N$ spanned by $\{m_1, \dots, m_k\}$ is free of rank $k$. So $\p [\FF/\NN]\tens{\p[[\FF]]} N\cong \p[\FF/\NN]^k$ is a $\p$-subspace of $\p [\FF/\NN]\tens{\p[[\FF]]} M$ of dimension $k|\FF: \NN|$, implying that 
\[\dim_{\p [\FF/\NN]}(M) \geq \frac{\dim_{\p} \left(\p [\FF/\NN]\tens{\p[[\FF]]} N \right)}{|\FF: \NN|}=k= \dim_{\p [[\FF]]}M.\]
This proves (\ref{des0}).

Furthermore, if $M$ is a finitely presented $\p [[\FF]]$-module, then \cite[Theorem 1.4]{And191} implies that 
\begin{equation} \label{des1} \dim_{\p [[\FF]]} M=\lim_{i\rar \infty} \dim_{\p [\FF/\NN_i]} M.\end{equation} 
With (\ref{des0}) and (\ref{des1}), we can apply \cref{conver} to conclude that for all finitely generated $\p[[\FF]]$-module $M$, 
\[\dim_{\p [[\FF]]}M=\lim_{i\rar \infty}  \dim_{\p [\FF/\NN_i]}M.\qedhere\]
\end{proof}

The next proposition will allow us to view $\D_{\p [[\HH]]}$ inside $\D_{\p [[\FF]]}$ whenever $\HH$ is an open subgroup of $\FF$.

\begin{prop} \label{HFprop} Let $\HH$ be a closed finitely generated subgroup of $\FF$. The following holds. 
\begin{itemize}
    \item[(a)]  Let $\D_{\HH}$ be the division closure of $\p[[\HH]]$ in $\D_{\p[[\FF]]}$. Then $\D_H$ and $\D_{\p[[\HH]]}$ are isomorphic $\p[[\HH]]$-rings.
    \item[(b)] If $M$ is a finitely generated $\p[[\HH]]$-module, then 
    \[\dim_{\p[[\HH]]} M=\dim_{\p[[\FF]]}\left( \p[[\FF]]\tens{\p[[\HH]]}M\right).\]
    \item[(c)] If $\HH$ is open, then 
    \[\dim_{\D_{\HH}} \D_{\p[[\FF]]}=|\FF: \HH|.\]
    \item[(d)] If $\HH$ is open, then 
    \[\D_{\p[[\FF]]}\cong \p[[\FF]]\tens{\p[[\HH]]}\D_{\p[[\HH]]}\]
    canonically as $(\p[[\FF]], \p[[\HH]])$-bimodules. 
\end{itemize}
\end{prop}

\begin{proof} By \cref{subfree}, the pro-$p$ group $\HH$ is free. Let $\FF=\NN_1>\NN_2>\cdots$ be a chain of open normal subgroups with trivial intersection. Define $\HH_i= \NN_i\bigcap \HH$. Then   $\HH=\HH_1>\HH_2>\cdots$ is a chain of open normal subgroups with trivial intersection. Let $M$ be a finitely generated $\pH$-module. 

{\bf (a)}\, \,  By  \cref{dimapp}, 
\begin{equation} \label{dimapp0} \dim_{\pH}M =\lim_{i\rar \infty} \frac{\dim_{\p} \left(\p [\HH/\HH_i]\tens{\p[[\HH]]} M \right)}{|\HH: \HH_i|}. \end{equation}
There is an injection of finite $p$-groups  $\HH/\HH_i\inc \FF/\NN_i$ and the image can be identified with $\HH\NN_i/\NN_i$. This yields to an isomorphism
\[(\FF/ \NN_i)/(\HH\NN_i/ \NN_i)\cong \FF/ \HH \NN_i. \]
Therefore, considering $\p[[\FF/\NN_i]]$ as a right $\pH$-module, we obtain that 
\[\p[\FF/\NN_i]\cong \p[\HH/\HH_i]^{|\FF: \HH \NN_i|}\]
as right $\pH$-modules. So 
\[\dim_{\p} \left(\p [\HH/\HH_i]\tens{\p[[\HH]]} M \right)=\frac{\dim_{\p} \left(\p [\FF/\NN_i]\tens{\p[[\HH]]} M \right)}{|\FF: \HH \NN_i|}.\]
Combining the latter equation; the application of \cref{dimapp} to the finitely generated $\pF$-module $\pF\tens{\pH} M$; and (\ref{dimapp0}), we conclude that
\begin{align*}
\dim_{\pF}\left(\pF\tens{\pH} M\right) =&\lim_{i\rar \infty} \frac{\dim_{\p} \left(\p [\FF/\NN_i]\tens{\p[[\FF]]} (\pF\tens{\pH} M)\right)}{|\FF: \NN_i|}\\
=& \lim_{i\rar \infty} \frac{\dim_{\p} \left(\p [\FF/\NN_i]\tens{\p[[\HH]]} M \right)}{|\FF: \NN_i|}\\
=& \lim_{i\rar \infty} \frac{\dim_{\p} \left(\p [\HH/\HH_i]\tens{\p[[\HH]]} M \right)}{|\HH: \HH_i|}\\
=&  \dim_{\pH} M. 
\end{align*}
On the other hand, since $\D_{\HH}$ is the division closure of $\pH$ in $\D_{\pF}$; then 
\begin{align*} 
\dim_{\D_{\pH}}\left(\D_{\pH}\tens{\pH} M\right)=&\dim_{\pH} M\\ 
=& \dim_{\pF}\left(\pF\tens{\pH} M\right)\\
=& \dim_{\D_{\pF}} \left(\D_{\pF}\tens{\pF} \left(\pF\tens{\pH} M\right) \right)\\
=& \dim_{\D_{\pF}} \left(\D_{\pF}\tens{\pH} M\right)\\
=& \dim_{\D_{\HH}} \left(\D_{\HH}\tens{\pH} M\right).
\end{align*}
By Cohn's criterion (\cref{Cohn1}), the division $\pH$-rings of fractions $\D_{\HH}$ and $\D_{\pH}$ are isomorphic. 

{\bf (b)}\, \, This was verfied during the proof of {\bf (a)}.

{\bf (c)}\, \, Given a $\pF$-module $M$, there are canonical isomorphisms of $\p$-vector spaces $\p[\HH/\HH_i]\tens{\pH} M\cong \p\tens{\p[[\HH_i]]} M$ and $\p[\FF/\NN_i]\tens{\pF} M\cong \p\tens{\p[[\NN_i]]} M$. Moreover, since $\HH$ is open then by compactness there exists $N$ such that, for all $i\geq N$, $\NN_i\sub \HH$ and hence $\NN_i=\HH_i$. Now using the latter observations and {\bf (a)} we get
\begin{align*}
    \dim_{\D_{\HH}} \D_{\pF} &= \dim_{\pH} \D_{\pF}\\
    &= \lim_{i\rar \infty} \frac{\dim_{\p} \left(\p [\HH/\NN_i]\tens{\p[[\HH]]} \D_{\pF}   \right)}{|\HH: \NN_i|}\\
    &= \lim_{i\rar \infty} \frac{\dim_{\p} \left(\p\tens{\p[[\NN_i]]} \D_{\pF} \right)}{|\HH: \NN_i|}
    \end{align*}
    \begin{align*}
    &= |\FF: \HH|\, \lim_{i\rar \infty} \frac{\dim_{\p} \left(\p\tens{\p[[\NN_i]]} \D_{\pF} \right)}{|\FF: \NN_i|}\\
    &= |\FF: \HH|\,  \frac{\dim_{\p} \left(\p[\FF/\NN_i]]\tens{\p[[\FF]]} \D_{\pF} \right)}{|\FF: \NN_i|}\\
    &=  |\FF: \HH|\,\dim_{\pF}  \D_{\pF}=|\FF: \HH|.
\end{align*}

{\bf (d)} \, \, We first prove this when we additionally have that $\HH$ is normal. In this case,  the automorphism-actions of $\FF$ on $\D_{\FF}$, induced by conjugation, leave the subring  $\pH$ invariant; and, then,  $\FF$ also normalizes  its division closure $\D_{\HH}$ in $\D_N$. As a consequence, the $\F_p$ subspace  $\FF \D_{\HH}$ spanned by all the products $\{f \cdot d: f\in \FF, a\in \D_{\HH}\}$ is a subring of $\D_{\pF}$. Moreover, it is clear that $\dim_{\D_H} (\FF \D_{\HH})\leq |\FF: \HH|<\infty$. This implies that $\FF \D_{\HH}$ is a domain and a finite dimensional $\D_{\HH}$-algebra; hence it is a division ring. Therefore, since $\FF \D_{\HH}$ contains $\pF$ ,  $\FF \D_{\HH}=\D_{\p \FF}$. We deduce that the canonical map of $(\pF, \D_{\HH})$-bimodules 
\[\al: \pF \tens{\pH} \D_{\HH}\lrar \D_{\pF}\]
can be read as a surjective homomorphism of right $\D_{\HH}$-modules of the same finite dimension by {\bf (b)} and {\bf (c)}. So $\al$ is an isomorphism. 

Having verified the property of $(d)$ for open normal $\HH\no \FF$, we want to extend this property to all open $\HH\lo \FF$ by induction on the index $|\FF: \HH|$. If $|\FF: \HH|\leq p$, then $\HH\no \FF$ satisfies the property. If $|\FF: \HH|>p$, then there exists an open normal $\HH\lo \HH_1\no \FF$ with $|\FF: \HH_1|=p$. By the inductive hypothesis, 
\[\D_{\pF}\cong \pF\tens{\pH}\D_{\HH}\cong \pF \tens{\p[[\HH_1]]} \left( \p[[\HH_1]]\tens{\pH} \D_{\HH}\right)\cong \pF\tens{\pH} \D_{\pH},\]
with canonical isomorphisms of $(\pF, \D_{\HH})$-modules in each step. The induction is complete. 
\end{proof}

From now on, we can identify $\D_{\HH}$ and $\D_{\p[[\HH]]}$ as $\pH$-rings without explicit mention of it.

\section{The ring $\D_{\p G}$}\label{DpGring}

There are some convenient features about group rings $KG$ of locally indicable groups. A result of Higman says that these rings are domains and that $(KG)^{\times}$ is $K^{*}\times G$.
We will talk about a specific type of division rings of fractions of these group rings, introduced by Hughes.

\begin{defi} Let $\Ga$ be a group and let  $K\Ga\inc \D$ be a division $K\Ga$-ring of fractions with the following property. For every finitely generated $S\unlhd T\leq \Ga$, such that $T=S\rtimes \Z$, where $\Z\leq T$ is generated by $t\in T$, we have that 
\begin{equation} \label{HFpair} \sum_{i\in \Z} K S \, t^i\end{equation}
is a direct sum of $K S$-submodules of $K T\sub \D$. 
Then we say that $\D$ is Hughes-free. 
\end{defi}
The following is proven in \cite{Hug70}.

\begin{thm}[Hughes] \label{HFunique} Let $\Ga$ be a locally indicable group. There exists at most one Hughes-free division $K \Ga$-ring of fractions up to $K \Ga$-isomorphism. 
\end{thm}

The following result is a consequence of \cite[Theorem 1.1.]{And19}.
\begin{thm}[Jaikin-Zapirain] \label{HFuniversal} Let $G$ be a residually-(torsion-free nilpotent) group and let $E$ be a division ring. Then there exists the universal division $EG$-ring of fractions $\D_{EG}$ and it is Hughes-free. 
\end{thm}

As a consequence of the last two results, we can prove that \textit{a priori} different $\p G$-division rings of fractions are $\p G$-isomorphic and universal, as the following corollary reflects. Let $G$ be an abstract subgroup of $\FF$, then $G$ is residually-(torsion-free nilpotent) by \cref{restorfree} and \cref{HFuniversal} applies. 

\begin{cor} \label{D12} Let $\FF$ be a free pro-$p$ group and let $G\leq \FF$ be a subgroup. 
\begin{enumerate}
    \item Let $\D_1$ be the division closure of $\p G$ in $\D_{\p [[\FF]]}$. 
    \item Take any group order $\leq$ on $G$ and let $\D_2$ be the division closure of $\p G$ in the corresponding $\p ((G))$. 
\end{enumerate}
Then $\D_1$ and $\D_2$ are both isomorphic, as $\p G$-rings, to $\D_{\p G}$. 
\end{cor}

We divide the verification of this corollary in two lemmas. The first lemma will correspond to the verification of the first part and relies on an improvement of the Hughes-freeness criterion (\cref{HFunique}) that is due to Sánchez \cite[Theorem 6.3]{San08}. An alternative approach of the first part of \cref{D12} based on Cohn's criterion (\cref{Cohn1})  is taken in \cite[Proposition 3.4]{And20}.

\begin{lemma} \label{D1HF} Let $\FF$ be a free pro-$p$ group of finite rank and let $G\leq \FF$ be a subgroup.  Let $\D_G$ be the division closure of $\p G$ in $\D_{\p [[\FF]]}$. Then $\D_G$ is Hughes-free. 
\end{lemma}
\begin{proof}
 Since $\p G$ has a Hughes-free division $\p G$-ring of fractions (\cref{HFuniversal}), then in order to conclude that $\D_G$ is isomorphic to $\D_{\p G}$ (and hence that it is Hughes-free),  it suffices to ensure by \cite[Remark 6.4(a)]{San08} the following property:  for all  finitely generated  $T\leq G$, there exists a subgroup $N_0\unlhd T$ with $T=\lan N_0, t\cong N_0\rtimes \Z$ such that (\ref{HFpair}) is satisfied; or, in other words, such that the canonical map of $(\p T, \D_{N_0})$-bimodules 
\begin{equation} \label{N0} \delta: \p[T]\tens{\p N_0} \D_{N_0} \inc \D_{T}\end{equation}
is injective.

Consider the pro-$p$ group $\ove{T}$. It is a closed subgroup of $\FF$. By \cref{subfree}, it is free of positive rank. As a consequence, there is a continuous epimorphism $f: \ove{T}\lrar \Z_p$. Its kernel $\KK$ is a closed subgroup of $\ove{T}$. We consider $K=\KK\bigcap T$. Notice that $T/K\cong \Z^m$ for some $m\geq 1$.

We want to prove that the canonical map 
of $(\p T, \D_{K})$-bimodules 
\begin{equation} \al: \p[T]\tens{\p K} \D_{K} \lrar \D_{T}\end{equation}
is injective.

Recall that $\p T$ is a free right $\p K$-module with a free basis given by a choice $\{c_i\}\sub T$ of representatives of the left-cosets of $K$ in $T$. As a consequence, $\p T\tens{\p K} \D_K$ is a free right $\D_K$-module of basis $\{c_i\otimes 1\}$. 
Similarly, let $\HH$ be an open subgroup of $\ove{T}$. Then $\p[[\ove{T}]]$ is a free (profinite) right $\p[[\HH]]$-module with basis given by a choice of left coset representatives $\{c_i'\}$ of $\HH$ in $\ove{T}$. In particular, viewing $\p[[\ove{T}]]\tens{\p[[\HH]]} \D_{\HH}$ as a right $\D_{\HH}$-module; the set  $\{c_i'\otimes 1\}$ is  $\D_{\HH}$-linearly independent.

Let us take some $0\neq x=\sum_{i=1}^n c_i \otimes d_i\in \p T\tens{\p K} \D_K$, where $n\geq 1$ and $0\neq d_i\in \D_K$. First of all, since $K=\KK\bigcap T$, the elements $\{c_i\}$ represent different left-cosets of $\KK$ in $\ove{T}$. For these fixed $c_1, \dots, c_n$; there must exist an open subgroup $\HH\leq \ove{T}$ such that the elements $\{c_i\}_{1}^n $ represent different left-cosets of $\HH$ in $\ove{T}$. 

There is a commutative diagram of canonical homomorphisms of $(\p T, \D_K)$-bimodules
\begin{equation*} 
\begin{tikzcd}
\p T\tens{\p K}D_{K} \arrow[r, rightarrow, "\al"] \arrow[d, "\beta"] & \D_T \arrow[d, hookrightarrow]\\
\p[[\ove{T}]] \tens{\p [[\HH]]} \D_{\HH} \ar[r, "\ga"] & \D_{\ove{T}},
\end{tikzcd}
\end{equation*}
where the map $\ga$ is an isomorphism by \cref{HFprop}. 

By the choice of $\HH$, $\be(x)=\sum_{i=1}^n c_i\otimes d_i\neq 0$. So, due to the commutativity of the above diagram, $\al(x)\neq 0$; hence $\al$ is injective. 

We select any $K\leq N\unlhd T$ such that $T/N\cong \Z$. 

There is a commutative diagram of canonical homomorphisms of $(\p T, \D_K)$-bimodules
\begin{equation*}
    \begin{tikzcd}
\p N\tens{\p K}D_{K} \arrow[r, rightarrow, "\al_1"] \arrow[d, "\beta_1"] & \D_N \arrow[d, hookrightarrow]\\
\p T \tens{\p K} \D_{K} \ar[r, "\al"] & \D_{T},
\end{tikzcd}
\end{equation*}
where $\al$ has been proven to be injective; and $\beta_1$ is injective because $\p N$ is a direct summand of $\p T$ as right $\p K$-modules. It follows that $\al_1$ is injective. 

We denote $\D_K N$ to be the $\p$-linear span of elements of the form $ab$ with $a\in \D_K$ and $b\in N$. This is a subring of $\D_{T}$ because the automorphism-actions of $N$ on $\D_N$, induced by conjugation, leave the subring  $\p K$ invariant; and, then,  $N$ also normalizes  its division closure $\D_K$ in $\D_N$.

We claim that $\D_K N$ is an Ore domain. 

Notice that there is a natural ring homomorphism from the crossed product $\p K* N/K\lrar \D_K N$. By the injectivity of $\al_1$, the previous is, in fact, an isomorphism. Since $\D_K N$ is a domain and  $\p K* N/K$ is Noetherian, due to $N/K$ being abelian; then $\D_K N$ is an Ore domain by Goldie's theorem (\cref{Goldie}). Since $\D_K N\inc \D_N $ is a division $\p N$-ring of fractions of the Ore domain $\D_K N$; then, by (\ref{Orecriterion}),  $\D_N$ is isomorphic to $\Qo (\D_K N)$. 

Observe that $\D_K N$ can also be seen as the isomorphic image of the $(\p N, \D_K)$-bimodule $\p N\tens{\p K} \D_K$ under  $\al_1$. There is a commutative diagram of $(\p N, \D_K)$-bimodules
\begin{equation*}
    \begin{tikzcd}
\p T \tens{\p K}\D_{K} \arrow[rd, hookrightarrow, "\al"] \arrow[d, "\cong"] & \\
\p T \tens{\p N} \D_{K}N \ar[r, "\al_2"] & \D_{T},
\end{tikzcd}
\end{equation*}
from which we deduce that $\al_2$ is injective. Recall that $\tens{\D_K N} \Qo (\D_K N)$ is an exact functor, so the induced map 
\[ \al_3:  \left(\p T \tens{\p N} \D_{K}N\right)\tens{\D_K N} \Qo (\D_K N) \lrar\D_{T}\tens{\D_K N} \Qo(\D_K N)\cong \D_T \]
is injective, too. We recall that $\Qo(\D_K N)\cong \D_N$ canonically; so we have an injective canonical homomorphism of $(\p T, \D_N)$-bimodules 
\begin{equation}  \al_4 : \p[T]\tens{\p N} \D_{N} \inc \D_{T}.\end{equation}
This proves that we can take $N_0=N$ on (\ref{N0}) and the proof is complete.
\end{proof}

\begin{remark} In the previous lemma, the use of the criterion of J. Sánchez to establish the Hughes-free property of $\D_G$ was crucial.  It might be tempting to directly prove this property for any pair $(S, T)$, with $S\n T\n G$ and  $T/S\cong \Z$;  by reducing it to the ``pro-$p$ Hughes-free property'' that was established in  \cref{HFprop}. However, one may not be able to separate $S$ and $T$ with pro-$p$ closed subgroups. In other words, it may happen that   $T\sub \ove{S}$. For instance, one can take $\FF=\Z_p$; $a=1+\sum_{k=0}^\infty p^{2k+2}$; $b=1+\sum_{k=0}^\infty p^{2k+1}$; $S=\lan a\ran$; and $T=\lan a, b\ran$. It is clear that $a, b$ are $\Z$-linearly independent in $\Z_p$, so $T/S\cong \Z$. Furthermore,  $\ove{S}=\ove{T}=\FF$.
\end{remark}

\begin{lemma} \label{D2HF} Let $G$ be an orderable group, with a choice of $\leq$, and let $\D_G$ be the division closure of $\p G$ in the division $\p G$-ring $\p ((G))$. Then for any $H\leq G$, the canonical map of $(\p G, \p H)$-bimodules 
\[\p G\tens{\p H} \D_H \lrar \D_G\]
is injective. 
\end{lemma}
\begin{proof}
The canonical map $\al: \p G\tens{\p H} \D_H\lrar \D_G$ is defined by 
$a\otimes b\mapsto ab$, for all $a\in \p G$, $b\in \D_H$. 
We take a collection of representatives $\{t_i: i\in I\}$ of the left cosets of $H$ in $G$. This way, $G$ is equal to the disjoint union $\sqcup_{i\in I} t_i H$ and $\p G$ admits, as right $\p H$-module, the decomposition $\p G=\oplus_{i\in I} t_i \p H.$ So $\p G\tens{\p H} \D_H$ admits, as right $\D_H$-vector space, the basis $\{t_i\otimes 1: i\in I\}$. We gave a description of the inverses of elements of $\p ((G))$ in (\ref{inverse}). In particular, it implies that any element of $\D_H$ is supported in $H$. Let us take $x=\sum_{i\in I} t_i\otimes x_i\in \p G\tens{\p H} \D_H$, whose image equals $\al(x)=\sum_{i\in I} t_i x_i$. If $x\neq 0$, then  $x_j\neq 0$ for some $j$. Write $x_j=\sum_{h\in H} \lambda_h \, h$, with $\lambda_h\in \p$ and some $\lambda_{h_0}\neq 0$.  Since $t_ix_i$ is supported in the coset $t_i H$, then the coefficient of the power series   $\sum_{i\in I} t_i x_i$ corresponding to $t_j h_0$ is equal to $\lambda_{h_0}\neq 0$. So $\al(x)\neq 0$ and $\al$ is injective. 
\end{proof}
With these two lemmas, \cref{D12} follows directly from previous results. 

\begin{proof}[Proof of \cref{D12}] By \cref{restorfree}, $G$ is residually-(torsion-free nilpotent) and hence, by \cref{HFuniversal}, there exists a universal division $\p G$-ring of fractions $\D_{\p G}$ and this is Hughes-free. 
Both $\D_1$ and $\D_2$ are Hughes-free division $\p G$-rings by lemmas \ref{D1HF} and \ref{D2HF}. We deduce by  \cref{HFunique} that these division $\p G$-rings  are isomorphic, since all of them are Hughes-free. 
\end{proof}

\begin{cor}[Strong Hughes-free property] \label{Hughes} Let $G$ be a residually-(torsion-free nilpotent) group. Then $\D_{\p G}$ is strongly Hughes-free, this is,   for any $H\leq G$, the canonical map of  $(\p G, \p H)$-bimodules
    \[\p G\tens{ \p H} \D_{\p H}\longrightarrow \D_{\p G}\]
    is injective.
\end{cor}
\begin{proof} By \cref{nilorder2}, $G$ is orderable. After taking an ordering $\leq$ on  $G$, we know from \cref{D2HF} that the division closure $D_G$ of $\p G$ in $\p ((G))$ is Hughes-free. By \cref{HFunique} and \cref{HFuniversal}, $\D_G$ and $\D_{\p G}$ are isomorphic $\p G$-rings. There is  a commutative diagram of canonical homomorphisms of $(\p G, \p H)$-bimodules \begin{equation*}
    \begin{tikzcd}
       \p G\tens{\p H} \D_H \ar[d, "\al"] \ar[r, "\cong"] & \p G\tens{ \p H} \D_{\p H} \ar[d, "\beta"]\\
       \D_G \ar[r, "\cong"] & \D_{\p G},
    \end{tikzcd}
\end{equation*}
where the horizontal arrows are isomorphisms. We already verified in \cref{D2HF} that $\al$ is injective. So $\beta$ must be injective, too. 
\end{proof}

We finish this section with a formula that allow us to estimate dimensions of $\p G$-modules by approximating it with the aid of finite-index subgroups of $G$. 

\begin{prop}[Lück-type approximation] \label{Luckmodp0} Let $G$ be an abstract finitely generated subgroup of a free pro-$p$ group of finite rank $\FF$. Let $\FF=\NN_1>\NN_2>\dots$ be a chain of open normal subgroups of $\FF$ with trivial intersection. Define $G_j=G\bigcap \NN_j$. Then for every finitely generated $\p G$-module $M$, 
\[\dim_{\p G} M=\lim_{i\rar \infty} \frac{\dim_{\p} \left(\p \tens{\p[G_i]} M\right)}{|G: G_i|}.\]
\end{prop}
\begin{proof}
Let $\D_G$ be the division closure of $\p G$ in $\D_{\pF}$. By \cref{D12}, $\D_G$ and $\D_{\p G}$ are isomorphic as $\p G$-rings. 
So 
\[\dim_{\p G} M=\dim_{\D_G} (\D_G \tens{\p G} M).\]
In addition, 
\[\dim_{\D_G} (\D_G \tens{\p G} M)= \dim_{\D_{\pF}} (\D_{\pF} \tens{\p G} M)=\dim_{\pF} (\pF\tens{\p G} M). \]
Since there is a canonical isomorphism $G/G_i\lrar \FF/\NN_i$, we have that 
\[\p \tens{\p G_i} M\cong \p[G/G_i]\tens{\p G} M \cong \p[\FF/\NN_i]\tens{\p G}M \cong  \p[\FF/\NN_i]\tens{\pF}\left(\pF\tens{\p G} M\right).\]
Thus, by \cref{dimapp}, 
\[\dim_{\p G} M=\lim_{i\rar \infty} \frac{\dim_{\p} \p[\FF/\NN_i]\tens{\pF}\left(\pF\tens{\p G} M\right)}{|\FF: \NN_i|} =\lim_{i\rar \infty} \frac{\dim_{\p} \left(\p \tens{\p[G_i]} M\right)}{|G: G_i|}\]
and the conclusion follows. 
\end{proof}

\section{Examples of $\D_{\p G}$-torsion-free modules }\label{torfreesection}

The importance of the existence of a universal division rings of fractions of a ring $R$, denoted by $\D_R$, is that it yields to a dimension theory of certain $R$-modules by extending them to scalars in $\D_R$. However, we need our $\D_R$-modules to be torsion-free for an statement of the type ``surjective homomorphisms between $R$-modules of the same dimension are isomorphisms''. 

We start by defining the concept of $\D_R$-torsion-free modules. We then give some general methods to construct them; and, finally, we give some explicit families of $\D_{\p G}$-torsion-free modules. 

\begin{defi}
Let $R\inc \D$ be an embedding of the ring $R$ into a division ring $\D$. Let $M$ be a $R$-module. We say that $M$ is $\D$-torsion-free if the canonical map $M\lrar \D\tens{R}M$ is injective. 
\end{defi}

The following provides a more flexible criterion for verifying whether a module is torsion-free. 
\begin{lemma}[Characterisations of torsion-freeness]\label{torfree} Let $M$ be a $\p G$-module. Then $M$ is $\D_{\p G}$-torsion-free if and only if there exists a $\D_{\p G}$-module $N$ and an injective homomorphism of $\p G$-modules $M\inc N$. 
\end{lemma}

We leave stated other two lemmas which are not hard to prove. We refer to \cite[Section 4.1]{And20} for details. 

\begin{lemma} \label{quotientdim}
Let $\D$ be a division $R$-ring and let $M$ be a $\D$-torsion-free $R$-module of finite $\D$-dimension. Let $L$ be a non-trivial $R$-submodule of $M$. Then $\dim_{\D}(M/L)<\dim_{\D} (M)$. Moreover, if $\dim_{\D} L=1$, then $\dim_{\D} (M/L)=\dim_{\D}M -1$. 
\end{lemma}

\begin{lemma}\label{sestorfree} Let $\D$ be a division $R$-ring and let $1\rar M_1\rar M_2\rar M_3\rar 0$ be an exact sequence of $R$-modules. Assume that 
\begin{enumerate}
    \item $M_1$ and $M_3$ are $\D$-torsion-free;
    \item $\dim_{\D}M_1$ and $\dim_{\D}M_3$ are finite;  and
    \item $\dim_{\D} M_1+\dim_{\D} M_3=\dim_{\D} M_2$. 
\end{enumerate}
Then $M_2$ is also $\D$-torsion-free.
\end{lemma}

We denote by $H\leq G$ subgroups of a free pro-p group $\FF$.
we alternatively classify what is it to be torsion free.

\begin{lemma} Let $M$ be a $\p H$-module. Then $M$ is $\D_{\p H}$-torsion-free if and only if $\p G\tens{\p H} M$ is $\D_{\p G}$-torsion-free. 
\end{lemma}
\begin{proof} Suppose that $M$ is $\D_{\p H}$-torsion-free. Then the canonical map 
\[M \lrar \D_{\p H}\tens{\p H} M\]
is injective. Notice that $\p G$ is a free right $\p H$-module. Thus the canonical map 
\[ \p G\tens{\p H} M \overset{ \al}{\longrightarrow} \p G\tens{\p H}\left( \D_{\p H}\tens{\p H} M\right)\]
is also injective. The right-hand side $\p G$-module is canonically isomorphic to 
\[\left(\p G\tens{\p H} \D_{\p H}\right) \tens{\p H} M.\]
We want to embed this $\p G$-module into a $\D_{\p G}$-module. To do so, we recall the strong Hughes-free property of \cref{Hughes}, which says that the canonical map of $(\p G, \p H)$-bimodules 
\[\p G\tens{ \p H} \D_{\p H}\longrightarrow \D_{\p G}\]
    is injective.
    Since $\D_{\p G}$ is a free right $D_{\p H}$-module and $D_{\p H}$ is a free right $\p H$-module, we deduce that both 
    $\p G\tens{ \p H} \D_{\p H}$ and $ \D_{\p G}$ are free right $\p H$-modules. It follows that the canonical map of $\p G$-modules 
    \[\left(\p G\tens{ \p H} \D_{\p H}\right)\tens{\p H} M \overset{ \be}{\longrightarrow} \D_{\p G}\tens{\p H} M.\]
    Composing $\al$ and $\be$, we construct an injective $\p G$-homomorphism between $\p G\tens{\p H} M$ and the $\D_{\p G}$-module $\D_{\p G}\tens{\p H} M$.

By \cref{torfree}, we  conclude  that $\p G\tens{\p H} M$ is  a $\D_{\p G}$-torsion-free module.
\end{proof}

\begin{prop} \label{Atorfree} Let $H\leq G\leq \FF$ and let $A$ be a maximal abelian subgroup of $H$. Then the $\F_pG$-module
\[\F_p I_H^G\, \big/\, \F_p I_A^G\]
is $\D_{\F_p G}$-torsion-free. 
\end{prop}

\begin{proof} Omitted. It is proven in \cite[Proposition 4.8]{And20}.
\end{proof}
This directly leads to the following consequence. 

\begin{cor} \label{utorfree} Let $u\in G$. Suppose that the cyclic group $\lan u\ran$ that $u$ generates is maximal abelian in $H$. Then the $\F_pG$-module 
\[\frac{\F_p I_H^G}{\F_p G (u-1)}\]
is $\D_{\F_p G}$-torsion-free. 
\end{cor}

The next proposition produces torsion-free modules that have the form of augmentation ideals of amalgamated products with cyclic amalgam. 

\begin{prop} \label{amaltor} Let $H_1$ and $H_2$ be two finitely generated subgroups of a finitely generated free pro-$p$ group $\FF$. Consider $A=H_1\bigcap H_2$ and suppose that $A$ is a maximal abelian subgroup of $H_1$. Let $G=\lan H_1, H_2\ran$ and let 
\[J=\{(x, -x): x\in \p I_A^{G}\}\leq \p I_{H_1}^G\oplus \p I_{H_2}^G.\]
Then the $\p G$-module 
\[M=\frac{\p I_{H_1}^G\oplus \p I_{H_2}^G}{J}\]
is $\D_{\p G}$-torsion-free and \[\dim_{\p G} M=\dim_{\p H_1} \p I_{H_1}+\dim_{\p H_2}\p I_{H_2}-1.\]
\end{prop}

Before giving the proof, we shall make an observation. If $A$ is finitely torsion-free abelian, as occurs in the previous proposition, then $A\cong \Z^n$ and the universal division $\p A$-ring of fractions is the field of fractions of the commutative domain $\p A$. Given a nontrivial finitely generated  ideal $I=(a_1, \dots, a_n)$ of  $ \p A$, we see that the $\Frac(\p A)$-vector space $\Frac(\p A)\tens{\p A} I$ can be generated by any $1\otimes a_i$ with $a_i\neq 0$. In fact, if $a_1\neq 0$, then $1\otimes a_i= a_1^{-1} a_1\otimes a_i = a_1^{-1} \otimes a_1 a_i=a_1^{-1}a_i\otimes a_1=a_1^{-1}a_i \,  (1\otimes a_1)$. This implies that 
\begin{equation}\label{dimA}
\dim_{\p A} I=\dim_{\Frac(\p A)} \Frac(\p A)\tens{\p A} I=1.
\end{equation}

\begin{proof}[Proof of \cref{amaltor}] We consider the $\p G$-submodule of $M$ defined by $L=(\p I_{A}^G\oplus \p I_{H_2}^G)/J$.

We want to apply \cref{sestorfree} to the short exact sequence $0\rar L\rar M\rar M/L\rar 0$ of $\p G$-modules. 

Since $L\cong \p I_{H_2}^G\leq \p G$, then it is $\D_{\p G}$-torsion-free and 
\[\dim_{\p G} L=\dim_{\p H_2} \p I_{H_2}.\] 
By the observation (\ref{dimA}), $1=\dim_{\p A} \p I_A=\dim_{\p G} \p I_A^G=\dim_{\p G} J$. So it is clear, by \cref{quotientdim}, that 
\[\dim_{\p G} M= \dim_{\p G} \p I_{H_1}^G+\dim_{\p G} \p I_{H_2}^G-1=\dim_{\p H_1}\p I_{H_1}+\dim_{\p H_2} \p I_{H_2}-1.\]

 On the other side, the quotient $M/L$ is isomorphic to $\p I_{H_1}^G/\p I_A^G$, which is $\D_{\p G}$-torsion-free by \cref{Atorfree}.
Again, by \cref{quotientdim},   \[\dim_{\p G} M/L=\dim_{\p G} \p I_{H_1}^G-\dim_{\p G}\p  I_A^G=\dim_{\p H_1} \p I_{H_1}-1.\] 

Therefore, the short exact sequence $0\rar L\rar M\rar M/L\rar 0$ of $\p G$-modules satisfies the requirements of \cref{sestorfree}; and the conclusion follows. 
\end{proof}

We now turn to the study of torsion-free modules that have the form of an augmentation ideal of a cyclic HNN extension. 

Our point is to prove that certain augmentation ideals, which may have the form 
\[R^m /R(u_1, \dots, u_m)\]
for some group ring $R$; 
are $\D_R$-torsion-free. When extending this $R$-module with coefficients in a bigger ring, some $u_i$ may become invertible. The following elementary lemma simply studies this scenario, which we shall encounter many times.

\begin{lemma} \label{freequo} Let $R$ be a unital ring and  $M$ be a $R$-module. Let $m_0\in M$ and let $u$ be a unit of $R$. Then there is an isomorphism of $R$-modules
\[\ga: \frac{M\oplus R}{(m_0, u)}\lrar M \]
given by 
\[\ga ([m, r])=m-ru^{-1}m_0,\]
with inverse 
\[\ga^{-1}(m)=[m, 0].\]
\end{lemma}

We can now state  the last proposition of this section. The proof is relatively long and, for convenience, it is divided in several intermediate lemmas. 

\begin{prop}\label{hnntor} Let $H\leq G$ be subgroups of $\FF$. Suppose that we can write $G=N\rtimes (t)=\lan N, t\ran$ for some $H\leq N\leq G$ and some $t\in G$. Let $u\in H$ be an element which generates a maximal abelian group $\lan u\ran $ in $H$ and suppose that $v=tut^{-1}\in H$. Then the $\p G$-module $M$ defined by 
\[M=\frac{\p I_H^G\oplus \p G}{\p G\,  (v-1-t(u-1), v-1)}\]
is $\D_{\p G}$-torsion-free. 
\end{prop}

\begin{proof} 
Since $G\cong N\rtimes (t)$, we can write $\p G\cong \F_p N[t^{\pm 1}, \sigma]$ as a skew-polynomial ring in one indeterminate $t$ with coefficients in $\p N$. Here $\sigma$ is an automorphism of the ring  $\D_{\p N}$ such that $\sigma (\p N)=\p N$, since it is induced from the automorphism of conjugation-by-$t$ in the subring $\p N$ of $\p G$. 

\begin{lemma} \label{isolem}
Let $R=\D_{\p N}[t^{\pm}, \sigma]$. Then $R$ is a (right) Ore domain and its (right) Ore division ring of fractions $Q_{\ore}(R)$ is isomorphic to $\D_{\p G}$.
\end{lemma}
To check this lemma, we start noting that $R$ is a domain. In addition, by the Hilbert's basis theorem, the ring $R$ is also Noetherian; and hence,  by Little Goldie's theorem \ref{Goldie}, $R$ is an Ore domain. We want to compare the resulting Ore division $R$-ring of fractions $R\inc \Qo (R)$ with $\D_{\p G}$. By \cref{D12}, we can identify $\D_{\p G}$ and $\D_{\p H}$ with the division closures of $\p G$ and $\p H$, respectively, in $\p ((G))$. In this setting, there is a natural monomorphism $\D_{\p N}\inc \D_{\p G}$. We want to extend this to a monomorphism of rings $R\inc \D_{\p G}$. Observe that, as $(\p G, \p H)$-bimodules, there is a canonical injection
\[R= \D_{\p N}[t^{\pm}, \sigma]\cong \p G\tens{\p N} \D_{\p N}\inc \D_{\p G}\]
according to \cref{D2HF}. In addition, we can observe that this map is a ring homomorphism. So we have an ring monomorphism $\iota : R\inc \D_{\p G}$. Since $\D_{\p G}$ is a division ring, we have, by the universal property (\ref{uni3}) of the embedding $\iota_{\ore}: R\inc \Qo (R)$, a commutative diagram of ring homomorphisms 
\begin{equation} 
\begin{tikzcd}
R \arrow[r, rightarrow, "\iota_{\ore}"] \arrow[dr, rightarrow, "\iota"] & \Qo (R) \arrow[d, dashrightarrow, "\hat{\iota}"]\\& D_{\p G},
\end{tikzcd}
\end{equation}
where $\iota$ and $\iota_{\ore}$ are injective. 

Moreover, $R$ contains $\p N[t^{\pm 1}, \sigma]=\p G$, so $\iota: R\inc \D_{\p G}$ is a division $R$-ring of fractions. Arguing as in \cref{Orecriterion},  $\hat{\iota}: \Qo (R)\lrar \D_{\p G}$ is an isomorphism of  $R$-rings. So \cref{isolem} is proven. 

There is a division $R$-ring $R\inc \D_{\p N} ((t^{\pm}, \sigma))$, where the latter is the ring of Skew Laurent series. 
The universal property of $\Qo (R)$ and the isomorphism $ \Qo(R)\cong \D_{\p G}$ of \cref{isolem} induces a commutative diagram of rings
\begin{equation} 
\begin{tikzcd}
R \arrow[r, hookrightarrow, "\iota_{\ore}"] \arrow[dr, hookrightarrow] & \D_{\p G} \arrow[d, dashrightarrow]\\&  D_{\p N}((t^{\pm}, \sigma)).
\end{tikzcd}
\end{equation}

In particular, we can upgrade this diagram to one of the following form

\begin{equation} \label{skewdia}
    \begin{tikzcd}
     \p G \ar[d, hookrightarrow] \ar[r, hookrightarrow] & \D_{\p N}[t^{\pm 1}, \sigma] \ar[d, hookrightarrow] \\
     \D_{\p G}  \ar[r, hookrightarrow] &  \D_{\p N}((t^{\pm 1}, \sigma)),
    \end{tikzcd}
\end{equation}
where all the arrows are injective. This induces the following commutative diagram of $\p G$-modules

\begin{equation} \label{skewdia0}
    \begin{tikzcd}
     \p I_H^G \ar[d, "\ga"] \ar[r, "\al"] & \D_{\p N}[t^{\pm 1}, \sigma]\tens{\p G} \p I_H^G \ar[d, "\be"] \\
     \D_{\p G} \tens{\p G} \p I_H^G \ar[r, "\delta"] &  \D_{\p N}((t^{\pm 1}, \sigma)) \tens{\p G} \p  I_H^G\\
    \end{tikzcd}
\end{equation}
\begin{lem} \label{lemskewdia0} The map $\beta$ of the diagram (\ref{skewdia0}) is injective. 
\end{lem}

We claim that both $\ga$ and $\delta$ are injective. On the one hand, 
the $ \p G$-module $\p I_H^G$ is a submodule of $\p G$, so it is $\D_{\p G}$-torsion-free and $\ga$ is injective. Furthermore,   $\D_{\p N}((t^{\pm 1}, \sigma))$ is a free right $\D_{\p G}$-module, since $\D_{\p G}$ is a division ring; so $\D_{\p G}$ is a direct summand of $\D_{\p N}((t^{\pm 1}, \sigma))$ as a right $\p G$-module. Hence $\delta$ is injective. Since the diagram (\cref{skewdia0}) is commutative, $\beta$ must be injective and \cref{lemskewdia0} is proven.

We will now view this injective map $\beta$ in another diagram. 

The commutative diagram of rings (\ref{skewdia})  induces the following commutative diagram of $\p  G$-modules

\begin{equation} \label{skewdia2}
    \begin{tikzcd}
     M\ar[d, "\ga'"] \ar[r, "\al'"] & \D_{\p N}[t^{\pm 1}, \sigma]\tens{\p G}  M \ar[d, "\be'"] \\
     \D_{\p G} \tens{\p G} M \ar[r, "\delta'"] &  \D_{\p N}((t^{\pm 1}, \sigma)) \tens{\p G} M  \\
    \end{tikzcd}
\end{equation}
Now let $R$ be either $ \D_{\p N}[t^{\pm 1}, \sigma]$ or   $\D_{\p N}((t^{\pm 1}, \sigma)) .$ Observe that 
\begin{equation*}
\begin{tikzcd}
     R\tens{\p G} M\ar[r, "\cong"] & \frac{\left(R\tens{\p G} \p I_{H}^G\right)\oplus R}{R\, (v-1-t(u-1), v-1)}\ar[r, "\ga \, \cong"] & R\tens{\p G} \p I_{H}^G,\\
\end{tikzcd}
\end{equation*} 
where both arrows are isomorphisms of $\F_p$-modules and $\ga$ is given by \cref{freequo}, since $v-1$ is invertible in $\D_{\p N}$ and so also in $R$. We now integrate these isomorphisms $\ga$ into the commutative diagram of (\ref{skewdia2}) as follows.

\begin{equation} \label{skewdia3}
    \begin{tikzcd}
     M\ar[d, "\ga'"] \ar[r, "\al'"] & \D_{\p N}[t^{\pm 1}, \sigma]\tens{\p G}  M \ar[d, "\be'"] \ar[r, "\ga_1 \, \cong"] &  \D_{\p N}[t^{\pm 1}, \sigma] \tens{\p G} \p I_H^G \ar[d, "\beta"]\\
     \D_{\p G} \tens{\p G} M \ar[r, "\delta'"] &  \D_{\p N}((t^{\pm 1}, \sigma)) \tens{\p G} M \ar[r, "\ga_2 \, \cong"] &  \D_{\p N}((t^{\pm 1}, \sigma)) \tens{\p G} \p I_H^G \\
    \end{tikzcd}
\end{equation}
We already know, from \cref{lemskewdia0}, that $\beta$ is injective. So $\be'$ must be injective since  $\ga_1$ and $\ga_2$ are both isomorphisms. 

We claim that $\al'$ is also injective. If we proved this, then, by inspecting again the commutative diagram (\ref{skewdia3}),  $\gamma'$ would also  be injective and hence we could conclude that $M$ is $\D_{\p G}$-torsion-free, as we want. 

So it rests to prove the following lemma. 
\begin{lemma} \label{al}
The map $\al'$ of the diagram (\ref{skewdia3}) is injective. 
\end{lemma} 

We will prepare a bit the ground before checking this lemma. 

By assumption, $\lan u\ran$ is a maximal abelian subgroup of $H$. By \cref{utorfree}, the $\p N$-module 
\[M_0=\frac{\p I_H^N}{\p N (u-1)}\]
is $\D_{\p N}$-torsion-free. This implies that the canonical map 
\[M_0\lrar \D_{\F_p N} \tens{\p N} M_0\]
is injective. This translates into the following equality of subsets of $\D_{\p N} \tens{\p N} \p I_H^N$,
\begin{equation} \label{int1}  1\tens{\p N} \p I_H^N\, \, \bigcap\, \,  \D_{\p N}\tens{\p N} (u-1)=1\tens{\p N} \p N (u-1).  \end{equation}

From this, we want to prove the following lemma. 
We denote $H^{t^n}=t^n Ht^{-n}$, which are also subgroups of $N$ having $\lan u^{t^n}\ran $ as maximal abelian subgroup. 

\begin{lemma}  For all $n\in \Z$, we have the following equality of subsets of $\D_{\p N}\tens{\p N}  \p I_G$,
\begin{equation} \label{int3} 1\tens{\p N} \p t^n I_H^N\, \, \bigcap \, \, \D_{\p N}\tens{\p N} t^n(u-1)=1\tens{\p N}\p t^n N(u-1).\end{equation}
\end{lemma}
The same way we had the equality (\ref{int1}), we can derive, for the same reasons, 
\begin{equation} \label{int11}  1\tens{\p N} \p I_{H^{t^{-n}}}^N\, \, \bigcap\, \,  \D_{\p N}\tens{\p N} (u^{t^{-n}}-1)=1\tens{\p N} \p N (u^{t^{-n}}-1).  \end{equation}
Notice that $\p I_G$ has a right $\lan t\ran$-module structure  by multiplication. This induces a right $\lan t\ran$-module structure on $\D_{\p N}\tens{\p N}  \p I_G$. Since $N\n G$, then $t$ normalises $N$. We have the equations of subsets of $\D_{\p N}\tens{\p N}  \p I_G$:
\[\left(1\tens{\p N} \p N (u-1)\right)\, t^m= 1\otimes \p \, t^m N (u^{t^{-m}}-1),\]
and 
\[\left(1\tens{\p N} \p I_H^N \right)\, t^{m} = 1\otimes \p \, t^m  I_{H^{t^{-m}}} (u^{t^{-m}}-1).\]
As a consequence, applying the multiplication-by-$t^n$ automorphism of $\D_{\p N}\tens{\p N}  \p I_G$ to the equation  (\ref{int11});  we get (\ref{int3}).

Notice the following decomposition of $\p N$-modules
\[\p I_H^G=\oplus_{n\in \Z}\, \, \p \, t^n I_H^N,\]
which yields to the following decomposition of $\p$-vector spaces 
\[\D_{\p N}\tens{\p N} \p I_H^G\cong \oplus_{n\in \Z} \, \D_{\p N} \tens{\p N}\, \p\, t^n  I_H^N.\]

\begin{lemma} 
We have the following equation of subsets of $\D_{\p N}\tens{\p N} \p I_H^G$,
\begin{equation} \label{int2} 1\tens{\p N} \p I_H^G\, \, \bigcap \, \, \D_{\p N} \tens{\p N} \p G(v-1-t(u-1))=1\tens{\p N} \p G(v-1-t(u-1)). \end{equation}
\end{lemma}

Let us take an element $w$ that belongs to the left-hand side. This element will have the form 
\[w=\sum_{k=n_1}^{n_2}c_k\otimes t^k(v-1-t(u-1)),\, \, \mbox{for some $c_k\in \D_{\p N}$,}\]
and will also belong to $1\tens{\p N} \p I_H^G$.  We rewrite 
\[w= c_{n_1}\otimes t^{n_1}(v-1)+\sum_{k=n_1+1}^{n_2} \left(c_k\otimes t^k(v-1) -c_{k-1}\otimes t^{k}(u-1)\right) +c_{n_2}\otimes t^{n_2+1}(u-1).\]
Since $w\in 1\tens{\p N} \p I_H^G $, we can look at the highest power $t^{n_2+1}$ to deduce that 
\[c_{n_2}\otimes t^{n_2+1}(u-1)\in  1\tens{\p N} \p t^{n_2+1} I_H^N\, \, \bigcap \, \, \D_{\p N}\tens{\p N} t^{n_2+1}(u-1).\]
By (\ref{int3}), this implies that 
\[c_{n_2}\otimes t^{n_2+1}(u-1)\in 1\tens{\p N}\p t^{n_2+1} N(u-1),\]
so $c_{n_2}\in \p N$. 
Let $n_1+1\leq k\leq n_2$. Inspecting again the expression of $w$ at the component with power $t^k$, we have that \[c_k\otimes t^k(v-1) -c_{k-1}\otimes t^{k}(u-1)\in  1\tens{\p N} \p t^k I_H^N.\]
If we knew that  $c_k\in \p N$, then it would follow that 
\[ c_{k-1}\otimes t^{k}(u-1)\in 1\tens{\p N} \p t^k I_H^N\, \, \bigcap \D_{\p N}\tens{\p N} t^{k}(u-1).\]
By (\ref{int3}), this means that 
\[c_{k-1}\otimes t^{k}(u-1) \in 1\tens{\p N}\p t^n N(u-1),\]
and this implies that $c_{k-1}\in \F_p N$. 

We have proven that if $c_k\in \F_p N$, for $n_1<k\leq n_2$, then $c_{k-1}\in \p N$. Since we also know that $c_{n_2}\in \p N$, an inductive argument gives that  $c_k\in \p N$ for every $k$, meaning that  \[w\in 1\tens{\p N} \p N[t^{\pm}, \sigma]\,  (v-1-t(u-1))= 1\tens{\p N} \p G(v-1-t(u-1)).\]
This proves that 
\[ 1\tens{\p N} \p I_H^G\, \, \bigcap \, \, \D_{\p N} \tens{\p N} \p G(v-1-t(u-1))\sub 1\tens{\p N} \p G(v-1-t(u-1)),\]
one inclusion of (\ref{int2}). The reverse inclusion is trivial, so equation (\ref{int2}) is proven. 

Observe that there is a canonical isomorphism of $\p N$-modules
\[\D_{\p N}[t^{\pm}, \sigma] \cong \D_{\p N} \tens{\p N} \p N[t^{\pm}, \sigma]=  \D_{p N} \tens{\p N} \p G.  \]
This extends to a canonical isomorphism of $\p N$-modules
\[\psi: \frac{\D_{\p N}[t^{\pm}, \sigma] \tens{\p G} \p I_H^G }{\D_{\p N}[t^{\pm}, \sigma]\tens{\p G} \p G(v-1-t(u-1))}\lrar \frac{\D_{\p N} \tens{\p N} \p I_H^G}{ \D_{\p N}\tens{\p N} \p G(v-1-t(u-1))}. \]
Lastly, there is a commutative triangle of canonical $\p N$-homomorphisms 
\begin{equation*}
\begin{tikzcd}
    \frac{\p I_H^G}{\p G (v-1-t(u-1))} \ar[d, "\al''"] \ar[dr, "\eta"] & \\
     \frac{\D_{\p N}[t^{\pm}, \sigma] \tens{\p G} \p I_H^G }{\D_{\p N}[t^{\pm}, \sigma]\tens{\p G} \p G(v-1-t(u-1))} \ar[r, "\psi"] & \frac{\D_{\p N} \tens{\p N} \p I_H^G}{ \D_{\p N}\tens{\p N} \p G(v-1-t(u-1))},\\
\end{tikzcd}
\end{equation*}
The canonical map $\eta$ is injective due to (\ref{int2}). Since $\psi$ is an isomorphism, this implies that $\al''$ is injective. 
From the injectivity of $\al''$, the injectivity of $\al'$ follows directly. In fact, let $(x, y)\in \p I_H^G\oplus \p G$ belong to the kernel of $\al$. Then $1\otimes (x, y)=c \, (v-1-t(u-1), v-1)$ for some $c\in \D_{\p G}$. Notice that 
\[1\otimes x=c\otimes (v-1-t(u-1))\in 1\tens{\p G} \p I_H^G\bigcap \D_{\p N}[t^{\pm}, \sigma] \tens{\p G} (v-1-t(u-1)).\]
From the injectivity of $\al''$, this implies that 
\[c\otimes (v-1-t(u-1))\in 1\tens{\p G} \p G(v-1-t(u-1)),  \]
so $c\in \p G$  and then $(x, y)\in \p G(v-1-t(u-1), v-1).$
Thus $\al'$ is injective and  \cref{al} is demonstrated. The proof of \cref{hnntor} is complete.
\end{proof}

\chapter{Homology and {$L^2$}-Betti numbers}\label{L2section}


$L^2$-Betti numbers have played an important role in the solution of many problems in group theory. They were originally introduced by M. Atiyah in the context of Riemannian manifolds in 1974. J. Dodziuk   extended the notion of $L^2$-Betti numbers to free cocompact
actions of discrete groups $G$ on CW-complexes $X$ in 1977. Later on, in 1998, W. Lück extended this notion to arbitrary $G$-CW-complexes.

We shall not discuss here how these analytical and topological notions are defined. We are more interested in the algebraic approach that P. Linnell \cite{Lin93}  introduced in 1993. If $G$ is residually-(torsion-free nilpotent), we can combine Linnell's definition and a description of $\D_{\p G}$ that is based on the work of A. Jaikin-Zapirain and D. López-Álvarez \cite{And192} to redefine the $L^2$-Betti numbers as follows. 

Let $k\geq 0$. We define the {\bf $k$-th $L^2$-Betti number of $G$}, denoted $b_k^{(2)}(G)$,  by 
\[b_k^{(2)}(G)=\dim_{\D_{\Q G}} H_k (G; \D_{\Q G}).\]
By analogy, we can define in characteristic $p$ the {\bf $k$-th mod $p$ $L^2$-Betti number of $G$} as 
\[\beta_k^{\text{mod} p}(G)=\dim_{\D_{\p G}} H_k (G; \D_{\p G}).\]

We are most interested in $\bp (G)$ and we will later see how these numbers are related to the computations of dimensions of the previous chapter. 

First, we discuss how the most elementary version of Betti numbers can be used to study discrete groups. 

\begin{defi} \label{Betti} Given a finitely generated group $G$, its first Betti number is 
\[b_1(G)=\dim_{\Q} \left( G/[G, G]\right)\tens{\Z}\Q .\]
\end{defi}

If $G$ is finitely generated, we understand very well the dimension of the previous $\Q$-vector space. Since $G/[G, G]$ is a finitely generated abelian group, it must have the form $\Z^r\oplus T$, where $T$ is finite (the torsion subgroup). Tensoring by $\Q$ forgets about the torsion part of the $\Z$-module $G/[G, G]$, so 
\[ \left( G/[G, G]\right)\tens{\Z}\Q\cong \Q^r.\]
In particular,  $r=b_1(G)$. Interestingly, one can find $b_1(G)$ just by looking at certain finite quotients of $G$. 

\begin{lem} \label{b1} Let $G$ be finitely generated and let $p$ be a prime. Then $b_1(G)$ is the biggest $m$ such that $G$ surjects $(C_{p^k})^m$ for every $k\geq 1$. 
\end{lem}

We will see how one can exploit this observation to estimate Betti numbers of dense subgroups of free pro-$p$ groups. Later, this will have consequences for the first $L^2$-Betti number of a parafree group.

\begin{lem} \label{b1estimation}
Let  $G$ and $G'$ be finitely generated and let $p$ be a prime. If there is a group homomorphism $h: G\lrar G_{\hat{p}}'$ of dense image, then $b_1(G)\geq b_1(G')$. 
\end{lem}
\begin{proof} Let $m=b_1(G')$. Recall that, given a finite $p$-group, there is a correspondence (\ref{correspondence}) between the epimorphisms $G'\lrar P$ and the epimorphisms $G_{\hp}'\lrar P$; given by $f\mapsto f_{\hp}$. Similarly, given any epimorphism $g: G'_{\hp}\lrar P$, there is an epimorphism $g\circ f: G\lrar P$. This is because $P$ is finite, $g$ is continuous and $f(G)$ is dense in $G'_{\hp}.$

Therefore, every finite $p$-quotient of $G'$ is also a finite $p$-quotient of $G'$; hence $b_1(G)\geq b_1(G')$ by \cref{b1}. 
\end{proof}

In particular, the first Betti  number is a pro-$p$ invariant. 

\begin{cor}
If $G$ and $G'$ are finitely generated and $G_{\hat{p}}\cong G_{\hat{p}}'$, then $b_1(G)=b_1(G')$. 
\end{cor}

\section{Estimations and computations of $\b (G)$}
 In order to compute the $L^2$-Betti numbers of a finitely presented residually finite group, it suffices to have information about the Betti numbers of its finite-index subgroups. This is precised by the following important result of Lück \cite{Luc94}.

\begin{thm}[Lück's Approximation Theorem] \label{Luck} Let  $G$ be finitely presented and let $G=N_1>N_2>\dots >N_m>\dots $ be a sequence of finite-index normal subgroups  with $\bigcap_m N_m=1$. Then 
\[\lim_{m\rar \infty} \frac{b_1(N_m)}{|G:N_m|}= b_1^{(2)}(G).\]
\end{thm}

A favourable family of groups for which we can estimate or compute their $L^2$-Betti numbers are residually finite groups for which we have some structural information of their finite-index subgroups in relation to their index. This is the case of fundamental groups of some topological spaces, such as free groups and surface groups.

\begin{eg}[First $L^2$-Betti numbers of surface groups] Finitely generated free groups $F_n$ have $\b (F_n)=n-1$. If $S$ is a closed surface of genus $g\geq 1$, then \[\b(\pi_1(S))=-\chi (S)=2g-2.\] 
We will use Lück's approximation theorem to verify this. The case of $F_n$ follows from Schreier's index-rank formula (\cref{freeindex}) and the fact that a $b_1(F_k)=k$ for every $k\geq 1$. Let  $F_n=N_1>N_2>\dots >N_m>\dots $ be a sequence of finite-index normal subgroups  with $\bigcap_m N_m=1$. The existence of this chain is ensured by because $F_n$ is residually finite. By \cref{Luck}, 
\[\b(F_n)=\lim_{m\rar \infty}\frac{b_1(N_m)}{|F_n:N_m|}=\lim_{m\rar \infty} \frac{1+|F_n:N_m|(n-1)}{|F_n:N_m|}=n-1. \]
Let $S$ be a closed orientable surface. As we verified in \cref{surfaceab}, if $S$ is orientable, then $b_1(\pi_1(S))=2-\chi(S)$, and, if $S$ is non-orientable, then $b_1(\pi_1(S))=1-\chi(S)/2$. The computation of their $\b$ is analogous and we only precise the case when $S$ is orientable. Let  $\pi_1(S)=N_1>N_2>\dots >N_m>\dots $ be a sequence of finite-index normal subgroups  with $\bigcap_m N_m=1$. The existence of this chain is ensured by \cref{surfaceresfinite}. By \cref{surfaceindex}, $b_1(N_m)=2-|\pi_1(S): N_m|\chi (S)$, so 
\[\b(\pi_1(S))=\lim_{m\rar \infty}\frac{b_1(N_m)}{|\pi_1(S):N_m|}=\lim_{m\rar \infty} \frac{2-|\pi_1(S): N_m|\chi (S)}{|\pi_1(S):N_m|}=-\chi(S)=2g-2.\]
We have also used that $|\pi_1(S): N_m|\rar \infty$ because $\pi_1(S)$ is infinite. 
\end{eg}

In the previous example we saw the relation between the Euler characteristic of a surface and the first $L^2$-Betti number of its fundamental group. This relation is extended to one-relator groups and surface-plus-one-relation groups by Dicks and Linnell \cite{Dic07}. The Euler characteristic also appears in the computation of the first $L^2$-Betti numbers of compact 3-manifold groups \cite{Lot95}. As opposed to the case of surface groups, their first $L^2$-Betti numbers are not typically nonzero. 

In the case of finitely generated groups that a priori are not finitely presented, one can still estimate its first $L^2$-Betti number in a similar way as Lück's approximation theorem. 
The following is a result of Lück and Osin \cite{Luc11}.

\begin{thm} \label{Bettiest} Let $G$ be finitely generated and let $G=N_1>N_2>\dots >N_m>\dots $ be a sequence of finite-index normal subgroups  with $\bigcap_m N_m=1$. Then
\[\limsup_{m\rar \infty} \frac{b_1(N_m)}{|G:N_m|}\leq  b_1^{(2)}(G).\]
\end{thm}

In order to understand finitely generated normal subgroups of parafree groups, we will also use the following theorem from of Gaborian \cite[Theorem 6.8]{Gab02}.
\begin{thm} \label{shortBetti} Suppose that 
\[1 \rar N\rar G \rar G' \rar 1\]
is an exact sequence of groups where $N$ and $G'$ are infinite. If $b_1^{(2)}(N)<\infty$, then $b_1^{(2)}(G)=0$. 
\end{thm}

The following result is due to Bridson and Reid \cite[Proposition 7.5]{Bri14}. 

\begin{prop} \label{Bettiest1}
Let $G$ be a finitely generated group and let $F$ be a finitely presented group that is residually-$p$ for some prime $p$. Suppose that there is an embedding $G\inc F_{\hat{p}}$ of dense image. Then $\b (G)\geq \b(F)$.
\end{prop}

This is proven by estimating $\b(F)$ and $\b(G)$ with Lück approximations' of theorems \ref{Luck} and \ref{Bettiest}, respectively. 

In particular, the first $L^2$-Betti number is a pro-$p$ invariant among finitely presented residually-$p$ groups. 

\begin{cor}
Let $G$ and $G'$ be finitely presented residually-$p$ groups such that $G_{\hat{p}}\cong G'_{\hat{p}}$. Then $\b(G)=\b(G')$. 
\end{cor}

Before proving \cref{Bettiest1}, we make a few observations about the pro-$p$ topology. 

In general, given abstract groups $\La\leq \Ga$; it may not happen that the canonical map  $\La_{\hp}\lrar \ove{\La}\sub \Ga_{\hp}$ is an isomorphism. We will implicitly search for this favourable scenario during \cref{pro-$p$ embeddings}. For the moment, we can give a simpler sufficient criterion to ensure that it holds. 

\begin{lemma} \label{pro-p subspaces} Let $ \La\inc \Ga$ be an inclusion of abstract groups, where $\La$ has finite $p$-power index in $\Ga$. Then the canonical map  $ {\La}_{\hp}\lrar \Ga_{\hp}$ is injective. In other words, by denoting $\iota_{\hp}: \Ga\lrar \Ga_{\hp}$, the canonical map ${\La}_{\hp}\lrar \ove{\iota_{\hp}(\La)}\leq \Ga_{\hp}$ is an isomorphism of pro-$p$ groups. 
\end{lemma}
\begin{proof}  By \cite[Proposition 3.2.6]{Rib00}, the conclusion of this lemma holds exactly when the pro-$p$ topology of $\Ga$ induces the pro-$p$ topology on its subspace $\La$. In addition, by  \cite[Lemma 3.1.4(a)]{Rib00}, the previous condition is ensured if $\La$ has $p$-power index in $\Ga$. 
\end{proof}

This observation can be used to compare first $L^2$-Betti numbers of abstract groups with the rank of their pro-$p$ completions.

\begin{proof}[Proof of \cref{Bettiest1}] Since $F$ is finitely generated, for every integer $D\geq 1$ there are only a finite number of subgroups of $F$ of index at most $D$. Let $F(d)$ denote the intersection of all normal subgroups of $F$ with $p$-power index at most $p^d$. Alternatively, $F(d)$ is the finite intersection of all kernels of surjective maps $F\lrar P$ for any finite $p$-group of order at most $p^d$. By the correspondence discussed in \cref{correspondence}, it is clear that the closure $\ove{F(d)}$ of $F(d)$ in $F_{\hp}$ is equal to the finite intersection of all the open normal subgroups of $F_{\hp}$ (here we have implicitly  used the fact that $\ove{X\bigcap Y}=\ove{X}\bigcap \ove{Y}$ for $p$-power index subgroups of $F$). In particular, each $\ove{F(d)}$ has finite index in $F$ and $\bigcap_d \ove{F(d)}=\{1\}$. 

We introduce $G(d)=G\bigcap \ove{F(d)}$, which does also verify that $\bigcap_d G(d)=\{1\}$. The canonical map $G\lrar F_{\hp}/\ove{F(d)}$ has dense image so it is surjective. By the choice of $G(d)$, it factors through an isomorphism $G/G(d)\lrar F_{\hp}/\ove{F(d)}$. On the other side, the canonical maps $F/F(d)\lrar (F/F(d))_{\hp}\lrar F_{\hp}/\ove{F(d)}$ are isomorphisms since $F/F(d)$ is a finite $p$-groups. This results in the equalities 
\[|G: G(d)|=|F_{\hp}: \ove{F(d)}|=|F: F(d)|.\]
By \cref{pro-p subspaces}, the canonical map $F(d)_{\hp}\lrar \ove{F(d)}\sub F_{\hp}$ is an isomorphism of pro-$p$ groups. Since $G$ is dense in $F_{\hp}$, then $G(d)$ is dense in $\ove{F(d)}\cong F(d)_{\hp}$, so $b_1(G(d))\geq b_1(F(d))$ for all $d\geq 1$ by \cref{b1estimation}. 

Lastly, we can compare $\b(F)$ and $\b(G)$ by using the sequences $\{F(d)\}_d$ and $\{G(d)\}_d$ by means of \cref{Luck} and \cref{Bettiest1}. Notice that 
\[\b(G)\geq \limsup_{d\rar \infty} \frac{b_1(G(d))}{|G: G(d)|}\geq \limsup_{d\rar \infty} \frac{b_1(F(d))}{|F: F(d)|}=\b(F).\qedhere\]
\end{proof}

We mentioned in the introduction that $\b (G)=\dim_{\D_{\Q G}} H_1 (G; \D_{\Q G}) $. The group ring $\Q G$ of zero characteristic admits some analytical tools that are not available in the regime of positive characteristic. However, for our purposes, we prefer \textit{a priori} to work with mod-$p$ $L^2$ Betti numbers $\bp$ instead of $\b$. Our group $G$ of study is a subgroup $G\leq \FF$ of a finitely generated free pro-$p$ group. As such, the ring $\D_{\p G}$ has other favourable properties in comparison to $\D_{\Q G}$; since it lies inside $\D_{\pF}$. Still, for parafree groups,  $L^2$-Betti numbers do not depend on the characteristic.

\section{The mod-$p$ version of $L^2$-Betti numbers $\bp (G)$}

Here we briefly describe the mod-$p$ analogue of the approximation results of $L^2$-Betti numbers that were described in the previous section. 

\begin{prop}[Mod-$p$ analogue of the Lück approximation theorem] \label{Luckmodp} Let $\FF$ be a finitely generated free pro-$p$ group and let $G$ be a finitely generated subgroup of $\FF$ of type $FP_k$ for some $k\geq 1$. Let $\FF=N_1>N_2>\dots$ be a chain of open normal subgroups of $\FF$ with trivial intersection. Define $G_j=G\bigcap N_j$. Then 
\[\beta_{k}^{mod  \,p} (G)=\lim_{j\rar \infty} \frac{\dim_{\F_p} H_k(G_j; \F_p)}{|G: G_j|}.\]
In particular, if $G$ is finitely generated, then
\[\bp (G)=\dim_{\p G} \p I_G-1.\]
\end{prop}
\begin{proof} If $G$ is of type $FP_k$, there exists a resolution of trivial  $\p G$-module  $\p$ of the form 
\[0\lrar R_k\lrar \p[G]^{n_{k-1}}\xrightarrow{f}\lrar \cdots \lrar \p[G]^{n_0}\lrar \p\lrar 0,\]
where $n_i$ are non-negative integers and $R_k$ is a finitely generated $\p G$-module. In order to compute the homology groups $H_k(G; *)$, we isolate the short exact sequence 
\[0\lrar R_k\lrar \p[G]^{n_{k-1}}\xrightarrow{f} \im f\lrar 0.\]
In these terms, we can write 
\[\dim_{\p G} H_k(G; \D_{\p G})= \dim_{\p G} R_k-n_{k-1}+\dim_{\p G}\im \phi,\]
and, similarly, 
\[\dim_{\p} H_k (G_j; \p)=\dim_{\p}(\p \tens{\p G_j}R_k)-n_{k-1}|G: G_j|+\dim_{\p}(\p \tens{\p G_j}\im f).\]
Putting these pieces together with \cref{Luckmodp0} directly leads to the desired conclusion. 
\end{proof}

\begin{cor} \label{Bettibound}
Let $G$ be a finitely generated dense group of a finitely generated free pro-$p$ group $\FF$. Then we have that $\bp(G)+1\geq d(\FF)$. Moreover, if $G$ is parafree then $\bp(G)+1\geq d(G_{\hp})$.
\end{cor}
\begin{proof}
Let $\FF=\NN_1>\NN_2>\dots$ be a chain of open normal subgroups of $\FF$ with trivial intersection. Define $G_j=G\bigcap \NN_j$. It is clear that the closure of $G_j$ in $\FF$ is $\NN_j$ and that $|G:G_j|=|\FF: \NN_j|$.  As a consequence, the natural map 
\begin{equation} \label{bet0} G_j/G_j^p[G_j, G_j]\lrar \NN_j/\NN_j^p[\NN_j, \NN_j]\end{equation}
has dense image. Since the $p$-abelianisation of $\NN_j$ is finite, the space  $ \NN_j/\NN_j^p[\NN_j, \NN_j] $ is discrete and we  deduce that the previous map is surjective.

We know view $\FF$ as the pro-$p$ completion of a free group $F$ of rank $d(\FF)$ and we consider the embedding $F\inc \FF$. We define $F_j=F\bigcap \NN_j$. Arguing as before, $\NN_j$ are the closure of $F_j$ and $|F:F_j|=|\FF: \NN_j|$. Since $F_j$ has $p$-power index in $F$, we can apply \cref{pro-p subspaces} to deduce that the canonical map $(F_j)_{\hp}\lrar \NN_j$ is an isomorphism of pro-$p$ groups. Hence, by Schreier's index-rank formula (\cref{freeindex}),  $\NN_j$ is a free pro-$p$ group of rank \[d(\NN_j)=d(F_j)= (d(F)-1)|F : F_j|+1= (d(\FF)-1)|\FF : \NN_j|+1.\] Then we obtain that
\begin{align*}
\dim_{\F_p} H_1(G_j; \F_p) =& \log_p |G_j: G_j^p[G_j, G_j]|\geq \log_p |\NN_j, \NN_j^p[\NN_j, \NN_j]| = d(\NN_j)\\
=& (d(\FF)-1)|\FF : \NN_j|+1=(d(\FF)-1)|G: G_j|+1.
\end{align*}
Thus, by \cref{Luckmodp}, we obtain that 
\[\beta_1^{mod\, p}(G)\geq d(\FF)-1.\]
Lastly, if $G$ is parafree, then $G_{\hp}$ is free by \cref{characterisations} and  we study the canonical embedding $G\inc G_{\hp}$. In this case, by \cref{pro-p subspaces}, the canonical map $(G_j)_{\hp}\lrar \NN_j$ is an isomorphism. By \cref{proabiso}, the surjective map of (\ref{bet0}) is actually an isomorphism. Re-doing the previous computations would lead to the equation 
\[\beta_1^{mod\, p}(G)= d(\FF)-1.\qedhere\]
\end{proof}

We conclude this chapter by writing down an explicit relation between several quantities that contain algebraic and topological information of parafree groups. 

\begin{cor} \label{parabetti}
Let $\Ga$ be a finitely generated parafree group. Then
\[\dim_{\D_{\p \Ga}} \p I_{\Ga}=\bp(\Ga)+1=\b(\Ga)+1=\abr (G)=d(\Ga_{\hp}).\]
\end{cor}
\begin{proof}
Since $\Ga_{\ab}$ is torsion-free, we know by \cref{Fra1}, \cref{Fra2} and \cref{proabiso};  that 
\[r_{\ab}(\Ga)=\dim_{\p} \Ga/[\Ga, \Ga]\Ga^p= \dim_{\p} \Ga_{\hp}/\Phi(\Ga_{\hp})= d(\Ga_{\hp}).\]
This is the last equality of the statement. The first equation is a consequence of \cref{Luckmodp}. Furthermore, by \cref{Bettibound}, 
$\bp(\Ga)=d(\Ga_{\hp})-1$. 
In addition,  $\bp (\Ga)\geq \b (\Ga)$ by \cref{Luckmodp} and \cite[Theorem 1.6]{Ers14}.
Lastly,  $\b(\Ga)\geq d(\ti{\Ga}_{\hp})-1$ by \cref{Bettiest1},. The conclusion follows. 
\end{proof}

 \chapter{Embbedings into free pro-{$p$} groups} \label{pro-$p$ embeddings}

In this section we detail a method for ensuring that a map from an abstract group $\ti{G}$ to a free pro-$p$ group is injective. We are particularly interested in the problem of producing families of parafree groups. Let  $\ti{G}$ be a candidate to being parafree, meaning that $\ti{G}_{\hat{p}}$ is free for every prime $p$. We want to study whether the canonical map $\ti{G}\lrar \ti{G}_{\hat{p}}$ is an embedding for some suitable prime $p$. This would establish the residual nilpotence of $\ti{G}$ and we could conclude that $\ti{G}$ is parafree.  

 An important fact that underlies the statement is that subgroups $G$ of free pro-$p$ groups $\FF$ are residually torsion-free nilpotent and hence $\F_p G$ admits a universal division ring of fractions $\D_{\F_p G}$.

\begin{thm}[Tool for constructing pro-$p$ embeddings] \label{pro-$p$ embedding lemma} Let $\ti{G}$ be a finitely generated abstract group, let $\FF$ be a finitely generated free pro-$p$ group and let $\phi : \ti{G}\rar  \FF$ be a group homomorphism. Suppose that we have the following conditions. 
\begin{enumerate}
    \item The image $G=\phi(\ti{G})$ is dense in $\FF$.
    \item The $\F_p G$-module $\p G\tens{\F_p \ti{G}} \p I_{\ti{G}}$ is a $\D_{\F_p G}$-torsion-free module of dimension  $d(\FF)$.
    \item The kernel of $\phi$ is free. 
\end{enumerate}
 Then the map $\phi$ is an embedding. 
\end{thm}
\begin{proof}
We will prove that the surjective map $\phi: \ti{G}\rar G$ verifies the assumption of \cref{kernelmodp} to deduce that $\ker \phi$ has trivial $p$-abelianisation. Since $\ker \phi$ is free, the latter would imply that $\ker \phi=1$ and the conclusion would follow. 

It is clear that the natural  homomorphism of $\F_p G$-modules \[\p G\tens{\F_p \ti{G}} \p I_{\ti{G}}\lrar \F_p I_{G},\]
defined by $a \otimes b\mapsto \phi (a)b$, is  surjective. If it was not injective, then, naming its kernel by $L$ and naming $M=\p I_{\ti{G}}\tens{\F_p \ti{G}} \F_p G$, we would deduce, after applying \cref{quotientdim}, that \[d(\FF)= \dim_{\D_{\F_p G}} M>\dim_{\D_{\F_p G}}(M/L)=\dim_{\D_{\F_p G}}(\F_p I_{G}).\]
This would contradict \cref{Bettibound}. 
\end{proof}

In the setting of our problem, that is, taking some $\ti{G}$ with free pro-$p$ completion and studying whether $\ti{G}\lrar \ti{G}_{\hat{p}}$
is injective, we make a few comments about the three conditions of \cref{pro-$p$ embedding lemma}. We take $G$ to be the image of $\ti{G}$ inside  $\ti{G}_{\hat{p}} $.

\begin{enumerate}
    \item The first condition will be naturally ensured. 
    \item The second condition will require to handle an augmentation ideal with tools from \cref{augsection}. Establishing torsion-freeness is the hardest part and requires  the most technical arguments, mostly from sections \ref{completedsection} and \ref{universalsection}.
    \item The third condition is natural from the point of view of Bass-Serre theory. If we take some free construction $\ti{G}$, meaning that $\ti{G}$ is the fundamental group of a graph of groups, and we ensure $\ker \phi$ intersects trivially every vertex group, then the $\ker \phi$ will necessarily be free. 
\end{enumerate}

We will use this theorem to construct families of parafree groups in subsequent sections by iterations of amalgamated products  and HNN extensions. These sources of examples will include  the examples of parafree groups $G_{i, j}, H_{i, j}, K_{i, j}, N_{p, q, r}$ described in \cref{parafamilies}.

\section{Amalgamated products}\label{paraamal}


The following statement describes exactly under which circumstances an amalgamated product of finitely generated  groups with abelian amalgams is parafree. 

\begin{thm}[{\bf Parafree amalgamated products of cyclic amalgam}] \label{amalgamparafree} Let $\ti{H_1}$ and $\ti{H_2}$ be finitely generated groups. Let $1\neq u_1\in \ti{H_1}$ and let $1\neq u_2\in \ti{H_2}$. Consider the following amalgamated product of cyclic amalgam
\[\ti{G}=\ti{H_1}\underset{u_1=u_2}{*}\ti{H_2}\cong \frac{\ti{H_1}* \ti{H_2}}{\lan \lan u_1u_2^{-1}\ran \ran}.\]
Then $\ti{G}$ is parafree if and only if the three following conditions hold:
\begin{enumerate}
    \item The groups $\ti{H}_1$ and $\ti{H}_2$ are parafree.
    \item The element $u_1u_2^{-1}$ of $\ti{H_1}* \ti{H_2}$ is not a proper power in the abelianisation.
    \item There is at least one $i\in \{1, 2\}$ such that $u_i$ is not a proper power in $\ti{H_i}$.
\end{enumerate} 
In this case, $\ti{G}$ is parafree of abelian rank $\abr(\ti{H_1})+\abr(\ti{H_2})-1$. 
\end{thm}
\begin{proof} We first observe that both conditions are necessary. 
\begin{enumerate}
    \item Let us suppose that $\ti{G}$ is parafree. Since both $\ti{H_1}$ and $\ti{H_2}$ are subgroups of $\ti{G}$, then they are residually nilpotent. In order to prove that they are, in fact, parafree, it would suffice to show, by \cref{characterisations}, that $\ti{H_i}_{\hp}$ is pro-$p$ free for all primes $p$. 
    By construction, it is clear that 
    \[d\left(\ti{G}/\ti{G}^p\, [\ti{G}, \ti{G}]\right)\geq d\left(\ti{H_1}/\ti{H_1}^p\, [\ti{H_1}, \ti{H_1}]\right)+d\left(\ti{H_2}/\ti{H_2}^p\, [\ti{H_2}, \ti{H_2}]\right)-1.\]
    So, by the equation (\ref{proabiso}), 
    \begin{equation}\label{amalnec1}
        d(\ti{G}_{\hp})\geq  d(\ti{H_1}_{\hp})+d(\ti{H_2}_{\hp})-1.
    \end{equation}
    By \cref{characterisations}, $\ti{G}_{\hp}$ is a free pro-$p$ group. Consider the canonical map $\phi: \ti{G}\lrar \ti{G}_{\hat{p}}$, which is injective, and name $H_1=\phi(\ti{H_1})$ and $H_2=\phi(\ti{H_2})$. By \cref{subfree}, $\ove{H_1}$ and $\ove{H_2}$ are free pro-$p$ groups. There is a canonical continuous homomorphism $\ti{H_i}_{\hp}\lrar \ove{H_i}$, which is clearly surjective. From this, it follows that
     \begin{equation}\label{amalnec3}
    d(\ti{H_i}_{\hp})\geq d(\ove{H_i}), \, \, \, \mbox{for $i\in \{1, 2\}$}.
    \end{equation}
    If we managed to prove that 
    \begin{equation}\label{amalnec2}
    d(\ti{H_i}_{\hp})=d(\ove{H_i}),  \, \, \, \mbox{for $i\in \{1, 2\}$},
    \end{equation}
    then we would conclude, by \cref{Hopf2}, that $\ti{H_i}_{\hp}$ is free, as we want. 
    By the universal property (\ref{unicop}) of the coproduct, there is a continuous homomorphism
\[f: \overline{H_1}\coprod \overline{H_2}\lrar \ti{G}_{\hat{p}}\]
which sends $\overline{H_i}$ to each corresponding copy in $\ti{G}_{\hat{p}}$. Notice that $\im f$ contains both $H_1$ and $H_2$, so $\phi(\ti{G})\leq \im f$, implying that $f$ is surjective. In addition, $f$ has non-trivial kernel since $\phi(u_j)\neq 1$ in $\overline{H_j}$ and $f(\phi (u_1)\phi(u_2)^{-1})=1$ in $\ti{G}_{\hat{p}}$. Here we have used that the canonical map $\iota: \overline{H_1}* \overline{H_2}\lrar \overline{H_1}\coprod \overline{H_2} $ is injective (see \cref{injectiveabs}). 
    
This verifies that $f$ is a non-injective and surjective map onto a free pro-$p$ group.  Hence, by \cref{Hopf2} and \cref{additivecop},
\[d(\overline{H_1})+d(\overline{H_2})=d\left( \overline{H_1}\coprod \overline{H_2}  \right)>d(\ti{G}_{\hat{p}}).\]
The previous inequality together with (\ref{amalnec1}) and (\ref{amalnec3}) directly imply (\ref{amalnec2}). 
    
    \item If the second condition did not hold then the abelianisation of $\ti{G}$ would not be torsion-free.
    \item Lastly, suppose that $u_1=v_1^{n_1}$ in $\ti{H_1}$, and that $u_2=v_2^{n_2}$ in $\ti{H_2}$; with $\min \{|n_1|, |n_2|\}\geq 2$.   Since $u_i\neq 1$ and $\ti{H_i}$ is torsion-free, then   $v_i\notin \lan u_i\ran$. By \cref{brittonam1}, 
    \[[v_1, v_2]={v_1}^{-1}v_2^{-1}v_1v_2\neq 1.\]
    Then the two-generated subgroup $\lan tv_1t^{-1}, v_2\ran$ of $\ti{G}$ would not be free because $[tv_1t^{-1}, v_2]\neq 1$ and $(tv_1t^{-1})^{n_1}=v_2^{n_2}$. By
    \cref{twogen}, we deduce $\ti{G}$ is not parafree.
\end{enumerate}

We now verify that the three given conditions are sufficient.

There is a canonical isomorphism 
\[\ti{G}_{\ab}\cong \frac{\ti{H_1}_{\ab}\oplus \ti{H_2}_{\ab}}{u_1-u_2}.\]
From the fact that $u_1-u_2$ is not a proper power in $\ti{G}_{\ab}$, we see that $\ti{G}_{\ab}$ is torsion-free of rank $\abr(\ti{H_1})+\abr(\ti{H_2})-1$. By (\ref{proabiso}), we rephrase this as
\[d(\ti{G}_{\hp})=d(\ti{H_1}_{\hp})+d(\ti{H_2}_{\hp})-1.\]

Let us fix an arbitrary prime $p$ from this point on. Since $u_1u_2^{-1}$ is not a proper power in $\ti{G}_{\ab}$, at least one of $u_1$ or $u_2$, say $u_j$;  is primitive in the $p$-abelianisation $\ti{G}/\ti{G}^p[\ti{G}, \ti{G}]$. 

Consider the canonical map $\phi: \ti{G}\lrar \ti{G}_{\hat{p}}$.

\begin{lem} \label{amal1} The pro-$p$ group $\ti{G}_{\hat{p}}$ is free and the element $\phi(u_j)$ is primitive in $\ti{G}_{\hat{p}}$.
\end{lem}
To prove this lemma, we name $\Gamma=\ti{H_1}*\ti{H_2}$ and we consider the canonical commutative diagram
\begin{equation}
    \begin{tikzcd}
    \Ga \ar[r] \ar[d] & \Ga_{\hat{p}}\ar[d] \ar[r, green, "\cong"] & \ti{H_1}_{\hat{p}}\coprod \ti{H_2}_{\hat{p}} \ar[d] \\
    \Ga/\Ga^p[\Ga, \Ga] \ar[r, blue, "\cong"] & \Ga_{\hat{p}}\big/ \Phi\left( \Ga_{\hat{p}}\right) \ar[r, red, "\cong"] &  {\ti{H_1}}_{\hat{p}}\big/ \Phi\left( {\ti{H_1}}_{\hat{p}}\right)\oplus  {\ti{H_2}}_{\hat{p}}\big/ \Phi\left( {\ti{H_2}}_{\hat{p}}\right),
    \end{tikzcd}
\end{equation}
with isomorphisms marked with $\cong$. The isomorphism in blue is proven in (\ref{proabiso}); the green one, in \cref{propprod}; and the red one, in \cref{Fra3}. 

By the choice of $u_j$,   $\{u_1u_2^{-1}, u_j\}$ are  $\F_p$-independent in $\Ga/\Ga^p[\Ga, \Ga]. $ This implies that the pair $\{\phi(u_1)\phi(u_2)^{-1}, \phi(u_j)\}$ is primitive in $\Ga_{\hat{p}}\big/ \Phi\left( \Ga_{\hat{p}}\right)$. So the latter pair is part of a set of free generators of $\Gamma_{\hat{p}}$ by \cref{freebasis}.

From \cref{propamal}, we know there is a canonical isomorphism of pro-$p$ groups 
\[\ti{G}_{\hat{p}}\cong \ti{H_1}_{\hat{p}}\coprod \ti{H_2}_{\hat{p}} \Big/ \, \overline{\lan \lan u_1 \, u_2^{-1}\ran \ran}. \]
So $\ti{G}_{\hat{p}}$ is a free pro-$p$ group and $\phi(u_j)$ is primitive in $\ti{G}_{\hat{p}}$. The \cref{amal1} is proven. 

\begin{lem} \label{amal2} The restrictions of $\phi$ to each $\ti{H_i}$ are injective. 
\end{lem}

To check this lemma, we consider each restriction \[\phi_i=\phi|_{\ti{H_i}}: \ti{H_i}\lrar H_i,\] where $H_i=\phi(\ti{H_i})$. The  subgroups  $\overline{H_1}$ and $\overline{H_2}$ of the free pro-$p$ group $\ti{G}_{\hat{p}}$  are closed. Hence they both are free pro-$p$ groups by \cref{subfree}. Since the induced ${\phi_{i}}_{\hat{p}}: \ti{H_i}_{\hat{p}}\lrar \overline{H_i}$ are surjective maps of free pro-$p$ groups, then, by \cref{Hopf2}, 
\begin{equation} \label{amal2eq1} d(\ti{H_i}_{\hat{p}})\geq d(\overline{H_i}), \, \, \, \mbox{for all $i\in \{1, 2\}.$}\end{equation}
Furthermore, by the universal property (\ref{unicop}) of the coproduct, there is a continuous homomorphism
\[f: \overline{H_1}\coprod \overline{H_2}\lrar \ti{G}_{\hat{p}}\]
which sends $\overline{H_i}$ to each corresponding copy in $\ti{G}_{\hat{p}}$. Notice that $\im f$ contains both $H_1$ and $H_2$, so $G\leq \im f$, implying that $f$ is surjective. In addition, $f$ has non-trivial kernel since $\phi(u_j)\neq 1$ in $\overline{H_j}$ and $f(\phi (u_1)\phi(u_2)^{-1})=1$ in $\ti{G}_{\hat{p}}$. Here we have used that the canonical map $\iota: \overline{H_1}* \overline{H_2}\lrar \overline{H_1}\coprod \overline{H_2} $ is injective (see \cref{injectiveabs}). 

This verifies that $f$ is a non-injective and surjective map onto a free pro-$p$ group. Hence, by \cref{Hopf2} and \cref{additivecop},
\[d(\overline{H_1})+d(\overline{H_2})=d\left( \overline{H_1}\coprod \overline{H_2}  \right)>d(\ti{G}_{\hat{p}}) =d(\ti{H_1}_{\hat{p}})+d(\ti{H_2}_{\hat{p}})-1.\]
The previous inequality together with (\ref{amal2eq1}) imply that 
\[d(\ti{H_i}_{\hat{p}})=d(\overline{H_i})\, \, \, \mbox{for all $i\in \{1, 2\}.$}\]
We now look at the commutative diagram of canonical arrows
\begin{equation*}
    \begin{tikzcd}
    \ti{H_i} \ar[r, "\phi_i"] \ar[d, "\iota_{\hp}"] & H_i \ar[d, hookrightarrow]\\
    \ti{H_i}_{\hat{p}}\ar[r, "{\phi_i}_{\hat{p}}"]& \overline{H_i}.\\
    \end{tikzcd}
\end{equation*}
The map ${\phi_i}_{\hat{p}}$ is a continuous surjection between finitely generated free pro-$p$ groups of the same rank, so it is an isomorphism by \cref{Hopf2}. Moreover, each $\ti{H_i}$ is parafree, so $\iota_{\hp}$ is injective, too. From this, we deduce that the surjective maps $\phi_i$ are also injective and \cref{amal2} is proven.

\begin{lemma} \label{amal3} The map $\phi: \ti{G}\lrar \ti{G}_{\hat{p}}$ is injective. 
\end{lemma}

We denote $G=\phi(\ti{G})$. We want to apply \cref{pro-$p$ embedding lemma} to the map $\phi: \ti{G}\lrar G$. On the one hand, notice that 
the kernel of $\phi$ intersects trivially the amalgam $\lan u\ran$ of $\ti{G}$, so, by \cref{amalfree}, the kernel is free. 

We consider the corresponding 
$\p G$-module 
\[M=\p G\tens{\p \ti{G}} \p I_{\ti{G}}.\]

Consider $A=\lan \phi(u_1)\ran$ and  \[J=\{(x, -x): x\in \p I_A^{G}\}\leq \p I_{H_1}^G\oplus \p I_{H_2}^G.\]
Then, by \cref{augamal}, 

\[M\cong \frac{\p I_{H_1}^G\oplus \p I_{H_2}^G}{J}.\] 

Without loss of generality, we suppose that $u_1$ is not a proper power in $\ti{H_1}$. 
Recall that $\phi_i: \ti{H_i}\lrar H_i$ is an isomorphism, so $H_i$ is parafree and  $\phi(u_1)$ is not a proper power in $\ti{H_1}$. By \cref{twogen}, $A=\lan u_1\ran $ is a maximal abelian subgroup of $H_1$.

By \cref{amaltor}, $M$ is $\D_{\p G}$-torsion-free with dimension 
\[ \dim_{\p G} M =  \dim_{\p H_1} \p I_{H_1}+\dim_{\p H_2}\p I_{H_2}-1.\]
In addition, by \cref{parabetti}, 
 \[  \dim_{\p H_1} \p I_{H_1}+\dim_{\p H_2}\p I_{H_2}-1=d({H_1}_{\hat{p}})+d({H_2}_{\hat{p}})-1 = d(\ti{G}_{\hat{p}}).\]
It follows that $\dim_{\p G} M=d(\ti{G}_{\hp})$.

We have verified the conditions of \cref{pro-$p$ embedding lemma} for the homomorphism $\phi: \ti{G}\lrar G\sub \ti{G}_{\hp}$. So \cref{amal3} is proven and $\ti{G}$ is residually nilpotent. We already know that each $\ti{G}_{\hat{p}}$ is free so, by \cref{characterisations}, $\ti{G}$ is parafree. 
\end{proof}

\section{HNN extensions}\label{parahnn}


The following result does not entirely describe which cyclic HNN extensions  are parafree, as in the case of amalgamated products, though it significantly reduces the complexity of the conditions to be verified. 

\begin{thm}[{\bf Parafree HNN cyclic extensions}] \label{hnnparafree} Let $\ti{H}$ be a finitely generated group. Let $u, v\in \ti{H}\setminus \{1\}$. Consider the following cyclic HNN extension of $\ti{H}$ \[\ti{G}=\ti{H} \hn{\lan u\ran}=\frac{\ti{H}* \lan t\ran}{\lan \lan tut^{-1}v^{-1}\ran \ran}. \]
Then $\ti{G}$ is  parafree if and only if the four following conditions hold:
\begin{enumerate}
    \item The group $\ti{H}$ is parafree.
    \item The image of $uv^{-1}$ is not a proper power in $\ti{H}_{\ab}$.
    \item At least one of $u$ or $v$ is not a proper power in $\ti{H}$.
    \item The image of $u$ is non-trivial in some finite nilpotent quotient of $\ti{G}$.
\end{enumerate}
In this case,  $\ti{G}$ is  parafree of the same abelian rank as $\ti{H}$. 
\end{thm}

\begin{proof} The given conditions are necessary. 
\begin{enumerate}
    \item Let us suppose that $\ti{G}$ is parafree. Since  $\ti{H}$ is a subgroup of $\ti{G}$, then it is residually nilpotent. In order to prove that it is parafree, it would suffice to show, by \cref{characterisations}, that $\ti{H}_{\hp}$ is pro-$p$ free for all primes $p$. 
    By construction, it is clear that 
    \[d\left(\ti{G}/\ti{G}^p\, [\ti{G}, \ti{G}]\right)\geq d\left(\ti{H}/\ti{H}^p\, [\ti{H}, \ti{H}]\right)\]
    So, by the equation (\ref{proabiso}), 
    \begin{equation}\label{hnnnec1}
        d(\ti{G}_{\hp})\geq  d(\ti{H}_{\hp}).
    \end{equation}
    By \cref{characterisations}, $\ti{G}_{\hp}$ is a free pro-$p$ group. Consider the canonical map $\phi: \ti{G}\lrar \ti{G}_{\hat{p}}$, which is injective, and name $H=\phi(\ti{H})$. By \cref{subfree}, $\ove{H}$ is a pro-$p$ group. There is a canonical continuous homomorphism $\ti{H}_{\hp}\lrar \ove{H}$, which is clearly surjective. From this, it follows that 
     \begin{equation}\label{hnnnec3}
    d(\ti{H}_{\hp})\geq d(\ove{H}).
    \end{equation}
    If we managed to prove that 
    \begin{equation}\label{hnnnec2}
    d(\ti{H}_{\hp})=d(\ove{H}),
    \end{equation}
    then we would conclude, by \cref{Hopf2}, that $\ti{H}_{\hp}$ is free, as we want. 
    By the universal property (\ref{unicop}) of the coproduct, there is a continuous homomorphism
\[f: \overline{H}\coprod \Z_p\lrar \ti{G}_{\hat{p}}\]
which sends $\overline{H}$ to its corresponding copy in $\ti{G}_{\hat{p}}$ and $\Z_p$ to the cyclic pro-$p$ group generated by $\phi(t)$. Notice that $\im f$ contains both $H$ and $\lan \phi(t)\ran$, so $\phi(\ti{G})\leq \im f$, implying that $f$ is surjective. In addition, $f$ has non-trivial kernel since $\phi(u_1)\neq 1$ in $\overline{H}$ and $f(\phi(t) (u_1)\phi(t)^{-1}\phi(u_2)^{-1})=1$ in $\ti{G}_{\hat{p}}$. Here we have used that the canonical map $\iota: \overline{H}* \Z_p\lrar \overline{H}\coprod \Z_p $ is injective (see \cref{injectiveabs}). 
    
This verifies that $f$ is a non-injective and surjective map onto a free pro-$p$ group.  Hence, by \cref{Hopf2} and \cref{additivecop},
\[d(\overline{H})+1=d\left( \overline{H}\coprod \Z_p  \right)>d(\ti{G}_{\hat{p}}).\]
The previous inequality together with (\ref{hnnnec1}) and (\ref{hnnnec3}) directly imply (\ref{hnnnec2}). 
    \item If the first condition did not hold then the abelianisation of $\ti{G}$ would not be torsion-free.
    \item Secondly, suppose that $u=w^{n_1}$ in $\ti{H_1}$, and that $v=w_2^{n_2}$ in $\ti{H_2}$; with $\min \{|n_1|, |n_2|\}\geq 2$.   Since $u, v\neq 1$   and $\ti{H_i}$ is torsion-free, then   $w_1\notin \lan u\ran$ and $w_2\notin \lan v\ran$. By \cref{brittonhnn1}, 
    \[[tw_1t^{-1}, w_2]=t{w_1}^{-1}t^{-1}w_2^{-1}tw_1t^{-1}w_2\neq 1.\]
    Then the two-generated subgroup $\lan tw_1t^{-1}, w_2\ran$ of $\ti{G}$ would not be free because $[tw_1t^{-1}, w_2]\neq 1$ and $(tw_1t^{-1})^{n_1}=w_2^{n_2}$. By
    \cref{twogen}, we would deduce that $\ti{G}$ is not parafree.
    \item If $\ti{G}$ were parafree then it would be residually-$p$;  hence $\ti{G}\inc \ti{G}_{\hp}$ would be an embedding and $u$ would not have trivial image in $\ti{G}_{\hp}$. 
\end{enumerate}

We are going to verify that these conditions are sufficient for $\ti{G}$ to be parafree. 

First of all, it is clear that 
\[\ti{G}_{\ab}\cong \frac{\ti{H}_{\ab}}{\lan u-v\ran}\oplus \lan t\ran.\]
From the fact that $u-v$ is not a proper power in $\ti{H}_{\ab}$, we see that $\ti{G}$ is torsion-free of the same rank as $\ti{H}$. 
By (\ref{proabiso}), we rephrase this as 
\[d(\ti{G}_{\hp})=d(\ti{H}_{\hp}).\]
In addition, we also see that $t$ is primitive in $\ti{G}_{\ab}$. Let us fix an arbitrary prime $p$. 

Consider the canonical map $\phi: \ti{G}\lrar \ti{G}_{\hat{p}}$.

\begin{lemma} \label{hnn0} The pro-$p$ group $\ti{G}_{\hp}$ is free and the element $\phi(t)$ is primitive in $\ti{G}_{\hp}$.
\end{lemma}

To prove this lemma, we name $\Gamma=\ti{H}*\lan t\ran$ and we consider the canonical commutative diagram
\begin{equation}
    \begin{tikzcd}
    \Ga \ar[r] \ar[d] & \Ga_{\hat{p}}\ar[d] \ar[r, green, "\cong"] & \ti{H}_{\hat{p}}\coprod \Z_p \ar[d] \\
    \Ga/\Ga^p[\Ga, \Ga] \ar[r, blue, "\cong"] & \Ga_{\hat{p}}\big/ \Phi\left( \Ga_{\hat{p}}\right) \ar[r, red, "\cong"] &  {\ti{H}}_{\hat{p}}\big/ \Phi\left( {\ti{H}}_{\hat{p}}\right)\oplus  \Z/p,
    \end{tikzcd}
\end{equation}
with isomorphisms marked with $\cong$. The isomorphism in blue is proven in (\ref{proabiso}); the green one, in \cref{propprod}; and the red one, in \cref{Fra3}. 
The pair $\{\phi(tut^{-1}v^{-1}), \phi(t)\}$ is are $\F_p$-linearly independent in $\Ga/\Ga^p[\Ga, \Ga]$. So this pair is also primitive in $\Ga_{\hat{p}}\big/ \Phi(\Ga_{\hat{p}})$. By \cref{freebasis}, they are part of a set of free generators of $\Ga_{\hp}$.

  By \cref{proppres}, we see that there is a canonical isomorphism 
\[\ti{G}_{\hat{p}}\cong \ti{H}_{\hat{p}} \coprod \Z_p\, \Big/ \, \overline{\lan \lan t u t^{-1} v^{-1}\ran \ran},\]
where the generator of $\Z_p$ has been denoted by $t$. 
So $\ti{G}_{\hp}$ is free; $\phi(t)$ is part of a topological generating set; and \cref{hnn0} is proven.

We name $H=\phi(\ti{H})$. 

\begin{lem}\label{hnn1} The restriction of $\phi$ to $\ti{H}$ is injective. 
\end{lem}

To verify this, consider the closed subgroup $\overline{H}\leq \ti{G}_{\hat{p}}$. Since $\ti{G}_{\hat{p}}$ is free, the pro-$p$ group $\overline{H}$ must be free by \cref{subfree}. We notice that the epimorphism $\phi: \ti{H}\lrar H$ induces a continuous epimorphism $\phi_{\hat{p}}: \ti{H}_{\hat{p}}\lrar \overline{H}$. 
In particular, by \cref{Hopf2}, that 
\begin{equation} \label{n>d} d(\ti{H}_{\hp})\geq d(\ove{H}).\end{equation} 
Furthermore, by the universal property (\ref{unicop}) of the coproduct, there is a continuous homomorphism
\[f: \overline{H}\coprod \Z_p \lrar \ti{G}_{\hat{p}},\]
which sends $\overline{H}$ to $\ti{G}_{\hat{p}}$, by inclusion; and   $\Z_p$ to the cyclic pro-$p$ group generated by $\phi(t)$. Since the image of $f$ contains both $H$ and $\phi(t)$, it follows that $\im f$ contains $\phi(\ti{G})$, so $f$ must be a surjective. In addition, it has a nontrivial kernel; since $\phi(t)\phi(u)\phi(t)^{-1}\phi(v)^{-1}=1$ and $\phi(u)\neq 1$, by assumption. Here we have used   that the canonical map $\iota: \ove{H}* \Z_p\lrar  \overline{H}\coprod \Z_p$ is injective (see \cref{injectiveabs}).

This verifies that $f$ is a non-injective and surjective continuous homomorphism  onto a free pro-$p$ group. Hence, by \cref{Hopf2} and \cref{additivecop},
\[d(\overline{H})+1=d(\overline{H})+d(\Z_p)=d\left( \overline{H}\coprod \Z_p \right)> d(\ti{G}_{\hat{p}}). \]
This, in addition to (\ref{n>d}), implies that $d(\ti{H}_{\hp})=d(\overline{H})$. So  $\phi_{\hat{p}}: \ti{H}_{\hat{p}}\rar \overline{H}$ is a continuous epimorphism between free pro-$p$ groups of the same rank. By \cref{Hopf2}, this restriction $\phi_{\hat{p}}|_{\ti{H}}$ must be an isomorphism. 

Looking at the commutative diagram
\begin{equation*}
    \begin{tikzcd}
    \ti{H} \ar[r, "\phi"] \ar[d, "\iota"] & H \ar[d, hookrightarrow]\\
    \ti{H}_{\hat{p}}\ar[r, "\phi_{\hat{p}}"]& \overline{H},\\
    \end{tikzcd}
\end{equation*}
we observe that $\phi: \ti{H}\lrar H$ is injective. This is because
$\phi_{\hat{p}}|_{\ti{H}}$ is injective, as we just proved; and $\iota_{\hp}$ is injective, due to $\ti{H}$ being parafree. 

 The verification of \cref{hnn1} is complete. 
 
\begin{lemma}\label{hnn3} The map $\phi: \ti{G}\lrar \ti{G}_{\hp}$ is injective.
\end{lemma}

We denote $G=\phi(\ti{G})$. We want to apply \cref{pro-$p$ embedding lemma} to the map $\phi: \ti{G}\lrar G\sub \ti{G}_{\hp}$. We already know, from \cref{hnn0}, that  $\ti{G}_{\hp}$ is free. We already On the one hand, notice that 
the kernel of $\phi$ intersects trivially the subgroup $\lan u\ran$ of $\ti{G}$. By \cref{hnntree}, the kernel is free.

By \cref{freebasis}, we can pick a set $\{a_2, \dots, a_n\}$ of $n-1$ elements of $\ti{H}$ such that  \[\{\phi(t), \phi(a_2), \dots, \phi(a_n)\}\] is a set of free topological generators of $\ti{G}_{\hat{p}}$. We define a continuous homomorphism  $q: \ti{G}_{\hat{p}}\lrar \Z_p$ such that $q(\phi(t))=1$; and $q(\phi(a_k))=0$ if $2\leq k\leq n$. The restriction $q|_G$ verifies that its its kernel $\ker q|_G$ contains $H$. This is because $\{\phi(a_2), \dots, \phi(a_n)\}$ belong to the kernel of $q$; and they generate $H$ modulo $[H, H]$. In addition, $t\notin \ker q|_{G}$, by construction. Since $G/\ker q|_G\cong \Z^m$ for some $m\geq 1$, we can take $N\unlhd G$ such that $\ker q|_G\leq N$, $G/N\cong \Z$ and $G=\lan N, t\ran\cong N\rtimes \lan t\ran$.

We now consider the $\p G$-module 
\[M=\p G\tens{\p \ti{G}} \p I_{\ti{G}}.\]
By \cref{aughnn1}, this $\p G$-module is isomorphic to 
\[M\cong \frac{\p I_{H}^G\oplus \p G}{\p G\, (\phi(v)-1-\phi(t)(\phi(u)-1),\,  \phi(v)-1)}.\]

Without loss of generality, we suppose that $u$ is not a proper power in $\ti{H}$. 
Since $\phi: \ti{H}\lrar H$ is an isomorphism, $H$ is parafree and $\phi (u)$ is not a proper power in $H$. From \cref{twogen}, we deduce that that  $\lan \phi(u)\ran$ is a maximal abelian subgroup of $H$. By \cref{hnntor}, this implies that $M$ is $\D_{\p G}$-torsion-free.  

Using \cref{freequo}, we have the following isomorphisms of $\D_{\p G}$-modules
\[\D_{\p G}\tens{\p G} M\cong \D_{\p G}\tens{\p G} \p I_H^G\cong \D_{\p G}\tens{\p H} \p I_H\cong \D_{\p G}\tens{\D_{\p H}} \left( \D_{\p H}\tens{\p H} \p I_H \right).\]
Combining these isomorphisms with \cref{parabetti}, yields to
\[\dim_{\p G} M= \dim_{\p H} \p I_H=d(\ti{H}_{\hat{p}})=d(\ti{G}_{\hat{p}}).\]

We have verified the conditions of \cref{pro-$p$ embedding lemma} for the homomorphism $\phi: \ti{G}\lrar G\sub \ti{G}_{\hp}$. So \cref{hnn3} is proven and $\ti{G}$ is residually nilpotent. We already know that very pro-$p$ completion $\ti{G}_{\hat{p}}$ is free; so, by \cref{characterisations}, $\ti{G}$ is parafree. 
\end{proof}

In the case of amalgamated products (\cref{amalgamparafree}), we could explicitly describe when the resulting group would be parafree. However, in the case of the HNN extensions treated in \cref{hnnparafree}, we gave a characterisation whose last condition (the image of $u$ being non-trivial in some finite nilpotent quotient) can be hard to check. Our last result will get rid of this condition in the case when the starting group $\ti{H}$ has abelian rank equal to 2.

\begin{cor}[{\bf Parafree cyclic HNN extensions of groups with abelian rank 2}] \label{hnnparafree2}
Let $\ti{H}$ be a finitely generated group of abelian rank 2. Let $u, v\in \ti{H}\setminus \{1\}$. Consider the following cyclic HNN extension of $\ti{H}$ \[\ti{G}=\ti{H} \hn{\lan u\ran}=\frac{\ti{H}* \lan t\ran}{\lan \lan tut^{-1}v^{-1}\ran \ran}. \]
Then $\ti{G}$ is  parafree if and only if the four following conditions hold:
\begin{enumerate}
    \item The group $\ti{H}$ is parafree.
    \item The image of $uv^{-1}$ is not a proper power in $\ti{H}_{\ab}$.
    \item At least one of $u$ or $v$ is not a proper power in $\ti{H}$.
    \item  The images of $u$ and $v$ generate a subgroup isomorphic to $\Z^2$ in $\ti{H}_{\ab}$. 
\end{enumerate}
In this case,  $\ti{G}$ is  parafree of abelian rank 2. 
\end{cor}
\begin{proof}
It is not hard to see that the four given conditions  imply the conditions of \cref{hnnparafree}. In fact, under these assumptions,  the image of $u$ is primitive in all the pro-$p$ completions of $\ti{G}$. 

In order to prove that they are necessary,  let us suppose that  $\ti{G}$ is parafree. By \cref{hnnparafree}, $(1), (2)$ and $(3)$ hold. We assume that $(4)$ does not hold. Then, by $(2)$, the image of $uv^{-1}$ in $\ti{H}$ generates a subgroup which contains both $u$ and $v$. So there exists $c\in \Z$ such that $u\equiv (uv^{-1})^c \mod [\ti{H}, \ti{H}]$. In addition, since $\ti{H}_{\hp}$ is free of rank 2 and the image of $uv^{-1}$ is primitive in $\ti{H}_{\hp}$; we can denote by $\NN$ the closed normal subgroup generated by $uv^{-1}$ in $\ti{H}_{\hp}$ and observe that $\ti{H}_{\hp}/\NN\cong \Z_p$ because $uv^{-1}$ has primitive image in $\ti{H}_{\ab}$. Thus $[\ti{H}, \ti{H}]\sub \NN$ and $u\in (uv^{-1})^c$. 

Denote by $\NN_1$ the closed normal subgroup of $\ti{G}_{\hp}$ generated by $uv^{-1}$. Then $u\in \NN_1$ by the previous argument. Now observe that $uv^{-1}=[u^{-1}, t^{-1}]$. We are going to verify that $u\in \bigcap_k \ga_k \ti{G}_{\hp}=\{1\}.$ This would result in the contradiction $u=1$, finishing the proof. We proceed inductively to check the claim. The base of the induction is trivial because $u\in \ga_1\ti{G}_{\hp}=\ti{G}_{\hp}$. If $u\in \ga_k\ti{G}_{\hp}$, then $uv^{-1}\in [\ga_k \ti{G}_{\hp}, \ti{G}_{\hp}]=\ga_{k+1} \ti{G}_{\hp}$. Since $\ga_k \ti{G}_{\hp}$ is a closed normal subgroup of $\ti{G}_{\hp}$, then $\NN_1\sub \ga_{k+1} \ti{G}_{\hp}$ and hence $u\in \ga_{k+1} \ti{G}_{\hp}$, completing the inductive step. 
\end{proof}

\section{General fundamental groups}
There is a way to extend our results characterising parafree amalgamated products and HNN extensions for more general fundamental groups.

\begin{cor}[Corollary 1.4 of \cite{And21}] \label{fundamentalparafreeintr}  Let $(\mathcal G, \Gamma)$ be a  graph of   groups over a finite 
graph $\Gamma$ whose edge morphisms are injective. Let  $W=\pi(\mathcal G, \Gamma)$ be its  fundamental group.  Assume that all  vertex subgroups  $\mathcal G(v)$ ($v\in V(\Gamma)$)   are finitely generated and all   edge subgroups $\mathcal G(e)$ ($e\in E(\Gamma)$)  are cyclic. Then $W$ is parafree if and only if the following four conditions hold.
 
 \begin{enumerate}
 \item All the vertex  subgroups $\mathcal G(v)$ ($v\in V(\Gamma)$)    are parafree.
 
 \item The abelianisation of $W$ is torsion-free   of rank 
 $$\abr(W)=\sum_{v\in V(\Gamma)}\abr(\mathcal G(v))-\sum_{e\in E(\Gamma)} \abr(\mathcal G(e)) -\chi(\Gamma),$$   where $\chi(\Gamma)=|V(\Gamma)|-|E(\Gamma)|$-1.

 \item All the centralisers of non-trivial elements in $W$ are cyclic.

 \item For  each non-trivial edge subgroup of  $\mathcal G(e)$  ($e\in E(\Gamma)$) there is a finite nilpotent quotient of $W$ where the image of this edge subgroup is non-trivial.
 \end{enumerate}

 \end{cor}

\backmatter

\bibliography{aomsample}

\providecommand{\bysame}{\leavevmode\hbox to3em{\hrulefill}\thinspace}
\providecommand{\noopsort}[1]{}
\providecommand{\mr}[1]{\href{http://www.ams.org/mathscinet-getitem?mr=#1}{MR~#1}}
\providecommand{\zbl}[1]{\href{http://www.zentralblatt-math.org/zmath/en/search/?q=an:#1}{Zbl~#1}}
\providecommand{\jfm}[1]{\href{http://www.emis.de/cgi-bin/JFM-item?#1}{JFM~#1}}
\providecommand{\arxiv}[1]{\href{http://www.arxiv.org/abs/#1}{arXiv~#1}}
\providecommand{\doi}[1]{\url{https://doi.org/#1}}
\providecommand{\MR}{\relax\ifhmode\unskip\space\fi MR }
\providecommand{\MRhref}[2]{%
  \href{http://www.ams.org/mathscinet-getitem?mr=#1}{#2}
}
\providecommand{\href}[2]{#2}
\begin{thebibliography}{10}

\bibitem{Bar16}
\bgroup\scshape{}L.~Bartholdi\egroup{}, Amenability of groups is characterized
  by {M}yhill's theorem,  \emph{Journal of the European Mathematical Society}
  \textbf{21} no.~10 (2019), 3191–3197, With an appendix by D. Kielak.

\bibitem{Bau62}
\bgroup\scshape{}G.~Baumslag\egroup{}, On generalised free products,
  \emph{Mathematische Zeitschrlft} \textbf{78} (1962), 423--–438.

\bibitem{Bau67}
\bgroup\scshape{}G.~Baumslag\egroup{}, Groups with the same lower central
  sequence as a relatively free group. i. the groups,  \emph{Transactions of
  the American Mathematical Society} \textbf{129} no.~2 (1967), 308--321.

\bibitem{Bau69}
\bgroup\scshape{}G.~Baumslag\egroup{}, Groups with the same lower central
  sequence as a relatively free group. ii. the properties,  \emph{Transactions
  of the American Mathematical Society} \textbf{142} (1969), 507--538.

\bibitem{Bau05}
\bgroup\scshape{}G.~Baumslag\egroup{}, \emph{Parafree Groups}, pp.~1--14,
  Birkh{\"a}user Basel, Basel, 2005.

\bibitem{Bau06}
\bgroup\scshape{}G.~Baumslag\egroup{} and \bgroup\scshape{}S.~Cleary\egroup{},
  Parafree one-relator groups,  \emph{Journal of Group Theory} \textbf{9} no.~2
  (2006), 191--201.

\bibitem{Bri14}
\bgroup\scshape{}M.~R. Bridson\egroup{} and \bgroup\scshape{}A.~W.
  Reid\egroup{}, Nilpotent completions of groups, {G}rothendieck pairs, and
  four problems of {B}aumslag,  \emph{International Mathematics Research
  Notices} \textbf{2015} no.~8 (2014), 1073--7928.

\bibitem{Bro82}
\bgroup\scshape{}K.~S. Brown\egroup{}, \emph{Cohomology of Groups},
  \emph{Graduate Texts in Mathematics}, Springer-Verlag New York, 1982.

\bibitem{Cle17}
\bgroup\scshape{}A.~E. Clement\egroup{}, \bgroup\scshape{}S.~Majewicz\egroup{},
  and \bgroup\scshape{}M.~Zyman\egroup{}, \emph{The theory of nilpotent
  groups}, Birkhäuser Basel, 2017.

\bibitem{Coh74}
\bgroup\scshape{}P.~M. Cohn\egroup{}, Localization in semifirs,  \emph{Bulletin
  of the London Mathematical Society} \textbf{6} (1974), 13--20.

\bibitem{Coh95}
\bgroup\scshape{}P.~M. Cohn\egroup{}, \emph{Skew fields. Theory of general
  division rings}, \emph{Enciclopedia of Mathematics and its applications}
  \textbf{57}, Cambridge University Press, 1995.

\bibitem{Coh06}
\bgroup\scshape{}P.~M. Cohn\egroup{}, \emph{Free Ideal Rings and Localization
  in General Rings}, \emph{New Mathematical Monographs}, Cambridge University
  Press, 2006.

\bibitem{Dic80}
\bgroup\scshape{}W.~Dicks\egroup{}, \emph{Groups, Trees and Projective
  Modules}, \emph{Lecture Notes in Mathematics}, Springer-Verlag Berlin
  Heidelberg, 1980.

\bibitem{Dic07}
\bgroup\scshape{}W.~Dicks\egroup{} and \bgroup\scshape{}P.~Linnell\egroup{},
  {$L^2$}-betti numbers of one-relator groups,  \emph{Annals of Mathematics}
  \textbf{337} (2007), 855–874.

\bibitem{Dix99}
\bgroup\scshape{}J.~D. Dixon\egroup{}, \bgroup\scshape{}M.~P.~F.
  Du~Sautoy\egroup{}, \bgroup\scshape{}A.~Mann\egroup{}, and
  \bgroup\scshape{}D.~Segal\egroup{}, \emph{Analytic Pro-P Groups}, 2 ed.,
  \emph{Cambridge Studies in Advanced Mathematics}, Cambridge University Press,
  1999.

\bibitem{Dix82}
\bgroup\scshape{}J.~D. Dixon\egroup{}, \bgroup\scshape{}E.~W.
  Formanek\egroup{}, \bgroup\scshape{}J.~C. Poland\egroup{}, and
  \bgroup\scshape{}L.~Ribes\egroup{}, Profinite completions and isomorphic
  finite quotients,  \emph{Journal of Pure and Applied Algebra} \textbf{23}
  no.~3 (1982), 227--231.

\bibitem{Ers14}
\bgroup\scshape{}M.~Ershov\egroup{} and \bgroup\scshape{}W.~Lück\egroup{}, The
  first {$L^2$}-{B}etti number and approximation in arbitrary characteristic,
  \emph{Doc. Math} \textbf{19} (2014), 313--332.

\bibitem{Gab02}
\bgroup\scshape{}D.~Gaboriau\egroup{}, Invariants {$\ell^2$} de relations
  d’équivalence et de groupes,  \emph{Institut des Hautes \'Etudes
  Scientifiques} \textbf{95} (2002), 93--150.

\bibitem{Gru57}
\bgroup\scshape{}K.~W. Gruenberg\egroup{}, Residual properties of infinite
  soluble groups,  \emph{Proceedings of the London Mathematical Society}
  \textbf{(3) 7} (1957), 29--62.

\bibitem{Hal76}
\bgroup\scshape{}M.~Hall\egroup{}, \emph{The theory of groups}, Chelsea
  Publishing Company, New York, 1976, Second edition.

\bibitem{Hem72}
\bgroup\scshape{}J.~Hempel\egroup{}, Residual finiteness of surface groups,
  \emph{Proceedings of the American Mathematical Society} \textbf{32} no.~1
  (1972), 323--323.

\bibitem{Her15}
\bgroup\scshape{}D.~Herbera\egroup{} and
  \bgroup\scshape{}J.~S{\'a}nchez\egroup{}, The inversion height of the free
  field is infinite,  \emph{Selecta Mathematica, New Series} \textbf{21} no.~3
  (2015), 883--929.

\bibitem{Hug70}
\bgroup\scshape{}I.~Hughes\egroup{}, Division rings of fractions for group
  rings,  \emph{Communications on Pure and Applied Mathematics} \textbf{23}
  (1970), 181--188.

\bibitem{And191}
\bgroup\scshape{}A.~Jaikin-Zapirain\egroup{}, An explicit construction of the
  universal division ring of fractions of {$E\lan \lan x_1, \dots, x_n\ran
  \ran$},  \emph{Journal of Combinatorial Algebra} (2020), 369–--395.

\bibitem{And212}
\bgroup\scshape{}A.~Jaikin-Zapirain\egroup{}, The finite and soluble genus of a
  finitely generated free group,  \emph{Preprint} (2021). Available at
  \url{https://matematicas.uam.es/~andrei.jaikin/preprints/profinitefree.pdf}.

\bibitem{And20}
\bgroup\scshape{}A.~Jaikin-Zapirain\egroup{}, Free {$\Q$}-groups are residually
  torsion-free nilpotent,  \emph{Preprint} (2021). Available at \url{http:
  //matematicas.uam.es/~andrei.jaikin/preprints/baumslag.pdf}.

\bibitem{And19}
\bgroup\scshape{}A.~Jaikin-Zapirain\egroup{}, The universality of {H}ughes-free
  division rings,  \emph{Selecta Mathematica} \textbf{27} no.~74 (2021).

\bibitem{And192}
\bgroup\scshape{}A.~Jaikin-Zapirain\egroup{} and
  \bgroup\scshape{}D.~López-Álvarez\egroup{}, The {S}trong {A}tiyah and
  {L}ück {A}pproximation conjectures for one-relator groups.,
  \emph{Mathematische Annalen} \textbf{376} (2019), 1741–1793.

\bibitem{And21}
\bgroup\scshape{}A.~Jaikin-Zapirain\egroup{} and
  \bgroup\scshape{}I.~Morales\egroup{}, Parafree fundamental groups of graph of
  groups with cyclic edge subgroups,  \emph{In preparation} (2021).

\bibitem{Joh97}
\bgroup\scshape{}D.~L. Johnson\egroup{}, \emph{Presentations of Groups}, 2 ed.,
  \emph{London Mathematical Society Student Texts}, Cambridge University Press,
  1997.

\bibitem{Kir95}
\bgroup\scshape{}R.~Kirby\egroup{}, Problems in low dimensional {T}opology,
  (1995).

\bibitem{Li20}
\bgroup\scshape{}H.~Li\egroup{}, Bivariant and extended sylvester rank
  functions,  \emph{Journal of the London Mathematical Society} (2021),
  222--249.

\bibitem{Lin93}
\bgroup\scshape{}P.~A. Linnell\egroup{}, Division rings and group von neumann
  algebras.,  \emph{Forum mathematicum} \textbf{5} no.~6 (1993), 561--576.

\bibitem{Lot95}
\bgroup\scshape{}J.~Lott\egroup{} and \bgroup\scshape{}W.~Lück\egroup{},
  {$L^2$}-topological invariants of compact 3-manifolds,  \emph{Inventiones
  Mathematicae} \textbf{120} (1995), 15--–60.

\bibitem{Lyn62}
\bgroup\scshape{}R.~C. Lyndon\egroup{} and \bgroup\scshape{}M.~P.
  Schützenberger\egroup{}, The equation $a^m=b^nc^p$ in a free group.,
  \emph{Michigan Mathematical Journal} \textbf{9} no.~4 (1962), 289 -- 298.

\bibitem{Lyn01}
\bgroup\scshape{}R.~C. Lyndon\egroup{} and \bgroup\scshape{}P.~E.
  Schupp\egroup{}, \emph{Combinatorial Group Theory}, \emph{Classics in
  Mathematics}, Springer, Berlin, Heidelberg, 2001.

\bibitem{Luc94}
\bgroup\scshape{}W.~Lück\egroup{}, Approximating {$L^2$}-invariants by their
  finite-dimensional analogues,  \emph{Geometric and Functional Analysis}
  (1994), 455--–481.

\bibitem{Luc11}
\bgroup\scshape{}W.~Lück\egroup{} and \bgroup\scshape{}D.~Osin\egroup{},
  Approximating the first {$L^2$}-betti number of residually finite groups,
  \emph{Journal of Topology and Analysis} \textbf{3} no.~2 (2011), 153--160.

\bibitem{May99}
\bgroup\scshape{}J.~P. May\egroup{}, \emph{A Concise Course in Algebraic
  Topology}, \emph{Chicago Lectures in Mathematics}, University of Chicago
  Press, 1999.

\bibitem{Neu49}
\bgroup\scshape{}B.~H. Neumann\egroup{}, On ordered division rings,
  \emph{Transactions of the American Mathematical Society} no.~1 (1949),
  202--252.

\bibitem{Nik06}
\bgroup\scshape{}N.~Nikolov\egroup{} and \bgroup\scshape{}D.~Segal\egroup{}, On
  finitely generated profinite groups, i: strong completeness and uniform
  bounds,  \emph{Annals of Mathematics} \textbf{165} (2006), 171--238.

\bibitem{Ore31}
\bgroup\scshape{}O.~Ore\egroup{}, Linear equations in non-commutative fields,
  \emph{Annals of Mathematics} \textbf{32} no.~3 (1931), 463--477.

\bibitem{Rib00}
\bgroup\scshape{}L.~Ribes\egroup{} and \bgroup\scshape{}P.~Zalesskii\egroup{},
  \emph{Profinite Groups}, Springer Berlin Heidelberg, 2000.

\bibitem{Sta65}
\bgroup\scshape{}J.~Stallings\egroup{}, Homology and central series of groups,
  \emph{Journal of Algebra} \textbf{2} no.~2 (1965), 170--181.

\bibitem{Tak51}
\bgroup\scshape{}M.~Stallings\egroup{}, Note on chain conditions in free
  groups,  \emph{Osaka Mathematical journal} \textbf{3} no.~2 (1951), 221--225.

\bibitem{San08}
\bgroup\scshape{}J.~Sánchez\egroup{}, Localization: On division rings and
  tilting modules,  \emph{Ph. D. Thesis, Universidad Autónoma de Barcelona}
  (2008). Available at
  \url{www.tdx.cat/bitstream/handle/10803/3107/jss1de1.pdf}.

\bibitem{Wil21}
\bgroup\scshape{}G.~Wilkes\egroup{}, Lecture notes on ``{Profinite groups and
  Group Cohomology}'', 2020-2021. Available at
  \url{https://www.dpmms.cam.ac.uk/~grw46/partiiiprofinite.html}.

\bibitem{Wil98}
\bgroup\scshape{}J.~S. Wilson\egroup{}, \emph{Profinite Groups}, \emph{London
  Mathematical Society Monographs}, Clarendon Press, 1998.

\end{thebibliography}
\bibliographystyle{aomplain}

\end{document}